\documentclass[12pt]{amsart} 
\usepackage{amssymb,latexsym,tikz}
\usepackage{graphics, tikz}
\usepackage[top=40mm, bottom=25mm, left=35mm, right=35mm]{geometry}

\usetikzlibrary{arrows,decorations.markings}
\tikzstyle arrowstyle=[scale=1]
\tikzstyle directed=[postaction={decorate,
decoration={markings,mark=at position .65 with {\arrow[arrowstyle]{stealth}}}}]

\usepackage{enumitem}
\usepackage{accents}
\usepackage{mathtools}
\usepackage{caption}

\usetikzlibrary{decorations.markings}
\usetikzlibrary{math, backgrounds,quotes, positioning, babel, 3d, graphs, graphs.standard}
\tikzset{mid vert/.style={/utils/exec=\tikzset{every node/.append style={outer sep=0.8ex}},
postaction=decorate,decoration={markings,
mark=at position 0.5 with {\draw[-] (0,#1) -- (0,-#1);}}},
mid vert/.default=0.75ex}

\usepackage{amscd}
\usepackage{xypic}

\usepackage[utf8]{inputenc}

\usepackage{graphicx}

\usepackage{subfig}

\usepackage{nicematrix}

\begin{document}

\title 
[Polarizations of Artin monomial ideals]
{Polarizations of Artin monomial ideals}

\author{Gunnar Fl{\o}ystad}
\address{Universitetet i Bergen \\
  Matematisk institutt\\       
         Postboks 7803\\
         5020 Bergen \\
         Norway}
       \email{nmagf@uib.no}

\author{Ine Gabrielsen}
\address{Universitetet i Bergen \\
  Matematisk institutt\\       
         Postboks 7803\\
         5020 Bergen \\
         Norway}
       \email{ine.gabrielsen@student.uib.no}
       
\author{Amir Mafi}
\address{Department of Mathematics \\ University Of Kurdistan \\ P.O. Box: 416,
  Sanandaj \\ Iran}
\email{a.mafi@uok.ac.ir}

\keywords{triangulation, ball, sphere, simplicial complex, cell complex,
  polarization,
  Artin monomial ideal, linkage, squeezed ball, Bier ball, letterplace
ideal}
\subjclass[2020]{Primary: 13F55; Secondary: 05E40, 05E45}
\date{\today}

\theoremstyle{plain}
\newtheorem{theo}{Theorem}[section]
\newtheorem{theorem}[theo]{Theorem}
\newtheorem{coro}[theo]{Corollary}
\newtheorem{lemm}[theo]{Lemma}
\newtheorem{lemma}[theo]{Lemma}
\newtheorem{prop}[theo]{Proposition}
\newtheorem{conj}[theo]{Conjecture}

\theoremstyle{definition}
\newtheorem{defi}[theo]{Definition}
\newtheorem{fact}[theo]{Fact}
\newtheorem{problem}{Problem}

\theoremstyle{remark}
\newtheorem{nota}[theo]{Notation}
\newtheorem{rema}[theo]{Remark}
\newtheorem{exam}[theo]{Example}

\newcommand{\psp}[1]{{{\bf P}^{#1}}}
\newcommand{\psr}[1]{{\bf P}(#1)}

\newcommand{\ini}[1]{\text{in}(#1)}
\newcommand{\gin}[1]{\text{gin}(#1)}
\newcommand{\kr}{{\Bbbk}}
\newcommand{\pd}{\partial}
\newcommand{\vardel}{\partial}
\renewcommand{\tt}{{\bf t}}

\newcommand{\coh}{{{\text{{\rm coh}}}}}

\newcommand{\modv}[1]{{#1}\text{-{mod}}}
\newcommand{\modstab}[1]{{#1}-\underline{\text{mod}}}

\newcommand{\sut}{{}^{\tau}}
\newcommand{\sumit}{{}^{-\tau}}
\newcommand{\til}{\thicksim}

\newcommand{\totp}{\text{Tot}^{\prod}}
\newcommand{\dsum}{\bigoplus}
\newcommand{\dprod}{\prod}
\newcommand{\lsum}{\oplus}
\newcommand{\lprod}{\Pi}

\newcommand{\La}{{\Lambda}}

\newcommand{\sirstj}{\circledast}

\newcommand{\she}{\EuScript{S}\text{h}}
\newcommand{\cm}{\EuScript{CM}}
\newcommand{\cmd}{\EuScript{CM}^\dagger}
\newcommand{\cmri}{\EuScript{CM}^\circ}
\newcommand{\cler}{\EuScript{CL}}
\newcommand{\clerd}{\EuScript{CL}^\dagger}
\newcommand{\clerri}{\EuScript{CL}^\circ}
\newcommand{\gor}{\EuScript{G}}
\newcommand{\cF}{\mathcal{F}}
\newcommand{\cG}{\mathcal{G}}

\newcommand{\cM}{\mathcal{M}}
\newcommand{\cE}{\mathcal{E}}

\newcommand{\cI}{\mathcal{I}}
\newcommand{\cP}{\mathcal{P}}
\newcommand{\cK}{\mathcal{K}}

\newcommand{\cS}{\mathcal{S}}
\newcommand{\cC}{\mathcal{C}}
\newcommand{\cO}{\mathcal{O}}
\newcommand{\cJ}{\mathcal{J}}
\newcommand{\cU}{\mathcal{U}}

\newcommand{\cQ}{\mathcal{Q}}
\newcommand{\cX}{\mathcal{X}}
\newcommand{\cY}{\mathcal{Y}}
\newcommand{\cZ}{\mathcal{Z}}
\newcommand{\cV}{\mathcal{V}}

\newcommand{\mm}{\mathfrak{m}}

\newcommand{\dlim} {\varinjlim}
\newcommand{\ilim} {\varprojlim}

\newcommand{\CM}{\text{CM}}
\newcommand{\Mon}{\text{Mon}}

\newcommand{\Kom}{\text{Kom}}

\newcommand{\EH}{{\mathbf H}}
\newcommand{\res}{\text{res}}
\newcommand{\Hom}{{\rm{Hom}}}
\newcommand{\inhom}{{\underline{\text{Hom}}}}
\newcommand{\Ext}{\text{Ext}}
\newcommand{\Tor}{\text{Tor}}
\newcommand{\ghom}{\mathcal{H}om}
\newcommand{\gext}{\mathcal{E}xt}
\newcommand{\id}{\text{{id}}}
\newcommand{\im}{\text{im}\,}
\newcommand{\codim} {\text{codim}\,}
\newcommand{\resol}{\text{resol}\,}
\newcommand{\rank}{\text{rank}\,}
\newcommand{\lpd}{\text{lpd}\,}
\newcommand{\coker}{\text{coker}\,}
\newcommand{\supp}{\text{supp}\,}
\newcommand{\Ad}{A_\cdot}
\newcommand{\Bd}{B_\cdot}
\newcommand{\Fd}{F_\cdot}
\newcommand{\Gd}{G_\cdot}

\newcommand{\sus}{\subseteq}
\newcommand{\sups}{\supseteq}
\newcommand{\pil}{\rightarrow}
\newcommand{\vpil}{\leftarrow}
\newcommand{\rpil}{\leftarrow}
\newcommand{\lpil}{\longrightarrow}
\newcommand{\inpil}{\hookrightarrow}
\newcommand{\pils}{\twoheadrightarrow}
\newcommand{\projpil}{\dashrightarrow}
\newcommand{\dotpil}{\dashrightarrow}
\newcommand{\adj}[2]{\overset{#1}{\underset{#2}{\rightleftarrows}}}
\newcommand{\mto}[1]{\stackrel{#1}\longrightarrow}
\newcommand{\vmto}[1]{\stackrel{#1}\longleftarrow}
\newcommand{\mtoelm}[1]{\stackrel{#1}\mapsto}
\newcommand{\bihom}[2]{\overset{#1}{\underset{#2}{\rightleftarrows}}}
\newcommand{\eqv}{\Leftrightarrow}
\newcommand{\impl}{\Rightarrow}

\newcommand{\iso}{\cong}
\newcommand{\te}{\otimes}
\newcommand{\into}[1]{\hookrightarrow{#1}}
\newcommand{\ekv}{\Leftrightarrow}
\newcommand{\equi}{\simeq}
\newcommand{\isopil}{\overset{\cong}{\lpil}}
\newcommand{\equipil}{\overset{\equi}{\lpil}}
\newcommand{\ispil}{\isopil}
\newcommand{\vvi}{\langle}
\newcommand{\hvi}{\rangle}
\newcommand{\susneq}{\subsetneq}
\newcommand{\sgn}{\text{sign}}

\newcommand{\xd}{\check{x}}
\newcommand{\ortog}{\bot}
\newcommand{\tL}{\tilde{L}}
\newcommand{\tM}{\tilde{M}}
\newcommand{\tH}{\tilde{H}}
\newcommand{\tvH}{\widetilde{H}}
\newcommand{\tvh}{\widetilde{h}}
\newcommand{\tV}{\tilde{V}}
\newcommand{\tS}{\tilde{S}}
\newcommand{\tT}{\tilde{T}}
\newcommand{\tR}{\tilde{R}}
\newcommand{\tf}{\tilde{f}}
\newcommand{\ts}{\tilde{s}}
\newcommand{\tp}{\tilde{p}}
\newcommand{\tr}{\tilde{r}}
\newcommand{\tfst}{\tilde{f}_*}
\newcommand{\empt}{\emptyset}
\newcommand{\bfa}{{\mathbf a}}
\newcommand{\bfb}{{\mathbf b}}
\newcommand{\bfd}{{\mathbf d}}
\newcommand{\bfl}{{\mathbf \ell}}
\newcommand{\bfx}{{\mathbf x}}
\newcommand{\bfm}{{\mathbf m}}
\newcommand{\bfv}{{\mathbf v}}
\newcommand{\bft}{{\mathbf t}}
\newcommand{\bbfa}{{\mathbf a}^\prime}
\newcommand{\la}{\lambda}
\newcommand{\bfen}{{\mathbf 1}}
\newcommand{\bfe}{{\mathbf 1}}
\newcommand{\ep}{\epsilon}
\newcommand{\en}{r}
\newcommand{\tu}{s}
\newcommand{\Sym}{\text{Sym}}

\newcommand{\ome}{\omega_E}

\newcommand{\bevis}{{\bf Proof. }}
\newcommand{\demofin}{\qed \vskip 3.5mm}
\newcommand{\nyp}[1]{\noindent {\bf (#1)}}
\newcommand{\demo}{{\it Proof. }}
\newcommand{\demodone}{\demofin}
\newcommand{\parg}{{\vskip 2mm \addtocounter{theorem}{1}  
                   \noindent {\bf \thetheorem .} \hskip 1.5mm }}

\newcommand{\lcm}{{\text{lcm}}}

\newcommand{\dl}{\Delta}
\newcommand{\cdel}{{C\Delta}}
\newcommand{\cdelp}{{C\Delta^{\prime}}}
\newcommand{\dlst}{\Delta^*}
\newcommand{\Sdl}{{\mathcal S}_{\dl}}
\newcommand{\lk}{\text{lk}}
\newcommand{\lkd}{\lk_\Delta}
\newcommand{\lkp}[2]{\lk_{#1} {#2}}
\newcommand{\del}{\Delta}
\newcommand{\delr}{\Delta_{-R}}
\newcommand{\dd}{{\dim \del}}
\newcommand{\Del}{\Delta}

\newcommand{\bb}{{\bf b}}
\newcommand{\cc}{{\bf c}}
\newcommand{\xx}{{\bf x}}
\newcommand{\yy}{{\bf y}}
\newcommand{\zz}{{\bf z}}
\newcommand{\mv}{{\xx^{\aa_v}}}
\newcommand{\mF}{{\xx^{\aa_F}}}

\newcommand{\Symm}{\text{Sym}}
\newcommand{\pnm}{{\bf P}^{n-1}}
\newcommand{\opnm}{{\go_{\pnm}}}
\newcommand{\ompnm}{\omega_{\pnm}}

\newcommand{\pn}{{\bf P}^n}
\newcommand{\hele}{{\mathbb Z}}
\newcommand{\nat}{{\mathbb N}}
\newcommand{\rasj}{{\mathbb Q}}
\newcommand{\bfone}{{\mathbf 1}}

\newcommand{\dt}{\bullet}

\newcommand{\disk}{\scriptscriptstyle{\bullet}}

\newcommand{\cxF}{F_\dt}
\newcommand{\pol}{f}

\newcommand{\Rn}{{\mathbb R}^n}
\newcommand{\An}{{\mathbb A}^n}
\newcommand{\frg}{\mathfrak{g}}
\newcommand{\PW}{{\mathbb P}(W)}

\newcommand{\pos}{{\mathcal Pos}}
\newcommand{\g}{{\gamma}}

\newcommand{\Vaa}{V_0}
\newcommand{\Bp}{B^\prime}
\newcommand{\Bpp}{B^{\prime \prime}}
\newcommand{\bbp}{\mathbf{b}^\prime}
\newcommand{\bbpp}{\mathbf{b}^{\prime \prime}}
\newcommand{\bp}{{b}^\prime}
\newcommand{\bpp}{{b}^{\prime \prime}}

\newcommand{\oLa}{\overline{\Lambda}}
\newcommand{\ov}[1]{\overline{#1}}
\newcommand{\ovv}[1]{\overline{\overline{#1}}}
\newcommand{\tm}{\tilde{m}}
\newcommand{\po}{\bullet}

\newcommand{\surj}[1]{\overset{#1}{\twoheadrightarrow}}
\newcommand{\Supp}{\text{Supp}}

\newcommand{\ZZ}{{\mathbb Z}}
\newcommand{\NN}{{\mathbb N}}
\newcommand{\RR}{{\mathbb R}}

\newcommand{\oR}{\overline{R}}

\newcommand{\bfu}{{\mathbf u}}
\newcommand{\nn}{{\mathbf n}}
\newcommand{\oa}{\overline{a}}
\newcommand{\cop}{\text{cop}}

\newcommand{\opos}{\text{op}\,}

\newcommand{\ngmi}{\text{neg}}
\newcommand{\up}{\text{up}}
\newcommand{\dw}{\text{down}}

\newcommand{\diw}[1]{\widehat{#1}}
\newcommand{\di}[1]{\diw{#1}}
\newcommand{\bo}{b}
\newcommand{\ub}{u}
\newcommand{\fs}{\infty}
\newcommand{\ifst}{\infty}
\newcommand{\mon}{{mon}}
\newcommand{\cl}{\text{cl}}
\newcommand{\intr}{\text{int}}
\newcommand{\ul}[1]{\underline{#1}}

\newcommand{\bipil}{\leftrightarrow}
\newcommand{\bfc}{{\mathbf c}}

\renewcommand{\mp}{m^\prime}

\newcommand{\np}{n^\prime}
\newcommand{\Mod}{\text{Mod }}
\newcommand{\Sh}{\text{Sh } }
\newcommand{\st}{\text{st}}

\newcommand{\hM}{\tilde{M}}
\newcommand{\hs}{\tilde{s}}
\newcommand{\ee}{\mathbf{e}}

\newcommand{\lex}{{\text{lex}}}
\newcommand{\ordGL}{\succeq_{\lex}}
\newcommand{\ordG}{\succ_{\lex}}
\newcommand{\ordML}{\preceq_{\lex}}
\newcommand{\ordM}{\prec_{\lex}}
\newcommand{\tLa}{\tilde{\Lambda}}
\newcommand{\tGa}{\tilde{\Gamma}}
\newcommand{\STS}{\text{STS}}
\newcommand{\ii}{{\rotatebox[origin=c]{180}{\scalebox{0.7}{\rm{!}}}}}

\newcommand{\jj}{\mathbf{j}}

\newcommand{\shmod}{\texttt{shmod}\,}
\newcommand{\hf}{\underline{f}}
\newcommand{\Glim}{\lim}
\newcommand{\Gcolim}{\colim}
\newcommand{\fm}{f^{\underline{m}}}
\newcommand{\fn}{f^{\underline{n}}}
\newcommand{\gn}{g_{\mathbf{|}n}}
\newcommand{\se}[1]{\overline{#1}}
\newcommand{\cB}{{\mathcal{B}}}
\newcommand{\fin}{{\text{fin}}}
\newcommand{\Poset}{{\text{\bf Poset}}}
\newcommand{\Set}{{\text{\bf Set}}}

\newcommand{\Homi}{\text{Hom}^{L}}
\newcommand{\Homb}{\text{Hom}_S}
\newcommand{\Homu}{\text{Hom}^u}
\newcommand{\bfnu}{{\mathbf 0}}

\newcommand{\semip}{up semi-finite }
\newcommand{\semim}{down semi-finite }
\newcommand{\dif}[1]{\di{#1}_{fin}}
\newcommand{\llin}{\raisebox{1pt}{\scalebox{1}[0.6]{$\mid$}}}
\newcommand{\promap}{\mathrlap{{\hskip 2.8mm}{\llin}}{\lpil}}
\newcommand{\dual}{{\widehat{}}}
\newcommand{\Bool}{{\bf Bool }}
\newcommand{\pro}{{pro}}
\newcommand{\ovF}{\overline{F}}
\newcommand{\TO}{{\mathrm to}}
\newcommand{\dT}{\overset{\rightarrow}{T}}

\newlength{\dhatheight}
\newcommand{\doublehat}[1]{%
    \settoheight{\dhatheight}{\ensuremath{\widehat{#1}}}%
    \addtolength{\dhatheight}{-0.35ex}%
    \widehat{\vphantom{\rule{1pt}{\dhatheight}}%
    \smash{\widehat{#1}}}}

\newcommand{\napo}{natural }

\newcommand{\colim}{\text{colim}}
\newcommand{\oPsi}{\overline{\Psi}}
\newcommand{\bfC}{{\mathbf C}}
\newcommand{\bfi}{{\mathbf i}}
\newcommand{\ik}{{\mathfrak k}}
\newcommand{\ic}{{\mathfrak c}}

\newcommand{\ijfr}{{\mathfrak j}}
\newcommand{\iifr}{{\mathfrak i}}

\newcommand{\iC}{{\mathfrak C}}

\newcommand{\Ab}{A_\bullet}
\newcommand{\Rb}{R_\bullet}
\newcommand{\oDo}{\overline{D^\opos}}
\newcommand{\bfA}{{\mathbf A}}
\newcommand{\het}{\text{ht}}

\newcommand{\bfal}{{\mathbf{ \alpha}}}
\newcommand{\al}{{\alpha}}
\newcommand{\alp}{{\alpha^\prime}}
\newcommand{\tci}{tight}
\newcommand{\bitci}{bitight}
\newcommand{\oD}{\overline{D}}
\newcommand{\Delst}{\vardel \Delta}
\newcommand{\mlin}{-}

\newcommand{\ben}{{\mathbf 1}}
\newcommand{\bfj}{{\mathbf j}}
\newcommand{\bfio}{\hat{\bfi}}
\newcommand{\pot}{{\rm{pos}}}
\newcommand{\rac}{r-acyclic}

\newcommand{\Ign}[1]{}

\begin{abstract}
  
  We show that any polarization of an Artin monomial ideal
defines a triangulated ball. This settles a conjecture of
A.Almousa, H.Lohne and the first author \cite{AFL}.

Geometrically, polarizations of ideals
strictly containing $(x_1^{a_1}, \ldots, x_n^{a_n})$ define full-dimensional
triangulated balls on the sphere which is the {\it join}
of boundaries of simplices of
dimensions $a_1-1, \cdots, a_n-1$. We prove that every full-dimensional
Cohen-Macaulay sub-complex of this joined sphere is of this kind, and 
these balls are constructible.

Such a triangulated ball has a dual cell complex which is a
sub-complex of the {\it product} of simplices of dimensions $a_1-1, \cdots , 
a_n-1$. We prove that this cell complex gives a cellular minimal
free resolution of
the Alexander dual ideal of the triangulated ball. When the product of
simplices is a hypercube, using these dual cell complexes
we classify in a range of examples
all polarizations of the Artin monomial ideal.

We also show that the squeezed balls of G.Kalai \cite{Ka} derive
from polarizations of Artin monomial ideals.
\end{abstract}

\maketitle
\setcounter{tocdepth}{1}

\tableofcontents

\maketitle

\section{Introduction}
A polarization of an Artin monomial ideal is a squarefree monomial
ideal and so corresponds to a simplicial complex via the Stanley-Reisner (SR)
correspondence.
We prove a fundamental theorem on the combinatorial geometry of monomial
ideals: 

\medskip
\fbox{\parbox{11.8cm}
{For any {polarization} of an {\it Artin} monomial ideal, the associated
simplicial complex has a topological realization which is a {\it ball}
(save when the Artin monomial ideal is a complete intersection, then
the polarization gives a sphere).}}

\medskip
This proves a conjecture by A.Almousa, H.Lohne, and the first author
\cite[Conjecture 3.2]{AFL}. We give further a detailed geometric
analysis and understanding of polarizations of Artin monomial
ideals in $k[x_1, \ldots, x_n]$.
For a finite set $V$, let $\RR^V$ be the real vector space with basis
the elements of $V$, and let $\Delta(V)$ be the geometric simplex
in $\RR^V$ with vertices $V$. For finite sets $A_1, A_2, \ldots, A_n$
the product polytope
\[ \Pi(A_\bullet) = \Delta(A_1) \times \Delta(A_2) \times \cdots
  \times \Delta(A_n) \]
(if each cardinality $|A_i| = 2$ this is a hypercube), is dual (polar) to
the simplicial polytope which is the join of boundaries
\[ \Delst(\Ab) = \pd \Delta(A_1) * \pd \Delta(A_2) * \cdots * \pd \Delta(A_n),\]
(this is a hyperoctahedron when each $|A_i| = 2$).
It is known that polarizations of Artin monomial ideals are precisely
the full-dimensional Cohen-Macaulay sub-complexes of $\Delst(\Ab)$.
We show that such sub-complexes are balls, by proving that they
are {\it constructible} simplicial complexes.

These sub-complexes have dual cell sub-complexes
of the product polytope $\Pi(\Ab)$, and have the property we
call {\it restriction-acyclic}, meaning that for 
restrictions to every sub-cell $\Pi(\Rb)$ of $\Pi(\Ab)$
(where $R_i \sus A_i$), all
homologies vanish. In many example cases this enables simple analysis
and classification of all polarizations of an Artin monomial ideal,
Section \ref{sec:hypercube}. 

\medskip
We proceed to recall the notions of separation and polarization,
with an example. We further give more detailed descriptions of the above.

\medskip 
\noindent {\bf Separation and polarization.}

For a field $k$ and a set $V$, let $k[V]$ be the polynomial ring with the
elements of $V$ as variables.
\begin{defi} \label{def:intro-sep}
Let $V^\prime \mto{p} V$ be a surjection of finite sets with the
cardinality of $V^\prime$ one more than that of $V$. Let $v_1$ and $v_2$
be the two distinct elements of $V^\prime$ which map to a single
element $v$ in $V$. Let $J$ be a monomial ideal in the polynomial ring
$k[V]$ and $J^\prime$ a monomial ideal in $k[V^\prime]$.
Then $J^\prime$ is a {\it simple separation} of $J$ if:
\begin{itemize}
\item[i.] The monomial ideal $J$ is the image of $J^\prime$ by the map
$k[V^\prime] \pil k[V]$.
\item[ii.] Both the variables $x_{v_1}$ and $x_{v_2}$ occur in some minimal
  generators of $J^\prime$ (usually in distinct generators, but in same
  generator is allowed).
\item[iii.] The variable difference $x_{v_1} - x_{v_2}$ is a non-zero divisor
in the quotient ring $k[V^\prime]/J^\prime$.
\end{itemize}

The minimal generators of $J$ and $J^\prime$ will then be
in bijection.
More generally, if $V^\prime \mto{p} V$ is a surjection of finite sets
and $J \sus k[V]$ and $J^\prime \sus k[V^\prime]$ are monomial ideals such
that $J^\prime$ is obtained by a succession of simple separations of $J$,
$J^\prime$ is a {\it proper separation} of $J$.
 If further $V^\prime \sus V^{\prime \prime}$ is an inclusion,
  the {\it extension} ideal $J^{\prime \prime} = (J^\prime) \sus k[V^{\prime \prime}]$
 is a {\it separation} of the starting ideal $J$.

\end{defi}

Note that $J^{\prime \prime}$ is obtained from $J$ by a sequence of
simple separations where we do not require ii. to hold. 

\begin{defi} \label{def:intro-pol}
  Let $J \sus k[V]$ be a monomial ideal and $V^\prime \pil V$ be a surjection
  of finite sets. An ideal $J^\prime \sus k[V^\prime]$ is a
  {\it polarization}
  of $J$ if $J^\prime$ is squarefree and a separation of $J$.
  It is a {\it proper polarization} if it is a squarefree proper
  separation of $J$. 
\end{defi}

\begin{exam} \label{ex:intro-m2}
  Consider the monomial ideal $\ijfr = (x_1, x_2, x_3)^2
  \sus k[x_1,x_2,x_3]$ (with $k$ a field). 
This is Artin as it contains a power of each
  variable. The ideal has generators:
  \[ \ijfr : x_1^2, x_1x_2, x_2^2, x_1x_3, x_2x_3, x_3^2. \]
  Up to isomorphism this ideal has two polarizations $J_1$ (the standard
  polarization) and $J_2$ (the box polarization)
  with generators:
\begin{align*}  J_1 &:  x_{10}x_{11}, x_{10}x_{20}, x_{20}x_{21},
                     x_{10}x_{30}, x_{20}x_{30}, x_{30}x_{31}, \\
  J_2 &:  x_{10}x_{11}, x_{10}x_{21}, x_{20}x_{21},
        x_{10}x_{31}, x_{20}x_{31}, x_{30}x_{31}.
\end{align*}

These are polarizations because the quotient ring $k[x_1,x_2,x_3]/\ijfr$
is obtained from $k[x_{10}, x_{11}, x_{20},x_{21}, x_{30}, x_{31}]/J_i$
for $i = 1,2$ by dividing out by the regular sequence of
variable differences:
\[ x_{10}-x_{11}, x_{20}-x_{21}, x_{30}-x_{31}.\]
The corresponding simplicial complexes of the two ideals are
given in Figure \ref{fig:HexPol}.
\begin{figure}[ht]
  \begin{center}
\begin{tikzpicture}[dot/.style={draw,fill,circle,inner sep=1pt},scale=1.5]
\draw [help lines,white] (-1.5,-1.5) grid (1.5,1.2);

  \foreach \l [count=\n] in {0,1,2,3,4,5} {
    \pgfmathsetmacro\angle{90-360/6*(\n-1)}
      \node[dot] (n\n) at (\angle:1) {};
    \fill (\angle:1)  circle (0.08);
  }
  \foreach \l [count=\n] in {0,1,2,3,4,5} {
    \pgfmathsetmacro\angle{90-360/6*(\n-1)}
      \node (m\n) at (\angle:1.3) {};
  }
  \draw[thick] (n1) -- (n2) -- (n3) -- (n4) -- (n5) -- (n6) -- (n1);
\draw[fill=gray, fill opacity=0.3] (0,1)--(0.855,0.5)--(0.855,-0.5)--
  (0,-1)--(-0.855,-0.5)--(-0.855,0.5)--(0,1);
  \draw[thick] (n1)--(n3)--(n5)--(n1);

\node[font=\footnotesize] at (m1) {$x_{11}$};
\node[font=\footnotesize] at (m2) {$x_{30}$};
\node[font=\footnotesize] at (m3) {$x_{21}$};
\node[font=\footnotesize] at (m4) {$x_{10}$};
\node[font=\footnotesize] at (m5) {$x_{31}$};
\node[font=\footnotesize] at (m6) {$x_{20}$};
\end{tikzpicture}
  \qquad \quad  
  \begin{tikzpicture}[dot/.style={draw,fill,circle,inner sep=1pt},scale=1.5]
\draw [help lines,white] (-1.5,-1.5) grid (1.5,1.2);

  \foreach \l [count=\n] in {0,1,2,3,4,5} {
    \pgfmathsetmacro\angle{90-360/6*(\n-1)}
      \node[dot] (n\n) at (\angle:1) {};
    \fill (\angle:1)  circle (0.08);
  }
 \foreach \l [count=\n] in {0,1,2,3,4,5} {
    \pgfmathsetmacro\angle{90-360/6*(\n-1)}
      \node (m\n) at (\angle:1.3) {};
  }
  \draw[thick] (n1) -- (n2) -- (n3) -- (n4) -- (n5) -- (n6) -- (n1);
\draw[fill=gray, fill opacity=0.3] (0,1)--(0.855,0.5)--(0.855,-0.5)--
  (0,-1)--(-0.855,-0.5)--(-0.855,0.5)--(0,1);
  \draw[thick] (n6)--(n2)--(n5)--(n3);
 
\node[font=\footnotesize] at (m1) {$x_{10}$};
\node[font=\footnotesize] at (m2) {$x_{30}$};
\node[font=\footnotesize] at (m3) {$x_{21}$};
\node[font=\footnotesize] at (m4) {$x_{31}$};
\node[font=\footnotesize] at (m5) {$x_{11}$};
\node[font=\footnotesize] at (m6) {$x_{20}$};
  
\end{tikzpicture}
\end{center}
\caption{Simplicial complexes of $J_1$ and $J_2$}
\label{fig:HexPol}
\end{figure} 
Topologically both are discs of dimension two, 
i.e. two-dimensional balls.
\end{exam}

\medskip
\noindent {\bf Geometric setting.}
For a set $V$ consider the (abstract)
combinatorial simplex $\Delta(V)$ on $V$. It is
the simplicial complex on $V$ consisting of all subsets of $V$.
Its geometric realization is topologically a ball. Its boundary is the
simplicial complex $\pd \Delta(V)$ whose
maximal faces are all sets $V \backslash \{v\}$ for $v \in V$.
The geometric realization of this boundary is topologically a sphere $S^{|V|-2}$
where $|V|$ is the cardinality of $V$.

Given non-empty disjoint
finite sets $A_1, A_2, \ldots, A_n$ consider the simplicial join
\begin{equation*}
  \Delst(\Ab) = \pd \Delta(A_1) * \pd \Delta(A_2) * \cdots * \pd \Delta(A_n).
 \end{equation*}
 Its topological realization is the topological joins of the
 topological realizations of the $\pd \Delta(A_i)$,
 and so is a sphere of dimension $\sum_{i = 1}^n |A_i| - 1 - n$.
 In particular, if all $A_i$ have cardinality two, $\Delst(\Ab)$
 is the boundary of an $n$-dimensional hyperoctahedron.
 Let $A$ be the (disjoint) union of the $A_i$. We note that the facets
 of $\Delst(\Ab)$ are indexed by elements
 \begin{equation} \label{eq:intro-bfa}
   \bfa = (a_1, a_2, \ldots, a_n)  \in A_1 \times A_2 \times
   \cdots \times A_n.
   \end{equation}
 To such an element corresponds the facet $A \backslash \{a_1, \ldots, a_n\}$. 

 \medskip
  We relate these geometric objects to Artin monomial ideals.
 Let $k[A]$ be the
 polynomial ring with variables from $A$. Denote by $m(A_i)$
 the monomial $\prod_{a \in A_i} a$, and by $\alpha_i$  the cardinality of $A_i$. 
In the polynomial ring $k[x_1, \ldots, x_n]$ we have the complete
 intersection Artin monomial ideal:
 \begin{align*}
   \ic = (x_1^{\alpha_1}, x_2^{\alpha_2}, \ldots, x_n^{\alpha_n}).
 \end{align*}
 It has (up to isomorphism) a unique polarization, the squarefree
 complete intersection ideal
 \[ C = (m(A_1), m(A_2), \ldots, m(A_n)), \]
 which is the Stanley-Reisner ideal of the simplicial sphere $\Delst(\Ab)$.

 \medskip

Our main result is:

\medskip
\noindent {\bf Theorem \ref{thm:rainbow-main}}
 
{\it Let $\ijfr$ be an Artin monomial ideal in $k[x_1, \ldots, x_n]$
  strictly containing the
  complete intersection $\ic$ (but the $x_i^{\alpha_i}$ need not be
  minimal generators of $\ijfr$).
  Let $A \pil \{x_1, x_2, \ldots , x_n\}$  be
  a map with fiber $A_i$ over $x_i$ of cardinality $\alpha_i$.
  
  Any polarization 
  $J \sus k[A]$ of $\ijfr$ defines, via the SR-correspondence,
  a  full dimensional  simplicial ball on the simplicial sphere $\Delst(\Ab)$.
  Moreover, any proper polarization of $\ijfr$ is, when followed by
  a suitable extension,
  isomorphic to such a $J$. }

\medskip
By \cite[Remark 2.5]{AFL} the polarization $J$ in $k[A]$ is proper precisely
when each $x_i^{\alpha_i}$ is
a minimal generator for $\ijfr$. Also, if $\ijfr = \ic$, the polarization
is not a ball, but the sphere $\Delst(\Ab)$.

 \definecolor{cof}{RGB}{219,144,71}
\definecolor{pur}{RGB}{186,146,162}
\definecolor{greeo}{RGB}{91,123,69}
\definecolor{greet}{RGB}{52,111,72}

 \begin{figure}
   \centering

\begin{tikzpicture}[thick,scale=4]

\coordinate(A1) at (0,0);
\coordinate (A2) at (0.6,0.2);
\coordinate (A3) at (1,0);
\coordinate (A4) at (0.4,-0.2);
\coordinate (B1) at (0.5,0.5);
\coordinate (B2) at (0.5,-0.5);

\draw[fill=green,opacity=0.4] (A1)--(B1)--(A4);
\draw[fill=green,opacity=0.4] (A3)--(A4)--(B1);
\draw[fill=green,opacity=0.4] (A3)--(A4)--(B2);
\draw[fill=green,opacity=0.3] (B1)--(A2)--(A3);
\draw[very thick, green](A1)--(B1)--(A3)--(B2);

\begin{scope}[thick,opacity=1]
\draw[dashed](A1)--(B2);
\draw[very thick,green](B2)--(A3);
\draw[dashed](A1)--(A2)--(B2);
\draw[very thick,green](B1)--(A4);
\draw[very thick,green](A3)--(A4);
\draw[very thick,green](A1)--(A4);
\draw[very thick,green](A4)--(B2);
\draw[very thick,green](B1)--(A2);
\draw[very thick,green](A2)--(A3);

\end{scope}

\node[font=\footnotesize, below] at (B2) {$x_{10}$};
\node[font=\footnotesize, above] at (B1) {$x_{11}$};
\node[font=\footnotesize, left] at (A1) {$x_{20}$};
\node[font=\footnotesize, above] at (A2) {$x_{30}$};
\node[font=\footnotesize, right] at (A3) {$x_{21}$};
\node[font=\footnotesize, above] at (A4) {$x_{31}$};

\end{tikzpicture}
\hskip 1.5cm   
\begin{tikzpicture}[scale=4]

\coordinate(A1) at (0,0);
\coordinate (A2) at (0.6,0.2);
\coordinate (A3) at (1,0);
\coordinate (A4) at (0.4,-0.2);
\coordinate (B1) at (0.5,0.5);
\coordinate (B2) at (0.5,-0.5);

\draw[fill=blue, opacity=0.4] (B1)--(A2)--(A1); 
\draw[fill=blue, opacity=0.4] (B1)--(A2)--(A3);
\draw[fill=blue, opacity=0.4] (B1)--(A3)--(A4);
\filldraw[fill=blue, opacity=0.4] (B2)--(A2)--(A1); 

\begin{scope}[thick,opacity=1]
\draw[dashed](A1)--(A4)--(B2)--(A3);

\draw[very thick, blue](A1)--(A2)--(A3);
\draw[very thick, blue](B1)--(A1);
\draw[very thick, blue](B1)--(A3);

\draw[very thick, blue](A2)--(B2);
\draw[very thick, blue](B1)--(A2);
\draw[very thick, blue](A1)--(B2);
\draw[very thick, blue](A4)--(B1);
\draw[very thick, blue](A4)--(A3);
\end{scope}

\node[font=\footnotesize, below] at (B2) {$x_{10}$};
\node[font=\footnotesize, above] at (B1) {$x_{11}$};
\node[font=\footnotesize, left] at (A1) {$x_{20}$};
\node[font=\footnotesize, above] at (A2) {$x_{30}$};x
\node[font=\footnotesize, right] at (A3) {$x_{21}$};
\node[font=\footnotesize, above] at (A4) {$x_{31}$};

\end{tikzpicture}

\caption{Triangulated balls defined by  $J_1$ and $J_2$}
\label{FigOct}
\end{figure}

\begin {exam}
   We continue Example \ref{ex:intro-m2}.
    In  Figure \ref{FigOct} we have the embeddings on the octahedron
   of the two hexagon triangulations in Figure \ref{fig:HexPol}.
 \end {exam}

 The essential element in proving Theorem \ref{thm:rainbow-main} is
 to show the following:

\medskip
\noindent {\bf Theorem \ref{thm:rain-const}}
{\it Any full-dimensional Cohen-Macaulay
proper sub-complex on the simplicial sphere
$\Delst(\Ab)$ is a constructible triangulated ball.
}

\medskip
In fact, such sub-complexes $X$ are precisely those that via the
Stanley-Reisner correspondence come from polarizations of Artin
monomial ideals. 
As noted above
  (see \eqref{eq:intro-bfa}), the facets of
  $X$ are then indexed by a subset $T \sus A_1 \times \cdots \times A_n$.
  Thus we may write $X = X(T)$, and its Stanley-Reisner ideal $J = J(T)$.

  It is noteworthy that by \cite{Hall} there are such sub-complexes $X$
  of $\Delst(\Ab)$ which
  are {\it not shellable}, see Proposition \ref{pro:link-hall}. But by
  the above, they are constructible. 

  \medskip
\noindent {\bf Dual geometric setting.}
  A standard fact (see Subsection \ref{subsec:pre-const})
  is that $J$ is CM iff its Alexander dual $I$ is an ideal with linear
  resolution. We write $I(T)$ for the Alexander dual of $J(T)$. 
    Considering each $A_i$ as a color class, each
  $\bfa \in A_1 \times \cdots \times A_n$ gives
  a {\it rainbow monomial} $m(\bfa) = \prod_{i=1}^n a_i$.
  The ideal $I(T)$ is generated by the $m(\bfa)$ for $\bfa \in T$.
  The subset $T \sus A_1 \times \cdots \times A_n$ identifies with a subset
of the vertices of the product polytope
  \[ \Pi(\Ab) = \Delta(A_1) \times \cdots \times \Delta(A_n). \]
  The vertices of $\Pi(\Ab)$ are naturally labelled by rainbow monomials.
  The subset $T$ induces a {\it saturated} cell sub-complex $C(T)$ of $\Pi(\Ab)$
  consisting of all sub-cells $\Pi(\Rb)$ of $\Pi(\Ab)$ (for non-empty subsets
  $R_i \sus A_i$) such that the product set
  $R_1 \times \cdots \times R_n \sus T$.
  We prove that when $I(T)$ has a linear resolution,
  it admits a cellular minimal
  free resolution supported by the cell sub-complex $C(T)$, Theorem
  \ref{thm:cell-linres}. As consequence, we get a homological
  criterion on $T$, Corollary \ref{cor:cell-homcrit},
  for $T$ to correspond to a polarization of an Artin
  monomial ideal. The maximal cells of $C(T)$ correspond bijectively
  to the components in the unique irreducible monomial
  decomposition of the Artin ideal $\ijfr$. Analysing then the possible
  sub cell-complexes $C(T)$ corresponding to this decomposition, enables in
  many example cases a full understanding of the possible polarizations
  of the Artin monomial ideal $\ijfr$, Section \ref{sec:hypercube}. 

\begin{exam}
  \[ \ijfr = (x^2, y^2, z^2, w^2) + (xz, xw, yz,yw, zw). \]
  The simplicial complex associated to the square-free part here is
  given in Figure \ref{fig:intro-sikx}.

\begin{figure}[ht]
\begin{center}  
\begin{tikzpicture}[thick, dot/.style={draw,fill,circle,inner sep=1.5pt},scale=1.5]

\node[dot, label=left:{$y$}] (n1) at (0,0) {};
\node[dot, label=left:{$x$}] (n2) at (0,1) {};
\node[dot, label=right:{$w$}] (n3) at (1,0) {};
\node[dot, label=right:{$z$}] (n4) at (1,1) {};

\draw (n1) edge[] (n2);
\end{tikzpicture}
\end{center}
\caption{Simplicial complex of ideal $(xz,xw,yz,yw,zw)$}
\label{fig:intro-sikx}
\end{figure}

The associated cell complex $C(T)$ will then consist of
a square and two edges (in general a $(d-1)$-simplex is turned into
a $d$-hypercube, when each $\alpha_i = 2$).
The cell complex must be acyclic. This gives, up to
isomorphism, the four
possibilities for $C(T)$ given in Figure \ref{fig:intro-lpp}.

\begin{figure}[ht] 
\begin{center}
\begin{tikzpicture}[thick,scale=1]

\coordinate(A1) at (0,0);
\coordinate(A2) at (0,1);
\coordinate(A3) at (1,1);
\coordinate(A4) at (1,0);
\coordinate(A5) at (1,2);
\coordinate(A6) at (2,1);

\draw[fill=lightgray, opacity=0.4] (A1)--(A2)--(A3)--(A4)--(A1);

\draw (A1) edge["$y$"] (A2);
\draw (A2) edge["$x$"] (A3);
\draw (A3) edge["$y$"] (A4);
\draw (A4) edge["$x$"] (A1);
\draw (A3) edge["$w$", swap] (A5);
\draw (A3) edge["$z$", swap] (A6);
\end{tikzpicture}
\hspace{0.5cm}
\begin{tikzpicture}[thick,scale=1]

\coordinate(A1) at (0,0);
\coordinate(A2) at (0,1);
\coordinate(A3) at (1,1);
\coordinate(A4) at (1,0);
\coordinate(A5) at (2,0);
\coordinate(A6) at (2,1);

\draw[fill=lightgray, opacity=0.4] (A1)--(A2)--(A3)--(A4)--(A1);

\draw (A1) edge["$y$"] (A2);
\draw (A2) edge["$x$"] (A3);
\draw (A3) edge["$y$"] (A4);
\draw (A4) edge["$x$"] (A1);
\draw (A4) edge["$z$", swap] (A5);
\draw (A3) edge["$w$"] (A6);
\end{tikzpicture}
\hspace{0.5cm}
\begin{tikzpicture}[thick,scale=1]

\coordinate(A1) at (0,0);
\coordinate(A2) at (0,1);
\coordinate(A3) at (1,1);
\coordinate(A4) at (1,0);
\coordinate(A5) at (-1,0);
\coordinate(A6) at (2,1);

\draw[fill=lightgray, opacity=0.4] (A1)--(A2)--(A3)--(A4)--(A1);

\draw (A1) edge["$y$"] (A2);
\draw (A2) edge["$x$"] (A3);
\draw (A3) edge["$y$"] (A4);
\draw (A4) edge["$x$"] (A1);
\draw (A1) edge["$z$"] (A5);
\draw (A3) edge["$w$"] (A6);
\end{tikzpicture}
\hspace{0.5cm}
\begin{tikzpicture}[thick,scale=1]

\coordinate(A1) at (0,0);
\coordinate(A2) at (0,1);
\coordinate(A3) at (1,1);
\coordinate(A4) at (1,0);
\coordinate(A5) at (1,2);
\coordinate(A6) at (2,2);

\draw[fill=lightgray, opacity=0.4] (A1)--(A2)--(A3)--(A4)--(A1);

\draw (A1) edge["$y$"] (A2);
\draw (A2) edge["$x$"] (A3);
\draw (A3) edge["$y$"] (A4);
\draw (A4) edge["$x$"] (A1);
\draw (A3) edge["$z$", swap] (A5);
\draw (A5) edge["$w$"] (A6);
\end{tikzpicture}
\end{center}
\caption{Acyclic cell complexes of one square and two line segments}
\label{fig:intro-lpp}
\end{figure}

These correspond to the four, up to isomorphism, possible polarizations
of the Artin ideal $\ijfr$. They can be worked out to be:
  \begin{align*}
    J_1& = (x_0x_1, y_0y_1, z_0z_1, w_0w_1)
         + (x_0z_0,x_0w_0, y_0z_0, y_0w_0, z_0w_0)  \\
    J_2 & =(x_0x_1, y_0y_1, z_0z_1, w_0w_1)
         + (x_0z_0,x_0w_0, y_1z_0, y_0w_0, z_0w_0)  \\
   J_3 & =(x_0x_1, y_0y_1, z_0z_1, w_0w_1)
         + (x_1z_0,x_0w_0, y_1z_0, y_0w_0, z_0w_0)  \\
 J_4 & =(x_0x_1, y_0y_1, z_0z_1, w_0w_1)
         + (x_0z_0,x_0w_0, y_0z_0, y_0w_0, z_1w_0)     
  \end{align*}
This is explained in more detail in Example \ref{ex:cell-lpp}.
\end{exam}

\medskip
Let us inform on a line of research leading to our results.
In \cite{Bier}, T.Bier 
 constructed simplicial spheres
 from a simplicial complex $X$ on a set $V$. Essentially
 the idea is that $X$ and its Alexander dual should glue together
 to give a simplicial sphere, as in the topological setting of
 Alexander duality, \cite{Hat}, see \cite[Sec.5.6]{Mat} for an argument.

 In \cite{BZPS} Bj\"orner, Pfaffenholz, Sj\"ostrand and Ziegler
 generalized this from
 down-sets in Boolean posets (that is, simplicial complexes) to
 down-sets in any poset.
 In addition they proved (by a counting
 argument) that Bier
 spheres were so numerous that most of them
 could not be realized as convex polytopes.  Bier spheres
 were considerably generalized by D'Alì, Nematbakhsh and the first
 author in \cite{DFN}: Given any poset $P$ let
 $\cI$ be a finite down-set in the set of order-preserving maps
 $\Hom(P, \NN)$.
 This gives rise to the so called letterplace ideals
 $L(\cI;P)$, see Subsection \ref{subsec:app-lp} or \cite{FGH}, 
 which are SR-ideals of  triangulated balls and whose boundaries
 are simplicial spheres, again with a simply defined SR-ideal.
 Our main result here is a further considerable generalization of
 letterplace ideals.
 In particular since even the double specialization to
 Bier spheres is so numerous that they do not in general
 have convex realizations (though recently shown
 to be star-shaped \cite{JTR}),
 the same holds for the simplicial spheres given here.

 \medskip
 The organization of this article is as follows.

 \medskip
 \noindent {\bf Part I, Sections \ref{sec:pre}-\ref{sec:polball}:
   Polarizations and triangulated balls.} We give background
 and prove the main
 Theorem \ref{thm:rainbow-main}.
 
 Section \ref{sec:pre} recalls basics on simplicial complexes,
 Stanley-Reisner rings, Alexander duals, joins, and constructible
 monomial ideals.
Section \ref{sec:rbideal} gives the geometric setting, the join of
boundaries $\Delst(\Ab)$. The facets of these joins
correspond bijectively to elements $(a_1, \ldots, a_n)$
of $A_1 \times \cdots \times A_n$. This gives a {\it rainbow monomial}
$a_1 \cdots a_n$. Subsets $T$ of $A_1 \times \cdots \times A_n$ then generate
a {\it rainbow ideal} $I(T)$. The Alexander dual
ideals
of polarizations are rainbow ideals with linear resolution.
In Section \ref{sec:polball} we give the main result, Theorem \ref{thm:rainbow-main},
that polarizations of Artin monomial ideals give triangulated balls.
The argument is by showing that rainbow ideals with linear resolution are
constructible ideals. 

\medskip
\noindent {\bf Part II, Sections \ref{sec:cell}-\ref{sec:hypercube}:
  Dual cell complexes.}
We study the dual cell complex $C(T)$, which is a sub-complex of
the product polytope $\Pi(\Ab)$.

Section \ref{sec:cell} defines $C(T)$. When $I(T)$ has
linear resolution, we show that it admits a
cellular minimal free resolution supported by 
$C(T)$. We derive the homological criterion for
$T$ to correspond to a polarization. Such $T$ are called
{\it restriction-acyclic}.
Section \ref{sec:mid} shows that the maximal cells in $C(T)$ are
in bijection with the components in the unique minimal
irreducible decomposition of the Artin monomial ideal $\ijfr$.

How can we find all polarizations of the Artin ideal $\ijfr$?
Section \ref{sec:hypercube} addresses this when $\ijfr$ contains
the second power $x_i^2$ of all variables $x_i$. Such $\ijfr$ are
in bijection to simplicial complexes on vertices $\{1,2, \ldots, n\}$.
The product polytope $\Pi(\Ab)$ is then a hypercube. We show in a range
of examples how we may classify the polarizations of such $\ijfr$.

\medskip
\noindent {\bf Part III, Sections \ref{sec:link}-\ref{sec:Artin}: Linkage
  and complement triangulations.} We study the complement set $T^c$,
the complement triangulation of $X(T)$, and relations to linkage theory.

Section \ref{sec:link} shows that $T$ corresponds to a polarization
iff $T^c$ does. The Alexander duals of $I(T)$ and $I(T^c)$ are then
linked Cohen-Macaulay ideals. They define complement triangulated
balls $X(T)$ and $X(T^c)$ on the join polytope $\Delst(\Ab)$. 
Section \ref{sec:bound} describes the boundaries of
these triangulated balls. These are triangulated spheres whose
Stanley-Reisner ideals are then Gorenstein ideals.
Section \ref{sec:Artin} shows how the linkage descends to linkage
of Artin monomial ideals.

\medskip
\noindent {\bf Part IV, Sections \ref{sec:squeezed}-\ref{sec:further}:
  Squeezed spheres and further problems.}
 Section \ref{sec:squeezed} is
     somewhat independent from the rest. 
     We show that the
     squeezed balls and spheres of G.Kalai \cite{Ka}
     have Stanley-Reisner ideals which
     come from letterplace ideals, these ideals being polarizations of
     Artin monomial ideals.
     The last Section \ref{sec:further} sketches open research directions.
     In Subsection \ref{subsec:further-regquot}
     we indicate how most simplicial balls are not of the type
     we construct here. However, as we see in examples,
     they might be {\it derived} from the simplicial
     balls here by joining variables (as we show in Section \ref{sec:squeezed}
     for {squeezed balls and spheres}).

     \medskip
\noindent {\it Acknowledgments}:  
     The work on this paper started while the third author visited the first author at
the University of Bergen, as TMF-Guest Professor with a grant from the
Trond Mohn Foundation. 
The third author thanks the Department of Mathematics of the University of Bergen for its hospitality.
This work was also partially supported by a grant from the IMU-CDC, 

We are grateful to anonymous referees for substantial input
to improving the article, and in particular for providing the general 
statement of Proposition \ref{pro:link-xyz} and the example of
Proposition \ref{pro:link-hall}.

\section{Preliminaries}
\label{sec:pre}
We recall basics on simplicial complexes,
the Stanley-Reisner correspondence, and Alexander duality.
For introduction to these, see \cite{HeHi, Mi-St} or \cite[Sec.II.5]{BH93}.

\subsection{Simplicial complexes and squarefree monomial ideals}
Let $V$ be a finite set. A simplicial complex on $V$ is
a family $X$ of subsets of $V$ such that if $F \in X$
and $G \sus F$, then also $G \in X$. The $F$ are called the
{\it faces} of $X$ and the dimension of $F$ is the cardinality $|F| -1$.
Maximal faces are called {\it facets}. If all facets have the same
dimension, $X$ is said to be {\it pure}.
If $X$ is the full power set of $V$, then
$X$ is called the {\it simplex} on $V$ and is denoted $\Delta(V)$. 
The {\it boundary $\pd X$} of a pure
simplicial complex $X$ is the smallest simplicial
complex containing all faces of $X$ which are on only one facet, and
are properly contained in this facet.
A simplicial complex $X$ has a natural topological realization $|X|$ in the
space ${\mathbb R}^V$, \cite{Mu}.

Let $k$ be a field and $k[V]$ the polynomial ring with elements of $V$
as variables. If $R \sus V$ denote by $m(R) = \prod_{v \in R} v$ the
associated monomial. Such a monomial is said to be {\it squarefree},
and if an ideal is generated by such monomials, it is also said to be
{\it squarefree}.

Given a simplicial complex $X$ we have its associated Stanley-Reisner (SR)
ideal $I= I_X$ which is the ideal generated by the monomials
$m(R)$ such that $R$ is {\it not} in $X$. The quotient ring $k[V]/I$ is
the Stanley-Reisner ring of $X$.

Simplicial complexes form a lattice with join $X \cup Y$, the union
of the families $X$ and $Y$, and meet $X \cap Y$.
The squarefree monomial ideals also form a lattice
with join the sum of ideals, and meet the intersection of ideals.

The Stanley-Reisner correspondence gives
a one-to-one correspondence between the lattice of
simplicial complexes on $V$
and the lattice of squarefree monomial ideals in $k[V]$, with order relation
reversed. Hence the join $X \cup Y$ corresponds to the meet of ideals
$I_X \cap I_Y$, and the meet $X \cap Y$ corresponds to the sum of ideals
$I_X + I_Y$.

\subsection{Alexander duality}
\label{subsec:AD}
The Alexander dual monomial ideal of $I (=I_X)$
is the monomial ideal $J = I^\vee$ generated
by all monomials $m$ which have a non-trivial common divisor with every monomial
in $I$. This is an involutive correspondence. The
squarefree monomials in $J$ are the $m(F^c)$ where the $F$ are faces
of $X$, and $F^c = V \backslash F$ its complement. In particular the
minimal generators
of $J$ are the monomials $m(F^c)$ of the complements of facets $F$
of $X$, see \cite[Lemma 1.5.3]{HeHi}.

Alexander duality gives a bijection between lattices of ideals
with order relation reversed. Hence the Alexander dual of the join
$I + L$ is the meet $I^\vee \cap L^\vee$, and the dual of $I \cap L$ is
$I^\vee + L^\vee$.

\subsection{External join of simplicial complexes}

Let $X$ be a simplicial complex on $V$, and $Y$ a simplicial complex
on a disjoint vertex set $W$. The (external) join of $X$ and $Y$ is
the simplicial complex $X * Y$ on $V \cup W$ given by:
\[ X * Y = \{ A \cup B \, | \, A \in X, B \in Y \}. \]
(This is not to be confused with the (internal) join in the lattice
of simplicial
complexes on a fixed set $V$, which is simply the union of $X$ and $Y$.
Henceforth join will always mean external join.)
The SR-ideal of the join $X * Y$ is the ideal
$(I_X) + (I_Y)$ generated by $I_X$ and $I_Y$ in $k[V \cup W]$. 

The topological realization $|X * Y|$ of the join, is obtained by connecting
points
in $|X|$ and $|Y|$ by lines. So it is the set
$|X| \times [0,1] \times |Y|$ modulo identifying pairs
$(x,0,y)$ and $(x,0,y^\prime)$ and identifying pairs
$(x,1,y)$ and $(x^\prime, 1, y)$.
The join of two spheres $S^n$ and $S^{m}$ of dimensions $n$ and $m$,
is the sphere $S^{n+m+1}$. If $|Y|$ is two points (corresponding
to the simplicial complex $\{\{v_1\}, \{v_2\}, \emptyset\}$, which
is topologically the sphere $S^0$), the
join $X * Y$ has topological realization the suspension $S|X|$.
The suspension of $S^n$ is $S^{n+1}$.

\subsection{Constructible squarefree monomial ideals}
\label{subsec:pre-const}
A simplicial complex of dimension $d$ is {\it constructible} if either
\begin{itemize}
\item $X$ is a simplex of dimension $d$, or
\item there is a decomposition $X = Y \cup Z$ where $Y$ and $Z$
  are constructible of dimension $d$ and $Y \cap Z$ is constructible
  of dimension $d-1$.
\end{itemize}
The following implications are well known:
\[ \text{vertex decomposable} \Rightarrow
  \text{shellable} \Rightarrow \text{constructible} \Rightarrow
  \text{Cohen-Macaulay/any field}. \]
Let $I$ be the Alexander dual ideal of the Stanley-Reisner ideal $I_X$.
Translating the above to conditions on the ideal $I$ we get the following.
A squarefree monomial ideal $I$ is {\it $e$-constructible}
(notion introduced in \cite{Olt-Const}) if:
\begin{itemize}
\item $I = (m)$ is generated by a monomial of degree $e$, or
\item there is a decomposition $I = J + K$ such that $J$ and $K$
  are $e$-constructible, and $J \cap K$ is $(e+1)$-constructible.
\end{itemize}

\begin{rema} If $I$ corresponds to the simplicial complex $Y$ via
  the SR-correspondence, one should note that $I$ being constructible
  does {\it not} correspond to $Y$ being constructible. Rather $I$ being
  constructible corresponds to 
  the Alexander dual simplicial complex $Y^\vee$ being constuctible
  (where  $Y^\vee$ has the Alexander dual ideal $I^\vee$ as
  Stanley-Reisner ideal).
\end{rema}

The Eagon-Reiner theorem, \cite{ER98}, states that a squarefree
monomial ideal $I$ has a linear resolution if and only if the Alexander dual $I^\vee$
is a Cohen-Macaulay ideal, see standard textbooks \cite[Thm.8.1.9]{HeHi},
\cite[Thm.5.56]{Mi-St}.
As a consequence constructible ideals have linear resolutions.
In general for ideals there are implications:
\[  \text{vertex splittable} \Rightarrow
  \text{linear quotients} \Rightarrow \text{constructible} \Rightarrow
  \text{linear resolution/any field}. \]
The notion of vertex splittable is introduced in \cite{MoKa}, and
linear quotients in \cite{HT-cones}.

\medskip
The statement of the  following lemma is analogous to the idea of
constructibility. We state it
as we need it later for an important example in Proposition \ref{pro:link-hall}.

\begin{lemma} \label{lem:pre-CM}
  Let $X$ and $Y$ be Cohen-Macaulay (CM) simplicial
  complexes of dimension $d$ on the same vertex set. Suppose
  the intersection $X \cap Y$ has dimension $< d$.
  Then the union $X \cup Y$ is CM iff
  the intersection $X \cap Y$ is CM of dimension $d-1$.
\end{lemma}

\begin{proof}
  A simplicial complex $X$ of dimension $d$ is CM iff
  the local cohomology groups $H^i_{\mathfrak m}(k[X])$ vanish for
  $i \leq d$, and is non-vanishing for $i = d+1$, \cite[A.7]{HeHi}.
  We show that there is an isomorphism
  \begin{equation} \label{eq:pre-hp}
    H^p_{\mathfrak m}(k[X \cap Y]) \iso H^{p+1}_{\mathfrak m}
    (k[X \cup Y]), \quad p \leq d-1, \end{equation}
  which will prove the claim (since for $e$ the dimension of $X \cap Y$,
  the group $H^{e+1}_{\mathfrak m}(k[X \cap Y])$ is non-vanishing).

  Consider the exact sequences:
  \begin{equation} \label{eq:pre-eksakt}
    0 \pil K_1 \pil k[X \cup Y] \pil k[X] \pil 0,
    \quad  0 \pil K_2 \pil k[Y] \pil k[X \cap Y] \pil 0.
    \end{equation}
  Here
  \[ K_1 = \frac{I_X}{ I_{X \cup Y}}
    = \frac{I_X}{ I_{X} \cap I_{Y}}
    \iso \frac{I_X + I_Y}{I_Y} =
    \frac{I_{X \cap Y}}{I_Y} = K_2.\]
    Whence taking long exact cohomology sequences of the two short exact
    sequences \eqref{eq:pre-eksakt}, we get \eqref{eq:pre-hp}.
  \end{proof}

\section{Polarizations and their Alexander duals}
\label{sec:rbideal}

In polarizations of Artin ideals $\ijfr$ in $k[x_1, \ldots, x_n]$,
the variables in the polynomial ring of the polarization $J$,
come in classes $A_1, A_2, \ldots, A_n$. The variables in $A_i$
map to $x_i$ and we call the $A_1, \ldots, A_n$ color classes. 
This polarization $J$ defines a full-dimensional sub-complex
of the polytope whose boundary is a join of simplices.

We show that the Alexander dual of $J$ is a {\it rainbow ideal}:
It is generated
by monomials which have one variable from each color class. These
ideals enables us in the next section to prove our main result.

\subsection{Joins of boundaries of simplices}
Recall that $\Delta(V)$ is the simplex on the finite set $V$. 
The SR-ideal generated by the single monomial $m(V) = \prod_{v \in V} v$
defines the
boundary $\pd \Delta(V)$ of $\Delta(V)$, which is topologically
a sphere $S^{|V|-2}$ where $|V|$ is the cardinality of $V$. 

Let  $A_1, A_2, \ldots, A_n$ be disjoint finite sets, and $A = \sqcup_{i=1}^n A_i$
their union. The complete intersection ideal
\begin{equation*} 
  C = (m(A_1), m(A_2), \ldots, m(A_n))  \sus k[A]
\end{equation*}
is the SR-ideal of the join
of the boundaries of simplices
\begin{equation} \label{eq:liaison-join}
  \Delst(\Ab) = \pd \Delta(A_1) * \pd \Delta(A_2) * \cdots * \pd \Delta(A_n),
 \end{equation}
which topologically is a join of spheres, and so also a sphere.
Its dimension is $\sum_{i = 1}^n |A_i| -1-n$. If each cardinality
$|A_i| = 2$, then each boundary $\pd \Delta(A_i)$ is two points,
and we are taking iterated suspensions. This gives a sphere
of dimension $n-1$: For $n = 2$, (the boundary of) the square,
for $n = 3$
(the boundary of) the octahedron, and in general
(the boundary of) the hyperoctahedron.

\subsection{Rainbow ideals}
Denote $\bfA = A_1 \times A_2 \times \cdots \times A_n$.
For $\bfa = (a_1, a_2, \ldots, a_n)$ in $\bfA$, denote
by $[\bfa ]$ the set $\{a_1, a_2, \ldots, a_n \}$. 
The faces of $\Delst(\Ab)$ are all subsets $S \sus A$ such
that no $A_i \sus S$. In other words it is all $S$ such that
each $A_i \backslash S$ is non-empty. Equivalently we 
may find $\bfa \in \bfA$ with $[\bfa] \cap S = \emptyset$.
The facets of $\Delst(\Ab)$ are then precisely
the complement sets $[\bfa]^c$ for $\bfa \in \bfA$.

The following is of course true for any triangulated sphere, but the argument
is explicit and informative in our setting.

\begin{lemm} \label{lem:link-cdto}
  Any codimension one face of $\Delst(\Ab)$ is on exactly two facets.
\end{lemm}

\begin{proof} A codimension one face of $\Delst(\Ab)$ is
  $F = [\bfa]^c \backslash \{a^\prime \}$ for some $a^\prime \in A$.
  So there is an $i$ such that $a^\prime \in A_i$, and write $a_i^\prime
:= a^\prime$. Let $\bfa^\prime = (a_1, \ldots, a_i^\prime, \ldots, a_n)$.
  Then $F$ is precisely on the two facets $[\bfa]^c$ and $[\bfa^\prime]^c$.
\end{proof}

For $\bfa \in \bfA$ let $m(\bfa)$ be the monomial 
$\prod_{i = 1}^n a_i$. Thinking
of $A_i$ as a color class, we call such a monomial a {\it rainbow
monomial}.  The Alexander dual ideal $B$ of the complete intersection
$C = (m(A_1), \ldots, m(A_n))$ is generated by the complements
of the facets of $\Delst(\Ab)$, and so it is the monomial
ideal generated by all the  rainbow monomials. Letting
$(A_i)$ be the ``color ideal'' generated by the elements of $A_i$ (considered
as variables), then $B$ is the product of these color ideals:
\[ B  = (m(\bfa) \, | \, \bfa \in \bfA) = \prod_{i = 1}^n(A_i). \]

\begin{defi}
For a subset $T \sus \bfA$, denote by $I(T)$ the ideal generated
by $m(\bfa)$ for $\bfa \in T$. Such an ideal is called a {\it rainbow ideal}.
For the complement set $T^c = \bfA \backslash T$, we
usually write $L(T)$ for $I(T^c)$, the {\it complement} rainbow ideal.
\end{defi}

\begin{rema}  The ideal $B = I(\bfA)$ is the squarefree monomial ideal
  defining the Cox irrelevant ideal for the toric variety which
  is the multiprojective space
  \[ {\mathbb P}^{|A_1|-1} \times \cdots \times {\mathbb P}^{|A_n|-1}. \]
  In this setting \cite{BESVirt} introduces {\it virtual resolutions},
  i.e. resolutions with homology annihilated by powers of $B$.
  Combinatorial aspects of this, in particular Stanley-Reisner rings
  that are virtually Cohen-Macaulay, are studied in \cite{BerComb}.
  \end{rema}
 
\begin{rema} For $T \sus \bfA$ the family of subsets $[\bfa] \sus A$
  where $\bfa \in T$ are precisely the collection of edges
  of $n$-partite $n$-uniform hypergraphs, also called
   $n$-partite $n$-uniform clutters, \cite{AN}.
   They are also called multipartite uniform hypergraphs.
   
    In \cite{NR09} a class of  $d$-partite $d$-uniform hypergraphs,
    Ferrer's hypergraphs
  are studied. The associated rainbow ideals are shown to have a linear
  resolutions, and an explicit form as a cellular resolution by
  the ``complex of boxes'' is given.
  These ideals are again a special case of co-letterplace ideals
  \cite{FGH}, see Subsection \ref{subsec:further-acyclic},
  which also have linear resolutions of an explicit form given
  in \cite{DFN}. \cite{AN} uses the terminology $d$-partite $d$-uniform
  clutters, and explicitly describes the first linear strand of
  rainbow ideals $I(T)$.
 \end{rema}
    
 \begin{rema}
    $d$-partite $d$-uniform hypergraphs have been a topic of interest
    in combinatorics and combinatorial geometry. In \cite{Rec-GS} they
    look at the associated embedding by generic points in ${\mathbb R}^d$, 
    relating it to the colored Tverberg's theorem \cite{Mat}.
    Matchings in such hypergraphs are studied in \cite{Mat-DKPU, Beck}.
    Mengerian $r$-uniform hypergraphs are shown to be $r$-partite in
    \cite{BeSc}, but the converse implication does not hold. From
    a commutative algebra perspective, \cite{AlBa}
    relates this to symbolic powers of ideals, and identify the
    subclass of tripartite tri-uniform hypergraphs for which $I^{(2)} = I^2$.  
\end{rema}
    
  \subsection{Polarizations}
  \label{subsec:pol}

  Recall from the introduction, Definitions \ref{def:intro-sep} and
  \ref{def:intro-pol} of separations and polarizations of monomial ideals.

  Let $\ijfr$ be an Artin monomial ideal
  in $k[x_1, \ldots, x_n]$, containing the Artin complete
  intersection
  \[ \ic = (x_1^{\al_1}, \ldots, x_n^{\al_n}). \]
  We do not require the powers $x_i^{\alpha_i}$ 
  to be minimal generators of the ideal $\ijfr$. 
  Let $A_1, \ldots, A_n$ be disjoint sets with $A_i$ of cardinality $\alpha_i$,
  and $A = \cup_i A_i$.
We have a map of polynomial rings
  $k[A] \pil k[x_1, \ldots, x_n]$ induced by the map 
  $A \pil \{x_1,x_2, \ldots, x_n \}$ sending elements of $A_i$
  to the variable $x_i$.
  By \cite[Remark 2.5]{AFL} any proper polarization $J^\prime$
  of an Artin monomial ideal
  containing the complete intersection $\ic$ is, after possibly
  extending $J^\prime$, isomorphic to a polarization $J$ contained in $k[A]$. 
  
  We also recall the following, which is essentially \cite[Prop.3.6]{AFL}.
  Due to its importance here,
  we include the proof for the benefit of the reader.
  
\begin{prop} \label{pro:pol-cm}
  Let $C = (m(A_1), \ldots, m(A_n))$ be the complete intersection.
  Polarizations $J \sus k[A]$ of Artin monomial ideals containing the complete
  intersection $\ic$
  are precisely the Cohen-Macaulay
  squarefree ideals $J \sus k[A]$ of height $n$  containing  $C$.

  Such ideals $J$ are in turn precisely 
    the Alexander duals of rainbow monomial ideals in $k[A]$ with
    linear resolution.
    \end{prop}

    \begin{proof} 
      i. Assume $J$ is a polarization of an Artin
      monomial ideal containing $\ic$. Then $J \sus k[A]$ contains $C$, as
      $m(A_i)$ must be the monomial mapping to $x_i^{\al_i}$.
      Since $C$ is cut down to an Artin monomial ideal, $J$ is also cut
      down by the regular sequence (for $k[A]/J$) to an Artin monomial ideal.
      An Artin ideal is Cohen-Macaulay, and so the same holds for $J$ by
      \cite[Thm.2.1.3.a]{BH93}.

      ii. Conversely assume $J$ is Cohen-Macaulay and squarefree of height $n$
      containing $C$. The sequence 
      of variable differences cuts $C$ down to an Artin monomial ideal,
      and so also cuts $J$ down to an Artin monomial ideal $\ijfr$. Thus
      starting from $k[A]/J$, which has dimension $|A|-n$, 
      at each step when we divide out by a variable difference,
      dimension is cut by one. By \cite[Thm.2.1.2.c]{BH93}
      the sequence is regular for $k[A]/J$ 
      and so $J$ is a polarization of $\ijfr$. 
      
      By the Eagon-Reiner theorem, see Subsection
      \ref{subsec:pre-const},
      $J$ is Cohen-Macaulay iff its Alexander dual $I$ has a linear resolution.
      Let $X$ be the simplicial complex associated to $J$. That $J$ is
      unmixed of height $n$ is equivalent to  the facets of $X$ being
      facets of $\Delst(\Ab)$, and so $I = I(T)$ for some $T \sus \bfA$.
    \end{proof}

 \begin {exam} \label{eks:polar-fem}
    The standard polarization of $(x_1^2, x_1x_2^2, x_2^3)$ is
    \[ (x_{10}x_{11}, x_{10}x_{20}x_{21}, x_{20}x_{21}x_{22}). \]
    (Standard polarization means that each power $x_i^r$
    is replaced by the product $x_{i0} x_{i1} \cdots x_{i, r-1}$.)
    It corresponds to the simplicial complex to the left in
    Figure \ref{pol:123}. To the right it is drawn in blue on
    $\vardel \Delta(A_1) * \vardel \Delta(A_2)$
where $|A_1| = 2$ and $|A_2| = 3$.
     It consists of five of the six facets. Only the lower front
     facet, with vertices $x_{20}, x_{21}$ and $x_{10}$,
     is not in the simplicial complex.
    
 \begin{figure}
\begin{center}   
\begin{tikzpicture}[dot/.style={draw,fill,circle,inner sep=1pt},scale=1.5]
  \draw [help lines,white] (-1.5,-1.5) grid (1.5,1.2);

\coordinate(A1) at (0,-0.3);
\coordinate (A2) at (2,-0.3);
\coordinate (B1) at (1,0.2);
\coordinate (B2) at (1,0.6);
\coordinate (B3) at (1,1);

\draw[thick] (A1) -- (A2) -- (B3) -- (A1);
\draw[thick] (A1) -- (B2) -- (A2);
\draw[thick] (A1) -- (B1) -- (A2);
\draw[thick] (B1) -- (B2) -- (B3);
);
\draw[fill=gray, fill opacity=0.3](A1) -- (A2) -- (B3) -- (A1);

\node[font=\footnotesize, left] at (A1) {$x_{20}$};
\node[font=\footnotesize, right] at (A2) {$x_{21}$};
\node[font=\footnotesize, above] at (B3) {$x_{10}$};
\node[font=\footnotesize, right] at (B2) {$x_{22}$};
\node[font=\footnotesize, right] at (B1) {$x_{11}$};

\coordinate(A1) at (4,0);
\coordinate (A2) at (6.7,-0.1);
\coordinate (A3) at (5.6,0.53);
\coordinate (B2) at (5.33,-0.8);
\coordinate (B1) at (5.47,1.33);

\draw[fill=blue,opacity=0.3] (A1)--(A2)--(B1);
\draw[fill=blue,opacity=0.35] (A2)--(B1)--(A3);
\draw[fill=blue,opacity=0.35] (A1)--(A3)--(B1);
\draw[fill=blue,opacity=0.2] (A1)--(B2)--(A3);
\draw[fill=blue,opacity=0.2] (A2)--(B2)--(A3);

\draw[blue](A1)--(A2)--(B1)--(A1);
\draw[blue] (A2)--(B1);
\draw[blue] (A1)--(A2)--(B2)--(A1);

\draw[dashed, black, opacity=0.6] (A1)--(A3)--(B2);
\draw[dashed, black, opacity=0.6] (A2)--(A3)--(B1);

\node[font=\footnotesize, below] at (A3) {$x_{22}$};
\node[font=\footnotesize, above] at (B1) {$x_{11}$};
\node[font=\footnotesize, left] at (A1) {$x_{20}$};
\node[font=\footnotesize, right] at (A2) {$x_{21}$};
\node[font=\footnotesize, right] at (B2) {$x_{10}$};

\end{tikzpicture}
\end{center}
\caption{Triangulated ball of Example \ref{eks:polar-fem}}
\label{pol:123}
\end{figure}
\end {exam}
    
    \begin {exam} \label{ex:pol-J2}
  The ideal $J_2$ of Example \ref{ex:intro-m2} is a polarization of
  $(x_1,x_2,x_3)^2$ containing the complete intersection $(x_1^2, x_2^2,x_3^2)$.
  We can let
  \[ A_1 = \{x_{10}, x_{11}\},\quad  A_2 = \{x_{20}, x_{21}\}, \quad A_3
    = \{ x_{30}, x_{31} \}. \]
  The Alexander dual ideal of $J_2$ is the rainbow ideal $I_2$ with generators:
  \[ x_{10}x_{20}x_{30}, x_{10}x_{20}x_{31}, x_{10}x_{21}x_{31},
    x_{11}x_{21}x_{31}. \]
  We may identify each of $A_1,A_2,A_3$ with $\{0,1\}$, and then the
  ideal $I_2$ is $I(T)$ where
  \[T = \{000, 001, 011, 111 \}. \]
\end {exam}
  
 \begin{rema}
    The polarizations of powers of graded maximal ideals $(x_1, \ldots, x_n)^m$
    are classified by Almousa, Lohne and the first author in \cite{AFL}.
     In particular
polarizations of $(x_1, \ldots, x_n)^2$ are in bijections with trees
with $n$ edges. These results in \cite{AFL} are  
     partially extended to strongly stable ideals generated in a
     single degree in 
     Hao, Luo and Saint Rain \cite[Prop.4.2, Thm.6.1]{HLSR}.
     Note that such strongly stable ideals are not Artin unless they are a power
     of the graded maximal ideal.
     \cite{AB} gives resolutions of polarizations of $(x_1, \ldots, x_n)^m$,
     connecting
     the combinatorics of \cite{AFL} with that of Hook tableaux.
   \end{rema}
   
\newcommand{\bfr}{{\mathbf r}}
\newcommand{\bfq}{{\mathbf q}}

\section{Polarizations define triangulated balls}
\label{sec:polball}

The previous section shows that
polarizations of
Artin monomial ideals in $k[x_1, \ldots, x_n]$
containing the complete intersection $\ic = (x_1^{\alpha_1},
\ldots, x_n^{\alpha_n})$
correspond to rainbow ideals in $k[A]$ with linear resolutions,
for $A$ a union of color classes
$A_1, \ldots, A_n$ and cardinalities $|A_i|= \alpha_i$.
Rainbow ideals
are $I(T)$ for a subset $T \sus \bfA = A_1 \times \cdots \times A_n$.

We investigate ideals of type
\[ I = b_1 I_1 + b_2 I_2 + \cdots + b_m I_m \]
where the $I_j$ do not involve the $b$-variables, and the situation
when such ideals have a linear resolution.
We show that rainbow ideals have
linear resolutions if and only if they
are constructible. As consequence we get our main theorem, that
polarizations of Artin monomial ideals define triangulated
balls by the Stanley-Reisner correspondence.

\subsection{The main result}

The following is the key result leading to our main statement.
  
 \begin{theo} \label{thm:linres-const} For a subset
   $T \sus A_1 \times \cdots \times A_n$,
   the rainbow ideal $I(T)$ has $n$-linear resolution if and
   only if it is $n$-constructible.
  \end{theo}

  By Subsection \ref{subsec:pre-const}, if  $I(T)$ is constructible, then
  $I(T)$ has a linear resolution. The challenge is to show the other direction.
  But let us first give the consequence which is our main result.
  First recall the following, stated in the overview article
  \cite[Thm.11.4]{Bj95} by A.Bj\"orner:

   \medskip
  {\noindent \bf Theorem.}
  {\it Let $X$ be a pure constructible simplicial complex
    such that every codimension one face is on one or two facets.
    Then the topological realization $|X|$, is a ball or a sphere
    (the latter when every codimension one face is on exactly two
    facets).} 
\medskip
  
Using this, Theorem \ref{thm:linres-const} implies Conjecture 3.2 in \cite{AFL}.

\begin{theo}[Main Theorem] \label{thm:rainbow-main}
Let $\ijfr$ be an Artin monomial ideal in $k[x_1, \ldots, x_n]$
  strictly containing the
  complete intersection $\ic$ (but the $x_i^{\alpha_i}$ need not be
  minimal generators of $\ijfr$).
  Let $A \pil \{x_1, x_2, \ldots , x_n\}$  have
  fiber $A_i$ over $x_i$ of cardinality $\alpha_i$.
  
  Any polarization 
  $J \sus k[A]$ of $\ijfr$ defines, via the SR-correspondence,
  a  full dimensional  simplicial ball on the simplicial sphere $\Delst(\Ab)$.
  Moreover, any proper polarization of $\ijfr$ is,
  when followed by
  a suitable extension,
  isomorphic to such a $J$. 
\end{theo}

\begin{proof}
  By \cite[Remark 2.5]{AFL}, any proper polarization of $\ijfr$ will have
  at most $\alpha_i$ variables mapping to $x_i$. Hence it is, after
  possible extension, isomorphic to a polarization $J \sus k[A]$. 
  By Proposition \ref{pro:pol-cm} such a polarization is of the
  form $J = J(T)$, and is Cohen-Macaulay. Denote the
  associated  simplicial complex by $X(T)$. 
  The Alexander dual  $I(T)$ of $J(T)$ has a linear resolution. By Theorem
  \ref{thm:linres-const}, $I(T)$ is constructible and so $X(T)$ is
  constructible. By Lemma \ref{lem:link-cdto} every codimension one
  face of $X(T)$ is on at most two facets. Björner's theorem above
  gives that $X(T)$ is topologically a ball.
\end{proof}

We also get the following.

\begin{theo} \label{thm:rain-const}
 Any full-dimensional Cohen-Macaulay
proper sub-complex on the simplicial sphere
$\Delst(\Ab)$ is a constructible ball. It is obtained
as in Theorem \ref{thm:rainbow-main} above.
\end{theo}

\begin{proof}
  This follows by Proposition \ref{pro:pol-cm} and Theorem
  \ref{thm:linres-const}.
  \end{proof}

  At first we tried to show that such a sub-complex is shellable, which
  failed. In fact, as pointed out by an anonymous referee, such
  a sub-complex may not be shellable, see Proposition \ref{pro:link-hall}.

\subsection{Ideals with linear resolutions and constructability}
 
We proceed to prove Theorem \ref{thm:linres-const}.
First we give a general fact. Let $k[B \cup X]$ be a polynomial ring,
and $\hele^B$ the free abelian group with basis $e_b$ for $b \in B$.
We consider the polynomial ring to be $\hele^B \oplus \hele$ graded,
where the variables  $b \in B$ have degree $e_b$ and the variables
in $X$ have degree $e = {\mathbf 0} \oplus 1$.

The following is a crucial fact: If $p$ is a homogeneous
polynomial in $k[B \cup X]$ (for the grading by $\ZZ^B \oplus \ZZ$),
of degree $\bfr \oplus s$, then $p$ is a product
$\prod_{b \in B}b^{r_b} \cdot p^\prime$ where $p^\prime$ is a homogeneous
polynomial in $k[X]$ of degree $s$.

Let $M$ be a $\ZZ^B \oplus \ZZ$-graded 
$k[B \cup X]$-module. If $m \in M$ is homogeneous of degree $\bfr \oplus s$,
write $\deg_B(m) = \bfr$.
For $\bfr \in \ZZ^B$, denote $M_{\leq \bfr}$ be the
submodule of $M$ generated by the homogeneous elements $m \in M$
with $\deg_B(m) \leq \bfr$.
Note that if $F = \oplus_{j \in J} S(-{\mathbf u}_j)$ is a free module,
then $F_{\leq \bfr} = \oplus_{\deg_B({\mathbf u}_j) \leq \bfr} S(-{\mathbf u}_j)$.

\begin{lemma} Let $\phi : N \pil M$ be a morphism of $\ZZ^B \oplus \ZZ$-graded
  $k[B \cup X]$-modules, with $M$ a torsion free module. Let $\bfr \in \ZZ^B$. 
  If $K$ is the kernel of $\phi$, then $K_{\leq \bfr}$ is the kernel of
  $\phi_{\leq r} : N_{\leq \bfr} \pil M_{\leq \bfr}$.
\end{lemma}

\begin{proof} Clearly $K_{\leq \bfr}$ is contained in the kernel of
  $\phi_{\leq \bfr}$. 
  Suppose $n$ is homogeneous and
  in the kernel of $\phi_{\leq \bfr}$. We may write 
  $n = \sum p_i n_i$ where $p_i \in k[B \cup X]$ and the $n_i \in N$
  with $\deg_B(n_i) \leq \bfr$. Let $\deg_B(n) = \bfq$, 
  and consider the monomial
  $u = \prod_{q_b > r_b}b^{q_b - r_b}$. Since each $\deg_B(n_i) \leq \bfr$,
  we have each $\deg_B(p_i) \geq \deg(u)$. But then $p_i$ is
  divisible by $u$. Thus we may write $n = u \cdot n^\prime$
  with $\deg_B(n^\prime) \leq \bfr$. 

  Then $0 = \phi(n) = u \phi(n^\prime)$ and since $M$ is torsion
  free, this gives $\phi(n^\prime) = 0$, and so $n^\prime \in \ker \phi = K$.
  Thus  $n \in K_{\leq \bfr}$.
\end{proof}

This gives the following (a generalization of \cite[Thm.2.1]{GaHiPe}).

\begin{coro}
  Let $M$ be a $\ZZ^B \oplus \ZZ$-graded torsion free $k[B \cup X]$-module,
  and $F_\bullet$ a minimal free resolution of $M$. Let $\bfr \in \ZZ^B$.
  Then $F_{\bullet, \leq \bfr}$ is a minimal free resolution of $M_{\leq \bfr}$.
\end{coro}

We say $F_{\bullet}$ is a {\it $d$-linear resolution}, if
$F_i = \oplus_{j \in J} S(-{\mathbf u}_{ij})$ and for each
$\deg(u_{ij}) = \bfq \oplus s$, its total degee $s + \sum_{b \in B} q_b$
equals $d + i$.

\begin{coro} \label{cor:lin-Frlin}
  If $F_{\bullet}$ is a $d$-linear resolution of $M$, then
  $F_{\bullet, \leq \bfr}$ is a $d$-linear resolution of $M_{\leq \bfr}$.
\end{coro}

In the following we let $B = \{b_1, b_2, \ldots, b_m\}$ have
cardinality $m$. We write $\ZZ^B = \ZZ^m$ and $e_i = e_{b_i}$ for
$i = 1, \ldots, m$, and $e_{m+1} = {\mathbf 0} \oplus 1$.

\begin{prop} \label{pro:lin-ItoIR}
  Let $I_1, \ldots, I_m$ be ideals in $k[X]$.
Consider the product ideals $(b_j) (I_j) \sus k[B \cup X]$, which
we simply denote as $b_j I_j$. Also for $R \sus \{1,2, \cdots, m\}$
denote
\[ I_R = \cap_{j \in R} I_j  \sus k[X]. \]
If the ideal
\[ I = b_1 I_1 + b_2 I_2 + \cdots + b_m I_m \]
has $(d+1)$-linear resolution, then each $I_R$ has $d$-linear resolution.
\end{prop}

\begin{proof}
  First we claim: For each   $R \sus \{1,2, \ldots, m\}$
  the sub-sum  $\sum_{j \in R} b_j I_j$ has $(d+1)$-linear resolution.
  Let $\bfr = \sum_{j \in R} e_j $. Then $I_{\leq \bfr}$ is 
 $\sum_{j \in R} b_j I_j$, and by Corollary \ref{cor:lin-Frlin}
it has a linear resolution. 
As a special case each  $b_j I_j$ has $(d+1)$-linear
resolution, and so each $I_j$ has $d$-linear resolution.

We now use induction on $m$. If $m = 1$ the statement holds.
Let $m = 2$. We have
an exact sequence
\[ 0 \pil b_1 I_1 \cap b_2I_2 \pil b_1I_1 \oplus b_2I_2 \pil
b_1 I_1 + b_2 I_2 \pil 0, \]
and a long exact $\Tor$-sequence:
\begin{align*}  \cdots \pil & \Tor_{p+1}(b_1I_1 + b_2 I_2, k)_e \pil
  \Tor_{p}(b_1b_2 (I_1 \cap I_2), k)_e \pil \Tor_p(b_1I_1,k)_e \oplus
                           \Tor_p(b_2I_2,k)_e \\
  \pil & \Tor_{p}(b_1I_1 + b_2 I_2, k)_e \pil \cdots .
\end{align*}
Since $I = b_1I_1 + b_2 I_2$ has $(d+1)$-linear resolution, the Tor-groups
to the left vanish for $e - (p+1) \geq d+2$ or equivalently for $e-p \geq d+3$.
Similarly, to the right above the Tor-groups vanish for $e-p \geq d+2$.
Thus
\[ \Tor_{p}(b_1b_2 (I_1 \cap I_2), k)_e = 0, \quad e-p \geq d+3.\]
As $I_1 \cap I_2$ is generated in degree $\geq d$ it must have a
$d$-linear resolution.

\medskip
Assume now $m = 3$. From
\[ 0 \pil (b_1 I_1 +  b_2I_2) \cap b_3I_3 \pil (b_1I_1 + b_2I_2)
  \oplus b_3 I_3 \pil
b_1 I_1 + b_2 I_2 + b_3I_3 \pil 0, \]
the same argument gives that $(b_1I_1 + b_2 I_2) \cap (I_3)$ has
$(d+1)$-linear resolution. Let $b_1 i_1 + b_2 i_2$ be a generator for
this ideal in degree $d+1$. It is also in $(I_3)$, but since $I_3$ is generated
in degree $d$ and its generators do not involve the variables $b_1$ and $b_2$,
both $i_1, i_2 $ are in $(I_3)$. 
Therefore
\[ (b_1 I_1 + b_2I_2) \cap (I_3) = b_1 I_1 \cap I_3 + b_2 I_2 \cap I_3. \]
By the $m = 2$ case, both $I_1 \cap I_3$ and $I_2 \cap I_3$ have
$d$-linear resolution, as well as their intersection $I_1 \cap I_2 \cap I_3$.

For larger $m$ we split as $(b_1I_1 + \cdots + b_{m-1}I_{m-1}) + b_mI_m$, and
proceed in the same way as above.
\end{proof}

The converse of Proposition \ref{pro:lin-ItoIR} is not true
in general, but for monomial ideals it is.
We now proceed to do this. The key is the following:

\begin{lemma}
  Let $I$ and $J_1, \ldots, J_m$ be monomial ideals.
  Then
  \[ I \cap (J_1 + \cdots + J_m) = I \cap J_1 + \cdots + I \cap J_m. \]
\end{lemma}

\begin{proof}
  It is enough to prove this for $m = 2$, as the general case easily follows
  by induction. Suppose $f \in I \cap (J_1 + J_2)$. Then every monomial
  in $f$ is in $I$ and in $J_1 + J_2$. Then the monomial is in either
  $J_1$ or $J_2$. Whence $f = f_1 + f_2$ where $f_1$ consists of monomials
  in $I \cap J_1$ and $f_2$ of monomials in $I \cap J_2$.
\end{proof}

\begin{prop}
  Let $I_1, \ldots, I_m$ be monomial ideals in $k[X]$ such that for every
  $R \sus [m]$, the ideal $I_R = \cap_{j \in R} I_j$ has $d$-linear resolution.
  Then the ideal $I = \sum_{j=1}^n b_jI_j$ in $k[B \cup X]$ has
  $(d+1)$-linear resolution.
\end{prop}

\begin{proof}
  Again we use induction on $m$. Write
  \begin{equation} \label{eq:linres-KI}
    I = (b_1I_1 + \cdots + b_{m-1}I_{m-1}) + b_m I_m = K_{m-1} + b_mI_m.
    \end{equation}
    By induction we may assume $K_{m-1}$ has $(d+1)$-linear resolution. Note
    also that 
  $I_m$ has $d$-linear resolution. We have a short exact sequence:
  \begin{equation} \label{eq:linres-KIsekv}
    0 \pil K_{m-1} \cap b_mI_m \overset{f}\pil K_{m-1}
      \oplus b_m I_m \pil I \pil 0.
\end{equation}      
For the first term:
\begin{align*}
  K_{m-1} \cap b_mI_m&  = (b_1I_1 + \cdots + b_{m-1}I_{m-1}) \cap b_m I_m \\
                     &  = b_m(b_1I_1 \cap I_m + \cdots + b_{m-1} I_{m-1} \cap I_m)
\end{align*}
By induction the last ideal above has $(d+2)$-linear resolution. So
$K_{m-1} \cap b_m I_m$ has $(d+2)$-linear resolution. The mapping cone
of $f$ in the exact sequence \eqref{eq:linres-KIsekv}
then becomes a $(d+1)$-linear resolution of $I$.
\end{proof}

In fact, almost the same argument as above may be applied for
the notion of constructability.

\begin{prop} \label{pro:lin-constr}
  Let $I_1, \ldots, I_m$ be monomial ideals in $k[X]$. If for every
  $R \sus [m]$, $I_R$ is $d$-constructible, then
  \[ I = b_1 I_1 + \cdots + b_m I_m \]
  is $(d+1)$-constructible.
\end{prop}

\begin{proof}
  We write $I = K_{m-1} + b_mI_m$ as in \eqref{eq:linres-KI}.
  Here $b_m I_m$ is $(d+1)$-constructible
  and by induction on $n$, $K_{m-1}$ is also $(d+1)$-constructible.
  Further their intersection is:
  \begin{align*}  K_{m-1} \cap b_mI_m
    & = (b_1I_1 + \cdots + b_{m-1}I_{m-1}) \cap b_m I_m \\
    & = b_m(b_1I_1 \cap I_m + \cdots + b_{m-1} I_{m-1} \cap I_m),
\end{align*}
  and by induction, this is $(d+2)$-constructible.
  Thus $I$ is $(d+1)$-constructible.
\end{proof}

\begin{proof}[Proof of Theorem \ref{thm:linres-const}.]
  When we have one color class,  $n = 1$, $I$ is generated by variables and
  so is constructible. Let $n \geq 2$, and $A_1 = \{a_1, a_2, \ldots, a_m \}$.
   We may then write
  \[ I = a_1I_1 + \cdots + a_m I_m \] where the $I_j = I(T_j)$ for
  various $T_j \sus A_2 \times \cdots \times A_n$. Since $I$ has $n$-linear
  resolution, for every $R \sus [m]$, $I_R = \cap_{j \in R} I_j$ has
  $(n-1)$-linear resolution. As all the $I_j$ are monomial ideals
  generated in degree $(n-1)$, we must have $I_R = I(T_R)$ where
  $T_R = \cap_{j \in R} T_j$. By induction on $n$, each of these are
  $(n-1)$-constructible. Then by Proposition \ref{pro:lin-constr}, $I$ is
  $n$-constructible.
\end{proof}

\newcommand{\cpA}{\cup \Ab}
\newcommand{\prA}{\times \Ab}
\newcommand{\Cgeo}{C}
\newcommand{\susst}{\subseteq^*}
\newcommand{\Ctbu}{\tilde{C}_{\bullet}}
\newcommand{\Ctalg}{\tilde{{\mathcal C}}_{\bullet}}
\newcommand{\Ctcal}{\tilde{{\mathcal C}}}
\newcommand{\Calg}{{\mathcal C}_{\bullet}}
\newcommand{\Falg}{{\mathcal F}_{\bullet}}
\newcommand{\Fal}{{\mathcal F}}
\newcommand{\Cal}{{\mathcal C}}

\section{Cell complexes and minimal free resolutions
  of rainbow ideals}
\label{sec:cell}

The subset $T \sus A_1 \times \cdots \times A_n$ determines
a canonical cell complex $C(T)$. This cell complex has a
canonical monomial labelling. By classical constructions
\cite[Section 4]{Mi-St}
we obtain a cellular complex $\Ctbu(T)$. In our case this
is a linear complex. The beautiful thing is that when $I(T)$ has
linear resolution, the resolution is precisely this cellular complex,
Theorem \ref{thm:cell-linres}.

In particular, this gives an explicit homological criterion on $T$,
Corollary \ref{cor:cell-homcrit},
for when $J(T)$ is a polarization of an Artin monomial ideal.

\subsection{The cell complex and free
  cellular complex associated
  to $T$}

A finite set $A$ gives rise to the geometric
simplex spanned by the elements of $A$:
\[ \Delta(A) = \left\{\sum_{a \in A} \gamma_a \cdot a \, | \, 0 \leq \gamma_a \leq 1,  \sum_{a \in A} \gamma_a = 1 \right\} \sus \RR^A. \]
Given finite non-empty sets $ A_1, A_2, \ldots, A_n$, 
we get a polytope: 
\[ \Pi(\Ab) = \Delta(A_1) \times \Delta(A_2) \times \cdots
\times \Delta(A_n) \sus \times_{i = 1}^n \RR^{A_i} = \RR^A. \]
Its dual (or polar) polytope
is the join
\[ \Delst(\Ab) =  \pd \Delta(A_1) * \pd \Delta(A_2) * \cdots * \pd \Delta(A_n).
\]
The set of vertices of $\Pi(\Ab)$ is the set $A_1 \times \cdots \times
A_n$.

\begin{nota}
  When  $R_i \sus A_i, \, i = 1, \ldots, n$ are {\it non-empty} subsets,
  we write
$\Rb \susst \Ab$.
\end{nota}

\begin{defi}
  The polytope $\Pi(\Ab)$ gives rise to a cell complex $\Cgeo(\Ab)$,
  actually a polyhedral cell complex,
whose cells are the $\Pi(\Rb)$ with $\Rb \susst \Ab$. 
The dimension of such a cell is $\sum_i |R_i| - n$.
A subset $T \sus A_1 \times A_2 \times \cdots \times A_n$
induces a cell subcomplex $C(T)$ of $\Cgeo(\Ab)$ consisting of all
$\Pi(\Rb)$ such
$R_1 \times R_2 \times \cdots \times R_n \sus T$.
We say $C(T)$ is the {\it saturated} cell sub-complex of
$\Cgeo(\Ab)$ associated to $T$. 
\end{defi}

The cell complex $\Cgeo(T)$ gives rise to two important algebraic
complexes.

\medskip
\noindent{{\bf The (augmented) chain complex:}}
\begin{equation} \label{eq:minfree-Ct}
  \Ctbu(T): C_{-1}(T) \vpil C_0(T) \vpil C_1(T) \vpil
  \cdots \vpil C_p(T) \vpil \cdots ,
\end{equation}
where $C_p(T)$ is the $k$-vector space with basis all cells
of $\Cgeo(T)$ of dimension $p$ (and $C_{-1}(T)$ is the one-dimensional
vector space with basis the empty cell), and standard differentials
(see \eqref{eq:minfree-diff} below). The homology groups are the
reduced homology groups. Note that $\tilde{H}_{-1}(\Ctbu(T)) \neq 0$ iff
$T$ is empty.

\medskip
Let $S$ be the polynomial ring $k[A]$ whose variables
are the elements of $A = A_1 \cup \cdots \cup A_n$. 
A vertex $\bfa = (a_1,a_2, \ldots, a_n)$ of $\Pi(\Ab)$
may be labelled by the monomial
$m(\bfa) = a_1a_2\cdots a_n$. The cell $\Pi(\Rb)$ is then labelled by
the least common multiple of the monomials on its vertices. This is
simply
\[ m(\Rb) = m(R_1)m(R_2) \cdots m(R_n), \text{ where }
  m(R_i) = \prod_{r \in R_i}r,\]
which is a monomial of degree $\dim \Pi(\Rb) + n$. 

\medskip
\noindent{{\bf The (augmented) free cellular complex:}
  For $T \sus A_1 \times \cdots \times A_n$, this induces a labelling
by monomials of the cell sub-complex $\Cgeo(T)$.
By \cite[Section 4]{Mi-St} it then gives a {\it free cellular complex}
\[ \Ctalg(T) : \Cal_{-1}(T) \vmto{d_0} \Cal_0(T) \vmto{d_1}
  \cdots \vpil \Cal_p(T) \vpil \cdots , \]
where
$\Cal_p(T)$ is the free $S$-modules with generators all cells $\Pi(\Rb)$ in
$\Cgeo(T)$ of dimension $p$ (and $\Cal_{-1}(T) = S$).
The image of $d_0$ is $I(T)$.
We may put a total order on each $A_i$. Together with the ordered
sequence $A_1, A_2, \ldots, A_n$, this gives a total order on the
union $A = \cup_1^n A_i$. Then $\Rb \sus \Ab$ inherits a total order on
$R = \cup_1^n R_i$. For $b \in R$ denote by $\pot(b)$ its position in
this total order. For $p \geq 1$
the differential $d^p$ is then given by
\begin{equation} \label{eq:minfree-diff}
  \Pi(\Rb) \mapsto \sum_{\begin{matrix} i = 1, \ldots, n\\
      b \in R_i, |R_i| \geq 2 \end{matrix}}  (-1)^{\pot(b)} b \cdot
  \Pi(R_1, \ldots, R_i\backslash \{b \}, \ldots, R_n).
\end{equation}

\begin{exam} Let $A_1 = \{a_1, b_1\}$ and $A_2 = \{ a_2,b_2\}$ and
  $T = A_1 \times A_2$. The product polytope $\Pi (\Ab)$ is the
  square labelled as in Figure \ref{fig:cell-square}.

\begin{figure}[ht]
\begin{center}

\begin{tikzpicture}[thick,scale=2]

\coordinate[label=-135:{$a_1a_2$}] (A1) at (0,0);
\coordinate[label=135:{$a_1b_2$}] (A2) at (0,1);
\coordinate[label=45:{$b_1b_2$}] (A3) at (1,1);
\coordinate[label=-45:{$b_1a_2$}] (A4) at (1,0);

\draw (A1) edge[] (A2); 
\draw (A2) edge[] (A3); 
\draw (A4) edge[swap] (A3); 
\draw (A1) edge[swap] (A4); 

\draw[fill=lightgray, opacity=0.4] (A1)--(A2)--(A3)--(A4)--(A1);

\end{tikzpicture}

\end{center}
\caption{The square $\Pi(\Ab)$, with monomial labels on vertices}
\label{fig:cell-square}
\end{figure}
  The augmented chain complex is
  \[ k \vpil k^4 \vpil k^4 \vpil k, \]
  which has no homology, it is acyclic.
  Let $S = k[a_1,b_1, a_2, b_2]$. The rainbow ideal is
  $I(T) = (a_1a_2, a_1b_2, b_1a_2, b_1b_2)$. The augmented  free cellular complex is
  \[ S \overset{d}{\vpil} S(-2)^4 \vpil S(-3)^4 \vpil S(-4). \]
  The image of $d$ is $I(T)$, and in this case,
  omitting the first term,  the complex is a minimal free resolution of $I(T)$. 
In general, if $T$ is the full $A_1 \times \cdots \times A_n$, the free cellular
complex is always a minimal free resolution of $I(T)$. Finding out for
which $T$ this free cellular
complex {\it is} a resolution, Theorem \ref{thm:cell-linres},
is an aspect of the main result of our article. 
\end{exam}

The generators of $\Cal_p(T)$ for $p \geq 0$ have degrees $p+n$,
so $\Ctcal_{\geq 0}(T)$ is a {\it linear}
free cellular complex.
The free cellular complex $\Ctalg(T)$ is $\NN^A$-graded. A family of subsets
$\Rb \susst \Ab$ gives a restriction
$T_{|R} = T \cap (R_1 \times \cdots \times R_n)$. 
It also gives a squarefree degree $\deg(\Rb) \in \NN^A$.
The free cellular complex  $\Ctalg(T)_{|\deg(\Rb)}$ in this specific degree
identifies as the augmented chain complex $\Ctbu(T_{|\Rb})$
(see \eqref{eq:minfree-Ct}).

\begin{fact} \label{fact:cell-hom} Let $\Rb \sus^* \Ab$. 
  The homology of the augmented free cellular complex $\Ctalg(T)$
  in degree $\deg(\Rb) \in \NN^A$ is the reduced homology
  \[ H_i({\Ctalg}(T))_{\deg(\Rb)} = {H}_i(\Ctbu(T_{|\Rb}),
    \, \, i \geq -1.\]
\end{fact}

\noindent This is standard and rather immediate from the above,
see \cite[Section 4]{Mi-St}.

\begin{prop} \label{pro:cellcx-homology}
  For $p \geq 0$, the homology group $ H_p({\Ctalg}(T))$
 can only be nonzero in degrees $\geq p+n+2$.
\end{prop}

\begin{proof}
  In short, this is a consequence of $C(T)$ being a saturated cell complex.
 All generators of $\Cal_p(T)$ are in degree $p+n$.
  Thus $\ker d_{p}$ is generated in degrees $\geq p+n+1$.
Let $u$ be a generator of $\ker d_p$ of degree $p+n+1$.
The free cellular complex is $\NN^A$-graded, and so we may assume $u$ is
  homogeneous for this grading.
  This degree  is also squarefree
  as it is a minimal generator. Its degree is then  $\deg(\Rb)$
  for some $\Rb \susst \Ab$.
  Consider the free cellular complex $\Ctalg(\Rb)$. We have subcomplexes
  \[ \Ctalg(T_{|\Rb}) \sus \Ctalg(\Rb), \quad \Ctalg(T_{|\Rb}) \sus \Ctalg(T).\]
The latter is an isomorphism in degree $\deg(\Rb)$. 
  Thus if $\ker d_{p}$ has a generator $u$ in degree $\deg(\Rb)$, due to
  $\ker d_p^{\Ctalg(T_{|R})} \hookrightarrow \ker d_p^{\Ctalg(\Rb)}$ we get
  a generator of the latter kernel (in degree $\deg(\Rb)$),
  of total degree $p+1+n$.
  But $\ker d_p^{\Ctalg(\Rb)} = \im d_{p+1}^{\Ctalg(\Rb)}$
and $\Cal_{p+1}(\Rb)$ is the free rank one $S$-module
  $S \! \overset{}{\cdot} \! {\Pi(\Rb)}$,  with the generator $\Pi(\Rb)$ of degree $p+1+n$.
  This maps to a sum involving every face of $\Pi(\Rb)$ of dimension
  $p$. But due to the commutativity of the inclusion maps and differentials
  of 
$\Ctalg(T_{|\Rb}) \hookrightarrow \Ctalg(\Rb)$, 
 the degree $p+n$ part $\Cal_p(T_{|\Rb})$ must have
 each of the generators corresponding to a face of $\Pi(\Rb)$
 of dimension $p$. 
  In particular all vertices of $\Pi(\Rb)$ must be in $T$.  Since
  $\Cgeo(T)$ is saturated, 
  the full cell $\Pi(\Rb)$ is in $\Cgeo(T)$.
  Then the generator 
  $u$ of $\ker d_p$ is the image of a cell in $\Cal_{p+1}(T)$. 
\end{proof}

\subsection{The criterion on $T$ to correspond to a polarization}

We now get a completely explicit description of the resolution of $I(T)$,
when its resolution is linear.

\begin{theo} \label{thm:cell-linres}
  Let $T \sus A_1 \times \cdots \times A_n$. If $I(T)$ has
  linear resolution, the free cellular complex  $\Ctcal_{\geq 0}(T)$ is its minimal
  free resolution. (And  $\Ctcal_{\geq 0}(T)$ is a minimal free resolution
  precisely in this case.)
\end{theo}

\begin{proof}
  Let $\Falg$ be a minimal free resolution of $I(T)$. As $I(T)$ is
  $\NN^A$-graded, the same holds for $\Falg$. Since $\Falg$ is a resolution,
  and $\Calg(T)$ consists of free modules, there is a morphism
  $\phi : \Calg(T) \pil \Falg$ lifting the identity map on $I(T)$.
  Both $\Fal_0$ and $\Cal_0(T)$ have free bases mapping bijectively to the
  minimal generators of $I(T)$, the monomials $m(\bfa)$ for $\bfa \in T$.
  So there is an isomorphism $\Fal_0 \iso \Cal_0(T)$.

  Assume $\phi_i$ is an isomorphism for $i \leq p$. Consider the diagram:
 \[ \xymatrix{ \Cal_p(T) \ar[d]^{\phi_p} &
      \Cal_{p+1}(T) \ar[l]^{d_{p+1}^{\Cal}} \ar[d]^{\phi_{p+1}} \\
      \Fal_p & \Fal_{p+1}  \ar[l]^{d_{p+1}^{\Fal}}.} \]
  Since $\phi_p$ is injective, and $d_{p+1}^{\Cal}$ is injective on
  the vector space of generators, which has degree $p+n+1$, 
  the right vertical map $\phi_{p+1}$ is also injective on
  the vector space of generators. 
    The map $\phi_{p+1}$ is then injective since
    both its domain an codomain are free with generators in degree $p+1+n$.

  Let $Q_{p+1}$ be the cokernel of the right vertical map. It is an $S$-module
  generated in degree $p+n+1$. By the long exact cohomology sequence of
  the short exact sequences of complexes
  \[ 0 \pil \Cal_{\leq p+1}(T) \pil \Fal_{\leq p+1} \pil Q_{p+1}[-p-1] \pil 0 \]
(where $[-p-1]$ denotes that $Q_{p+1}$ is in homological degree $p+1$)
  we see that $Q_{p+1}$ maps surjectively to the homology $H_p(\Calg(T))$.
  So the latter module is generated in degree $p+n+1$. But it has no
  homology in this degree, by Proposition \ref{pro:cellcx-homology}.
 Hence this homology module vanishes.
  Then
  \[ H_{p+1}(\Fal_{\leq p+1}) \pil H_{p+1}(Q_{p+1}[-p-1]) = Q_{p+1} \]
  is surjective, where the left module is in fact the kernel
  $\ker d_{p+1}^{\Fal}$.
  Since this kernel is generated in degree $p+n+2$, this again implies
  that $Q_{p+1} = 0$. Whence $\phi_{p+1}$ is an isomorphism. 
\end{proof}

\begin{rema}
  The above result has a more general framework. 
  K.Yanagawa \cite{Ya99} shows that
  when a squarefree monomial ideal $I_{\Delta}$ has a linear resolution, it
  may be constructed directly from the canonical module $\omega_{\Delta^\vee}$
  of the (Cohen-Macaulay) Stanley-Reisner ring of the Alexander dual
  simplicial complex $\Delta^\vee$.
  In fact given a squarefree module $M$, there is a linear complex
  $F_\bullet(M)$ with terms given by $F_i(M) = \oplus_{|\bfr| = i}
  M_{\ben- \bfr} \te S(-\bfr)$ (where $\bfr$ has entries $0$ or $1$ and $\ben = (1,1,\ldots, 1)$). 
  The linear resolution of $I_{\Delta}$ is
  $F_\bullet(\omega_{k[\Delta^\vee]})$.
  
  When $\Delta$ is a full-dimensional triangulated ball
  on a simplicial polytope $P$ (as is the case here),
  by \cite[Thm.5.7.2]{BH93} the canonical module is given as
  $\omega_{k[\Delta^\vee]} = I_\Sigma/I_\Delta$ where $\Sigma$
  is the boundary of $\Delta$, 
  and the latter identifies
  as $I_{\Delta^c}/I_P$ where $\Delta^c$ is the complement simplicial complex
  on $P$. One works out that $F_{\bullet}(I_{\Delta^c}/I_P)$
  identifies as the dual of
  a free cellular complex supported on
  the relative simplicial complex $(P, \Delta^c)$.
  Letting $P^\vee$ be the dual polytope (in our case this is the
  product polytope  $\Pi(\Ab)$), the linear resolution identifies again as
  a subcomplex of the free cellular complex associated to $P^\vee$.
\end{rema}

The following is the general homological criterion on a subset
$T \sus A_1 \times \cdots \times A_n$ for when
$I(T)$ has a linear resolution.

\begin{coro} \label{cor:cell-homcrit}
  $I(T)$ has a linear resolution if and only if all homologies of
  all restrictions vanish: 
  \begin{equation} \label{eq:cellcx-homvan} {H}_i(\Ctbu(T_{|\Rb})) = 0
    \end{equation}
 for all $i \geq 0$ and $\Rb \sus^* \Ab$. 
\end{coro}

\begin{proof} Use Fact \ref{fact:cell-hom}.
\end{proof}

\begin{defi}
  A subset $T \sus A_1 \times \cdots \times A_n$ is {\it acyclic}
  if the vanishing \eqref{eq:cellcx-homvan} holds for all $i \geq 0$
  and $\Rb \susst \Ab$. We say the cell complex $C(T)$
  is {\it restriction-acyclic} or {\it \rac}. 
\end{defi}

\begin{rema} Note that  $C(T)$ may be acyclic (all homologies vanish),
but not be restriction-acyclic. Thus being r-acyclic is 
  not an intrinsic property
  of $C(T)$, but depends on how it is embedded in $\Pi(\Ab)$, see
  Remark \ref{rem:cube-acyclic}. Section \ref{sec:hypercube} has
  many examples of r-acyclic cell complexes.
\end{rema}

\begin{coro} $J(T)$ is a polarization of some Artin monomial ideal
  in $k[x_1, \ldots, x_n]$ containing $\ic = (x_1^{\alpha_1}, \ldots,
  x_n^{\alpha_n})$ iff $T$ is acyclic.
\end{coro}

\subsection{Inductive construction of sets $T$}

When $T$ is acyclic, the constructability established in
Section \ref{sec:polball} suggests $T$ may be built up inductively.
Here we detail this.

The subset $T \sus A_1 \times \cdots \times A_n$ gives
the saturated cell sub-complex $X = \Cgeo(T)$ of $\Pi(\Ab)$.
Vertices of the latter are $\bfa = (a_1,a_2, \ldots, a_n)$ where
$a_i \in A_i$. 
For $a \in A_1$ let:
\begin{itemize}
\item $T_a = \{\bfa \in T \, | \, a_1 = a \}$,
  \item $T_{-a} = \{ \bfa \in T \, | \, a_1 \neq a\}$,
  \item $T_{/a} = \{ \bfa \in T \, | \, a_1 \neq a
    \text{ and both } \bfa \text{ and } (a, a_2, \ldots, a_n) \in T \}$,
  \item $T^a = T_a \cup T_{/a}$.
  \end{itemize}

  Furthermore we have the associated cell sub-complexes of $X$:
\begin{itemize}
\item $X_a = \Cgeo(T_a)$,
\item $X_{-a} = \Cgeo(T_{-a})$,
\item $X_{/a} = \Cgeo(T_{/a})$,
\item $X^a = \Cgeo(T^a)$.
\end{itemize}

Note that $X_{/a}$ consists of the cells in $X_{-a}$ which may be extended to a
cell with a vertex in $T_a$.
Also, $X^a$ consists of the sub-cells 
of those cells of $X$ which have some vertex in $T_a$.
Note further that
\begin{equation} \label{eq:cellcx-X} X = X^a \cup X_{-a},
  \quad X_{/a} = X^a \cap X_{-a}.
  \end{equation}

\begin{lemm}
  $X_a$ and $X^a$ are homotopy equivalent.
\end{lemm}

\begin{proof}
  There is an inclusion $X_a \sus X^a$. For any set
  $B$ and $a \in B$, there is a retract of the simplex
  $\Delta(B)$ on $\Delta(\{a \})$ given by 
  \[ \sum_{b \in B}  \alpha_b \cdot b \mapsto
    (\sum_{b \in B}  \alpha_b \cdot b)(1-t) + ta. \]
  By letting $B = A_1$ we get a retract of
  $\Pi(\Ab)$ on $\Delta(\{a \}) \times \Delta(A_2) \times \cdots
  \times \Delta(A_n)$, and this restricts to a retract of
  $X^a$ on $X_a$. 
\end{proof}

The following shows how we can inductively build up acyclic
sets $T$. When $\Rb \susst \Ab$,
  for simplicity denote $X_{|\Pi(\Rb)}$ by $X_{|R}$. 

\begin{prop} \label{pro:cell-Xac}
 A saturated cell sub-complex $X = C(T)$ is r-acyclic iff
\begin{itemize}
\item  $X_a, X_{-a}$ and $X_{/a}$ are \rac,
\item For every $\Rb \susst \Ab$,
  whenever $X_{a|R}$ and $X_{-a|R}$ are non-empty, then
  $X_{/a|R}$ is non-empty.
\end{itemize}
In this situation,  $X^a$ is also \rac. 
\end{prop}

\begin{proof}
By \eqref{eq:cellcx-X}
we have a Mayer-Vietoris sequence (using reduced homology):
  \begin{equation*} \tH_{i+1}(X_{|R}) \pil \tH_i(X_{/a |R}) \pil
                   \tH_i(X^a_{|R}) \oplus \tH_i(X_{-a |R}) \pil
 \tH_i(X_{|R}) \pil \tH_{i-1}(X_{/a | R}) \pil \cdots .
\end{equation*}
Assume $X$ is \rac.
Since $X_a$ and $X_{-a}$ are restrictions,
they are \rac. Further $X^a$ is homotopy equivalent to $X_a$, and
similarly $X^a_{|R}$ is homotopy equivalent to $X_{a | R}$. 
Thus
\[ \tH_i(X^a_{|R}), \quad \tH_i(X_{-a |R}), \quad \tH_i(X_{|R})\]
all vanish for $i \geq 0$. Thus $X_{/a \,|R}$ is \rac.
Also, when  $X_{a|R}$ and $X_{-a|R}$ are non-empty, the sequence above,
when $i = -1$, gives that $\tH_{-1}(X_{/a|R})$ vanishes, and so
$X_{/a|R}$ is non-empty.

Conversely assume  $X_a, X_{-a}$ and $X_{/a}$ are all \rac.
The sequence above gives that all $\tH_i(X_{|R})$ vanish for $i \geq 0$,
except possibly
in the following case: $i = 0$ and  $X_{/a | R}$ is empty and
$X_{|R}$ is not empty.
But in this case at least one of $X_{-a |R}$ or $X_{a | R}$ is empty,
which implies the vanishing of $\tH_0(X_{|R})$.
\end{proof}

\newcommand{\nrel}{\not \sim}
\newcommand{\psiz}{\text{psize}\,}
\newcommand{\bfp}{{\mathbf p}}

\section{Minimal irreducible decomposition and
  maximal cells}
\label{sec:mid}

Given an Artin monomial ideal $\ijfr$.
The question we address here is what does $T$ and $C(T)$ look like
for polarizations $J(T)$ of $\ijfr$? How can we make various polarizations
of $\ijfr$ by working with the geometry of the cell complex $C(T)$?

We show there is a one-to-one correspondence between the maximal
faces of $C(T)$ and the components of the
minimal irreducible decompositions of the
Artin monomial ideal $\ijfr$. This is what guides us to
find the possible $T$'s which gives polarizations $J(T)$ of $\ijfr$,
see Section \ref{sec:hypercube}.

\subsection{Maximal cells in $C(T)$}
We first need some essential structural properties of the maximal cells
in $C(T)$ for an acyclic $T \sus A_1 \times \cdots \times A_n$. 

\begin{lemm} \label{lem:mid-lka}
  Let $\Pi(R_\bullet)$ be a maximal cell in $C(T)$. Suppose $R_1$ has at
  least two elements, and that $a \in R_1$ is one of them.
  Let $\overline{R}_\bullet$ be $R_\bullet$ but with $R_1$ replaced by
  $R_1 \backslash \{a\}$.
  Then $\Pi(\overline{R}_\bullet)$ is a maximal cell in $C(T_{/a})$.
  In particular, since $T_{/a} \sus T_{-a}$,
  there is a maximal cell in $C(T_{-a})$ containing
  $\Pi(\overline{R}_\bullet)$. 
\end{lemm}

\begin{proof}
  Let  $\Pi(R^\prime_\bullet)$ be a maximal cell in $C(T_{/a})$ with
  $R^\prime_i \supseteq {R}_i$ for $i \geq 2$,
  and $R^\prime_1 \supseteq R_1 \backslash \{a \}$. By definition of
  $T_{/a}$, adding $a$ to $R_1^\prime$ and keeping the
  rest of the $R^\prime_i$ unchanged, gives a new $R_\bullet^{\prime \prime}$
  and cell $\Pi(R_\bullet^{\prime \prime})$. This cell must be maximal in $C(T)$
  as $\Pi(R^\prime_\bullet)$ is maximal in $C(T_{/a})$. But then
  $\Pi(R_\bullet^{\prime \prime})$ must be $\Pi(R_\bullet)$, and so we must have each
  $R_i^\prime = \overline{R}_i$.
\end{proof}

\begin{prop} \label{pro:minirr-max}
  Let $\Pi(R_\bullet)$ and $\Pi(S_\bullet)$ be distinct maximal cells in
  an r-acyclic cell complex $C(T)$. Then there are indices $i$ and $j$
  such that we have strict inclusions
  \begin{equation} \label{eq:mid-RS}
    R_i \subset S_i, \quad S_j \subset R_j.
    \end{equation}
\end{prop}

\begin{proof}
  Since $C(T_{-a})$ is r-acyclic when $C(T)$ is, 
  we may by successive restrictions for $a \in A_i \backslash (R_i \cup S_i)$,
  assume that $R_i \cup S_i = A_i$ for each $i$.

  We want to show we cannot have the negation of \eqref{eq:mid-RS}.
  So assume the negation. Write $R_i \nrel S_i$ if none of them
  is included in the other. We may then assume that there is $0 \leq p \leq n$
  such that $R_i \supseteq S_i$ for $i = 1, \ldots, p$ and $R_i \nrel S_i$
  for $i = p+1, \ldots, n$.

  Our proof is by induction on $n$ and the size of the $A_i$.
  As an induction base
  we note two cases: 
  i) If $n = 1$ they cannot both be maximal. ii)
  We also cannot have all $|R_i| = 1$,
  since then $\Pi(R_\bullet)$ is a point and
  would not be connected to the rest, contradicting
  that $C(T)$ is \rac.

\noindent{\bf Case $p \geq 1$}: So $R_1 \supseteq S_1$. 
Assume $R_1 = S_1$. Let $a \in S_1$. If the cardinalities $|R_1| = |S_1|$
is $\geq 2$, we may by Lemma \ref{lem:mid-lka} get 
  smaller maximal sets in $C(T_{/a})$. Since $T_{/a}$ is acyclic,
  this is against induction. If $R_1 = S_1$ is of
  cardinality $1$, we may simply drop the first factor and reduce $n$. 

  So assume $R_1$ is strictly bigger than $S_1$.
 Let  $a \in R_1 \backslash S_1$ (so $|R_1| \geq 2$). We then get
  a maximal $\Pi(R^\prime_\bullet)$ in $C(T_{-a})$ with $R^\prime_i \supseteq R_i$
  for each $i \geq 2$. Also $R^\prime_1 \supseteq S_1$.
  Again since $T_{-a}$ is acyclic,
  we get a contradiction by induction, since we have no strict
  $R^\prime_i \subset S_i$.

  \noindent {\bf Case $p = 0$}: So we have $R_i \nrel S_i$ for each $i$.
  If some $R_i \cap S_i $ is non-empty (then both $|R_i|, |S_i| \geq 2$),
  we may use Lemma \ref{lem:mid-lka}, and get a contradiction by induction.
  So we may assume every $R_i \cap S_i$ is empty.

  At least one of the $R_i$ has more than one element, say $R_1$.
  Let $a \in R_1$. By restricting away from $a$, we get
  a maximal $\Pi(R^\prime_\bullet)$ in $C(T_{-a})$ with $R^\prime_i \supseteq R_i$
  for each $i \geq 2$, and so none $R^\prime_i \subseteq S_i$ for $i \geq 2$.
  We cannot have $R_1^\prime \sus S_1$, since $R_1$ and $S_1$ had
  empty intersection, and $|R_1| \geq 2$.
  By induction we get a contradiction.
\end{proof}

\subsection{Regular linear subspaces}
Let $V$ be a finite dimensional $k$-vector space and $\ell \sus V$
a linear subspace.

\begin{defi}
  For a polynomial ring $S = k[V]$ let $I \sus S$ be an ideal.
  If $\ell$ has a basis which is a regular sequence for $S/I$, we call
  $\ell$ a {\it regular linear space} for $S/I$ (or for $I$).
\end{defi}

\begin{lemm} \label{lem:reg-tor} Let $\ell$ be a linear subspace of $V$
  and $L = (\ell) \sus S$ the ideal generated by $\ell$.
  Then $\ell$ is a regular linear space for $S/I$
  iff ${\rm Tor}_i(S/I,S/L) = 0$ for
  all $i > 0$.
\end{lemm}

\begin{proof}
  This is by \cite{Eis} Corollary 17.5 and Theorem 17.6.
  \end{proof}

\begin{lemm} \label{lem:mid-IJL}
Let $I,J \sus S$ be ideals. Let $\ell \sus V$ be a regular linear space for
the ideals $I,J$ and $I+J$. Then $\ell$ is a regular linear space for $I\cap J$
and in the quotient ring $S/L$ (where $L = (\ell)$) we have:
\begin{equation} \label{eq:mid-IJL}
\frac{(I \cap J) + L}{L} = \left ( \frac{I+L}{L} \right )
\cap \left ( \frac{J+L}{L} \right ). \end{equation}
\end{lemm}

\begin{proof}
We have a short exact sequence
\[ 0 \pil S/(I \cap J) \pil S/I \oplus S/J \pil S/(I + J) \pil 0. \]
By tensoring with $S/L$, using the long exact Tor-sequence, and
Lemma \ref{lem:reg-tor}, we get that $\ell$ is
a regular linear space for $I \cap J$. We also get an exact sequence:
\[ 0 \pil \frac{S/L}{((I \cap J) + L)/L} \pil \frac{S/L}{(I+L)/L}
\oplus \frac{S/L}{(J+L)/L} \pil \frac{S/L}{(I + J+L)/L}  \pil 0, \]
which gives Equation \eqref{eq:mid-IJL}.
\end{proof}

The following rather trivial result will also be used.

\begin{lemma} \label{lem:mid-Ka}
Let $K, K_u, u \in U$ be graded ideals in $S = k[V]$, and suppose
$K = \cap_{u \in U} K_u$.  Let $a$ be a new variable.
These ideals generate ideals $(K), (K_u)$ in $k[V \oplus (a)]$.

\begin{itemize}
\item[a.] We have $(K) = \cap_{u \in U} (K_u)$ in $k[V \oplus (a)]$.
\item[b.] In the quotient ring $k[V] \iso k[V \oplus (a)]/(a)$ we have
  \begin{equation} \label{eq:mid-Ka} \frac{(\cap_u (K_u))  + (a)}{(a)}
    = \frac{(K) + (a)}{(a)} =
  \bigcap_{u \in U}\left ( \frac{(K_u) + (a)}{(a)} \right ).
  \end{equation}
\end{itemize}

\end{lemma}

\begin{proof}
  Note that all ideals above are $\ZZ^2$-graded, with  $(d, d_a)$ denoting
  the degree of a homogeneous polynomial where $d$ is the degree of
  terms in $k[V]$, and  $d_a$ is the power of $a$
  occurring.

  \noindent{a.} It is clear that $(K) \sus \cap_u (K_u)$.
  Let $p = p_0 + ap_1 + a^2p_2 \cdots$ be in the latter intersection.
  Due to the $\ZZ^2$-grading each $p_i$ is in each $(K_u)$ and so in $K_u$.
  Thus each $p_i$ is in $K$ and so $p$ is in $(K)$.

  \noindent {b.} Let $p \in k[V]$. Write $\overline{p}$ when
  considered as element in $k[V \oplus \langle a \rangle]/(a)$.
  If $\overline{p}$ is in the intersection on the right of
  Equation \eqref{eq:mid-Ka}, 
  then for each $u \in U$, $p + aq_u$ is in $(K_u)$ for some $q_u$.
  Due to $\ZZ^2$-grading,
  $p \in K_u$ for each $u$, and so $p$ is in $K$ and so
  $\overline{p}$ in the left side of Equation \eqref{eq:mid-Ka}.
\end{proof}

\subsection{Minimal irreducible decompositions}
Let now $T \sus A_1 \times \cdots \times A_n$ be an acyclic subset.
It gives the ideal
$I(T)$ with Alexander dual ideal $J(T)$, a polarization of an Artin
monomial ideal.

Writing $A_i = \{a^i_1, a^i_2, \ldots \}$, 
for each $i$ the variable differences $a^i_k - a^i_l$
generate a linear subspace $\ell_i$.
Let $\ell = \ell_1 + \ell_2 + \cdots
+ \ell_n$. It generates an ideal $L = (\ell) \sus S$, and
\[ k[A_1 \cup A_2 \cup \cdots \cup A_n]/L \iso k[x_1,x_2, \ldots, x_n].\]
Let $\Rb \sus^* \Ab$. It gives $U = R_1 \times \cdots \times R_n$.
We write $I(\Rb), J(\Rb)$ for $I(U), J(U)$. 
Then
\[ J(U) = J(\Rb) = (m(R_1), \ldots, m(R_n)). \]
Thus
\begin{equation} \label{eq:mid-JTL}
  \frac{J(U) + L}{L}  \sus \frac{k[A_1 \cup \cdots \cup A_n]}{L}
  \end{equation}
identifies as
\begin{equation} \label{eq:mid-xr}
  (x_1^{r_1}, x_2^{r_2}, \cdots, x_n^{r_n}) \sus k[x_1, \ldots, x_n],
  \quad r_i = |R_i|. \end{equation}

\begin{theorem} \label{thm:mid-mid}
   Let $T \sus A_1 \times \cdots \times A_n$ be an acyclic subset, 
  which is equivalent to $J(T)$ being a polarization of an Artin monomial ideal $\ijfr$ (by Corollary \ref{cor:cell-homcrit} and Proposition \ref{pro:pol-cm}).
Let $\Pi(R^k_\bullet), k = 1, \ldots, m$ be the maximal cells in $C(T)$.
  Denote $\ijfr_k = (x_1^{r^k_1}, x_2^{r^k_2}, \ldots, x_n^{r^k_n})$,
  the irreducible monomial ideal where $r^k_i$ is the cardinality
  of $R^k_i$.
  \begin{itemize}
  \item[a.] \[ J(T) = J(R^1_\bullet) \cap J(R^2_\bullet) \cap \cdots \cap
      J(R^m_\bullet). \]
  \item[b.] The minimal irreducible decomposition of $\ijfr$ is
    \begin{equation} \label{eq:mid-irr}
      \ijfr = \ijfr_1 \cap \ijfr_2 \cap \cdots \cap \ijfr_m.\end{equation}
  \end{itemize}
  So the components of the minimal irreducible decomposition of $\ijfr$
  are in bijection to the maximal cells of $C(T)$.
\end{theorem}

\begin{proof}
  Part a. We have
  \[ I(T) = I(R^1_\bullet) + I(R^2_\bullet) + \cdots + I(R^m_\bullet). \]
The Alexander dual is then
\[ J(T) = J(R^1_\bullet) \cap J(R^2_\bullet) \cap \cdots \cap J(R^m_\bullet). \]

\medskip
Part b. We must
show:
\begin{equation} \label{eq:mid-JL}
 \ijfr =  \frac{J(T) + L}{L} =
  \frac{J(R^1_\bullet) + L}{L} \cap \frac{J(R^2_\bullet) + L}{L} \cap
\cdots \cap \frac{J(R^m_\bullet) + L}{L}
\end{equation}  
which shows \eqref{eq:mid-irr}. In addition we must show
the minimality of the decomposition. 
\Ign{
We now write
\[J_a, J_{-a}, J_{/a}, J^a, J(\Rb)
  \text{ for the Alexander duals of }
  I_a, I_{-a}, I_{/a}, I^a, I(\Rb). \]
}
We have
\begin{equation} \label{eq:mid-IT1}
  I(T)  = I(T^a) + I(T_{-a}), \quad I(T_{/a}) = I(T^a) \cap I(T_{-a}).
\end{equation}
The first is clear. To see the second, since $I(T^a) = I(T_a) + I(T_{/a})$,
the intersection
\begin{equation} \label{eq:mid-IT2}
    I(T^a) \cap I(T_{-a}) =
    I(T_a) \cap I(T_{-a}) + I(T_{/a}) \cap I(T_{-a}).
  \end{equation}
  The last term is $I(T_{/a})$. As for the middle term, by the
 last arguments of Proposition
  \ref{pro:lin-ItoIR} it has $(n+1)$-linear resolution. A monomial generator
  of this intersection must have the form $a a_1 a_2 \cdots a_n$ where
  $a_i \in A_i$. Then $a_1 \cdots a_n \in I(T_{-a})$ and $aa_2 \cdots a_n
  \in I(T_a)$. So $(a_1, a_2, \ldots, a_n) \in T_{/a}$. Thus the left
  side of \eqref{eq:mid-IT2} is $I(T_{/a})$. 

On the Alexander dual side, \eqref{eq:mid-IT1} becomes:  
\[ J(T) = J(T^a) \cap J(T_{-a}), \quad J(T_{/a}) = J(T^a) + J(T_{-a}). \]
If $C(T)$ has only one maximal cell we are done by Equations
\eqref{eq:mid-JTL} and \eqref{eq:mid-xr}.
Suppose now $C(T)$ has at least two maximal cells. Let $a$ be an element
of some $A_i$, say $a \in A_1$ which is on at least one maximal cell
but not all.
Then $J(T^a)$ and $J(T_{-a})$ are ideals properly containing $J$, and
both $T^a, T_{-a}$ and $T_{/a}$, are acyclic subsets of $A_1 \times
\cdots \times A_n$, by Proposition \ref{pro:cell-Xac}.
Thus by Lemma \ref{lem:mid-IJL}
\[ \frac{J(T)+L}{L} = \left ( \frac{J(T^a) + L}{L}\right ) \cap
\left ( \frac{J(T_{-a}) + L}{L} \right ). \]
Both  $T^a$ and $T_{-a}$ are properly included in $T$.
Furthermore a maximal cell $\Pi(R^j_\bullet)$ of $C(T)$ is either in 
$C(T^a)$ (if $a \in R^j_1$) or in $C(T_{-a})$ (if $a \not \in R^j_1$).
Let $M^a$ denote the index set of $j$ in the former,
and $M_{-a}$ the $j$ in the latter. Then by induction and 
Lemma \ref{lem:mid-Ka}:
\[ \frac{J(T)+L}{L} = \cap_{j \in M^a} 
\left ( \frac{J(R^j_\bullet) + L}{L}\right ) \bigcap
\cap_{j \in M_{-a}} \left ( \frac{J(R^j_\bullet) + L}{L} \right ). \]
This show Equation \eqref{eq:mid-JL}.

To show minimality, suppose there was some $p \in [m]$ such
that 
\begin{equation} \label{eq:mid-jji}
 \cap_{i \in [m] \backslash \{p\}} \ijfr_i \sus \ijfr_p.
\end{equation}
Let $\Mon(n)$ be the monomials in the polynomial ring $k[x_1, \ldots, x_n]$,
and if $K \sus k[x_1, \ldots, x_n]$
is a monomial ideal, write $\Mon(n) = N(K) \sqcup \Mon(K)$ as a disjoint 
decomposition, where $\Mon(K)$ are the monomials in $K$ and $N(K)$ are
the normal monomials, those not in $K$.
Expression \eqref{eq:mid-jji} corresponds to 
\[ \cup_{i \in [m] \backslash \{p\}} N(\ijfr_i) \supseteq N(\ijfr_p). \]
But writing $\ijfr_p = (x_1^{r^p_1}, \ldots, x_n^{r^p_n})$,
the maximal element (for the division order in $\Mon(n)$) of $N(\ijfr_p)$ is
the monomial $x_1^{r^p_1-1}x_2^{r^p_2 -1} \cdots x_n^{r^p_n-1}$. If this
is in $N(\ijfr_i)$, then each $r^i_k \geq r^p_k, k = 1, \ldots, n$. 
But this contradicts Proposition \ref{pro:minirr-max}.
Thus we cannot have the situation of \eqref{eq:mid-jji}. 
\end{proof}

For $S \sus A = A_1 \cup \cdots \cup A_n$, write $S_i = S \cap A_i$.
 We say $S$ {\it covers} $T \sus A_1 \times \cdots \times A_n$ if
 for every maximal cell $\Pi(\Rb)$ in $C(T)$, there is some $i$ such
 that $R_i \sus S_i$. Note in particular that each $A_i \sus A$ covers
 $T$. 

 \begin{coro} \label{cor:mid-cover} The polarization $J(T)$ is generated by the monomials
   $m(S_\bullet)$ such that $S$ covers $T$.
 \end{coro}

 \begin{proof} It is clear that $S$ covers $T$ iff $m(S_\bullet)
   \in J(\Rb)$ for each maximal cell $\Pi(\Rb)$ of $C(T)$.
 \end{proof}

 For an {\it irreducible} Artin monomial ideal $\ijfr = (x_1^{r_1},x_2^{r_2},
 \ldots, x_n^{r_n})$ define its {\it power size} to be
 \[ \psiz (\ijfr) = \sum_{i = 1}^n (r_i - 1). \]
 (It is the number of variable powers $x_i^{s_i}$ in the
 quotient ring $k[x_1, \ldots, x_n]/\ijfr$.)
 The following will also be useful.

\begin{lemm} \label{lem:mid-dim}
  For any maximal cells $\Pi(\Rb^i)$ and $\Pi(\Rb^j)$
  in $C(T)$, the intersection $\Pi(\Rb^i) \cap \Pi(\Rb^j)$ has
  dimension $\leq \psiz (\ijfr^i + \ijfr^j)$. 
\end{lemm}

\begin{proof}
  Let $P_k = R^i_k \cap R^j_k$. Then the intersection
  $\Pi(\Rb^i) \cap \Pi(\Rb^j) = \Pi(P_\bullet)$. Thus
  $I(\Rb^i) \cap I(\Rb^j) \supseteq I(P_\bullet)$.
  Taking Alexander duals this gives
  \[ J(\Rb^i) + J(\Rb^j) \sus J(P_\bullet). \]
  Then for the associated irreducible Artin ideals in $S/L$ we have:
  \[ \ijfr^i + \ijfr^j \sus \ijfr^\prime, \quad \ijfr^\prime =
    \frac{J(P_\bullet) + L}{L}. \]
  Whence
  \[ \psiz (\ijfr^i + \ijfr^j)  \geq \psiz(\ijfr^\prime) = \dim \Pi(P_\bullet). \]
\end{proof}

  \subsection{Standard polarization}
Fix a total order on each $A_i$. So $A_i = \{a^i_1 < a^i_2 < \cdots
< a^i_{\alpha_i} \}$. For each $0 \leq r_i \leq \alpha_i$ define
$m^{ r_i}(A_i) = \prod_{j = 1}^{r_i} a^i_j$.
Each $\bfr = (r_1,r_2, \ldots, r_n)  \in \NN^n$ gives a monomial
$m(\bfr)  = x_1^{r_1} \cdots x_n^{r_n}$. The {\it standard polarization} of
the monomial $m({\bfr})$ is
\[ m^{\bfr}(A_\bullet) = \prod_{i = 1}^n m^{r_i}(A_i). \]
If $\ijfr = (x^{\bfr^1}, \ldots, x^{\bfr^N})$, its {\it standard polarization} is
\[ J = (m^{\bfr^1}(A_\bullet), \ldots, m^{\bfr^N}(A_\bullet)). \]

\begin{exam}
  Consider $\ijfr = (x_1^3x_2, x_2^4x_3^2, x_1x_3^3)$. Rename
  $A := A_1, B := A_2$ and $C := A_3$, and write $b_i$ instead of $a^2_i$
  for the elements of $B$, and similarly for $C$. The standard polarization
  of $\ijfr$ is then:
  \[ J^{st} = (a_1a_2a_3b_1, b_1b_2b_3b_4c_1c_2, a_1c_1c_2c_3). \]
\end{exam}

It is well-known that the standard polarization of any monomial
ideal is a polarization, \cite[Prop.1.6.2]{HeHi}

\begin{prop}
  Let $\ijfr$ be an Artin monomial ideal, with
  $\ijfr = \cap_{k = 1}^m \ijfr_k$ its minimal irreducible decomposition,
  where $\ijfr_k = (x_1^{r^k_1}, x_2^{r^k_2}, \ldots, x_n^{r^k_n})$. 
  Then the standard polarization of $\ijfr$ is
  \[ J^{st} = \cap_{k = 1}^m J^{st}_k. \]
  where $J^{st}_k = (m^{r^k_1}(A_1), \ldots, m^{r^k_n}(A_n))$.
\end{prop}

\begin{proof}
  That $J^{st} \sus J^{st}_k$ is immediate to see. Let then $m$
  be a monomial in the intersection. Each $J_k^{st}$ has a generator
  $m^{r^k_{i_k}}(A_{i_k})$ which divides $m$. Then $m$ is
  divisible by the least common multiple of these monomials. This 
  least common multiple  is $m^{\bfr}(A_\bullet)$ for a suitable $\bfr$.
  Furthermore $x^{\bfr}$ is the least common multiple of the
  $x_{i_k}^{r_{i_k}}$. So this is in the intersection $\cap_{k = 1}^m \ijfr_k$,
  and so $x^{\bfr}$ is in $\ijfr$. Write $x^{\bfr} = x^{\bfr^\prime} x^{\bfp}$
  where $x^{\bfp}$ is a minimal generator for $\ijfr$. Then
  $m^{\bfp}(A_\bullet) \in J^{st}$, and $m^{\bfr}(A_\bullet)$ is a multiple
  of $m^{\bfp}(A_\bullet)$ and so $m$ is in $J^{st}$.
\end{proof}

\begin{coro} \label{coro:mid-cellstd}
  Let $I^{st}$ be the Alexander dual of the standard polarization
  $J^{st}$, so $I^{st} = I(T^{st})$ for a subset $T^{st} \sus A_1 \times
  \cdots \times A_n$.
  The maximal cells of $C(T^{st})$ are the $\Pi(R^k_\bullet)$ where
  \[ R^k_\bullet = \{ (a^1_{i_1}, \ldots, a^n_{i_n}) \, | \,
    1 \leq i_j \leq r^k_j \text{ for } j = 1, \ldots, n \}. \]
  So this cell is
  \[ \Delta(R^k_1) \times \cdots \times \Delta(R^k_m) \]
  where $R^k_i$ is the initial segment of $A_i$ with $r^k_i$ terms.
  (Note that $(a^1_1, a^2_1, \ldots, a^n_1)$ is a common vertex of
  all such cells.)
\end{coro}

\begin{rema}
  As polarizations preserve all graded Betti numbers, by Theorem
  \ref{thm:cell-linres} every cell complex $C(T)$  has the same
  $f$-vector as the cell complex $C(T^{st})$ from the standard
  polarization. This is useful in the examples of Section \ref{sec:hypercube}.
\end{rema}

\newcommand{\imm}{{\mathfrak m}}
\usetikzlibrary{math, backgrounds,quotes, positioning, babel, 3d, graphs, graphs.standard}

\section{Simplicial complexes and hypercubes}
\label{sec:hypercube}

We look at the case when $\Pi(\Ab)$ is a hypercube.
The Artin monomial ideal then contains all squares of the variables.
The nice feature is that such ideals are equivalent to give a square
free monomial ideal and so a simplicial complex. In a number of
examples, for various simplicial complexes, we describe the possible
associated \rac cell complexes $C(T)$. This enables the description
of all possible polarizations, up to isomorphism.

\subsection{Hypercubes and ideals
  containing the squares of variables}
  
A subset  $S \sus V$, gives rise to a geometric simplex
in $\RR^V$:
\[ \Delta(S) = \left\{ \sum_{s \in S} \delta_s \cdot s \, | \,
  0 \leq \delta_s \leq 1, \, \sum_s \delta_s = 1 \right\}.
 \]
It also gives rise to a (hyper)cube in $\RR^V$.
\[ \square(S) = \left\{  \sum_{s \in S} \gamma_s \cdot s \, | \,
 0 \leq \gamma_s \leq 1 \right\}. \]

From  a combinatorial simplicial complex $X$ on a set $V$,
we get a geometric realization
$\Delta(X) \sus \RR^V$ consisting of all $\Delta(S)$ for $S \in X$.
We also get a cubical cell complex $\square(X)$ consisting of all
$\square(S)$ for $S \in X$. (Note that $\square(X)$ is contractible,
as each $\square(S)$ contains the origin.)

\begin{rema} The graph which is the one-skeleton of $\square(X)$
  is called the {\it simplex graph} in the literature. The book
  \cite{GC} considers graphs which may be embedded on hypercubes.
\end{rema}

Consider now those Artin monomial ideals
$\ijfr$ which contain the squares of all variables.
We may then write
\[ \ijfr = \imm_{\square} + \ijfr_{\Delta}, \quad \imm_\square = (x_1^2, x_2^2,
  \ldots, x_n^2), \] where
$\ijfr_\Delta$ is a squarefree monomial ideal corresponding to
a simplicial complex $X$. Each $A_i = \{a_i, b_i\}$ has then
cardinality two. 
The hypercube 
\begin{equation} \label{eq:mid-hypercube}
  \Delta(A_1) \times \Delta(A_2) \times \cdots \times \Delta(A_n).
\end{equation}
may be identified with $[0,1]^n$.
For $S \sus [n]$, a vertex $e_S = \sum_{i \in S} e_i$ of $[0,1]^n$ corresponds
to $(u_1, \ldots, u_n)$ in \eqref{eq:mid-hypercube}
with $u_i = b_i$ if $i \in S$ and $u_i = a_i$ otherwise.
The hypercube \eqref{eq:mid-hypercube}
has a natural monomial labelling, where a vertex may be
labelled with the monomial $u_1u_2 \cdots u_n$, see Figure
\ref{fig:scx-moncube}.

\begin{figure}[ht]
\begin{center}
\begin{tikzpicture}[thick, scale=2]
\coordinate[label=45:{$a_1b_2a_3$}] (B1) at (0,0,0);
\coordinate[label={$a_1b_2b_3$}] (B2) at (0,1,0);
\coordinate[label={$b_1b_2b_3$}] (B3) at (1,1,0);
\coordinate[label=right:{$b_1b_2a_3$}] (B4) at (1,0,0);
\coordinate[label=below:{$a_1a_2a_3$}] (B5) at (0,0,1);
\coordinate[label=left:{$a_1a_2b_3$}] (B6) at (0,1,1);
\coordinate[label=135:{$b_1a_2b_3 \hskip -2.2mm$}] (B7) at (1,1,1);
\coordinate[label=below:{$b_1a_2a_3$}] (B8) at (1,0,1);

\begin{scope}[dashed]
\draw[] (B1)--(B2);
\draw[] (B1)--(B5);
\draw[] (B4)--(B1);
\end{scope}
\draw[] (B2)--(B3);
\draw[] (B3)--(B4);
\draw[] (B6)--(B7);
\draw[] (B7)--(B8);
\draw[] (B2)--(B6);
\draw[] (B3)--(B7);
\draw[] (B4)--(B8);
\draw[] (B5)--(B6);
\draw[] (B8)--(B5);
\end{tikzpicture}
\end{center}
\caption{The cube with monomial labelling of vertices}
\label{fig:scx-moncube}
\end{figure}

In the following, for a subset $S$ of $[n]$ write $x(S)$ for
the ideal $(x_i)_{i \in S}$ and $x^2(S) = (x_i^2)_{i \in S}$.
It is well known that the minimal primary decomposition of
$I_X $ is
\begin{equation} \label{eq:mid-prim}
  I_X = \underset{\text{facet} F \text{ of } X} \bigcap  x(F^c),
  \end{equation}
see \cite[Lemma 1.5.4]{HeHi} or \cite[Theorem 1.7]{Mi-St}.

\begin{prop}
  The minimal irreducible decomposition of $\ijfr = \imm_\square + \ijfr_\Delta$
  is
  \[ \ijfr = \underset{{\rm{facet }} F {\rm{ of }} X}\bigcap (x(F^c) + x^2(F)). \] 

  In particular in any polarization of $\ijfr$, the maximal cells
  of $C(T)$ are in bijection with the facets $F$ of $X$.
  The dimension of a maximal cell is the cardinality $|F|$, one
  more than the dimension of the face $F$.  Each maximal cell
  is a sub-hypercube of the hypercube \eqref{eq:mid-hypercube}.
\end{prop}

\begin{proof}
  It is clear that we have
  \[ \ijfr = \imm_\square + \ijfr_\Delta \sus \bigcap_{\textrm{facet } F \textrm{ of } X}
    (x(F^c) + x^2(F)). \]
  It is clear that any monomial containing the square of a variable is contained
  in the left side. If $m$ is a squarefree monomial contained in
  the right side, it is also clear by \eqref{eq:mid-prim}, that $m$
  is in the left side.
  Whence we have equality above.
\end{proof}

\begin{prop}
  Let $J^\st = J(T^{st})$ be the standard polarization of $\ijfr$,
  where $\ijfr_\Delta$
  corresponds to a simplicial complex $X$. The cell complex
  $C(T^{st})$ is the canonical cubical complex $\square(X)$. With 
  its natural monomial labelling, it gives a free cellular resolution of
  the Alexander dual rainbow ideal $I^{st}$ of $J^{st}$.
\end{prop}

\begin{proof}
  This follows by Theorem \ref{thm:cell-linres} and Corollary
  \ref{coro:mid-cellstd}. 
  \end{proof}

\subsection{Examples}
In the following $X$ will be the simplicial complex whose
Stanley-Reisner ideal is $\ijfr_\Delta$. In our examples the
dimension of $X$ will be $\leq 1$. In a polarization with
associated cell complex $C(T)$, a vertex of $X$ corresponds
to an edge of $C(T)$ and an edge of $X$ gives a square in $C(T)$.
This correspondence is one-to-one between facets of $X$ and
maximal cells of $C(T)$ (but not in general
between faces and cells).

Furthermore, if two squares in $C(T)$ have an edge in common, the
corresponding edges in $X$ have a corresponding vertex in common,
by Lemma \ref{lem:mid-dim}. But conversely, if two edges in $X$ have
a vertex in common, the corresponding squares in $C(T)$ may or may
not have an edge in common, in fact the squares may even be disjoint
(see Example \ref{ex:scx-star}). 
The essential thing concerning $C(T)$ is that it is r-acyclic.

Instead of writing the variables as $x_1, x_2, \ldots$,
we write $x,y, \ldots$. We 
write $A_x = \{ x_a, x_b\}$ for the variables mapping to $x$,
  and similarly $A_y$ and so on.

\begin{exam}
  Let $\ijfr = (x^2, y^2) \sus k[x,y]$. It has a unique polarization $J$
  with Alexander dual $I$:
  \[ J = (x_ax_b, y_ay_b) \sus S = k[x_a, x_b, y_a, y_b],
    \quad I = (x_ay_a, x_ay_b, x_by_a, x_by_b). \]
  The associated $T$ is the whole $A_x \times A_y$. So
  $C(T)$ is the square with vertices labelled by monomials as in
  the left in Figure \ref{fig:scx-square}.

\begin{figure}[ht]
\begin{center}
       \begin{NiceTabular}{ccc}
       \Block[t,l]{1-1}{
\begin{tikzpicture}[thick,scale=2]
\coordinate[label=-135:{$x_ay_a$}] (A1) at (0,0);
\coordinate[label=135:{$x_ay_b$}] (A2) at (0,1);
\coordinate[label=45:{$x_by_b$}] (A3) at (1,1);
\coordinate[label=-45:{$x_by_a$}] (A4) at (1,0);

\draw (A1) edge[] (A2); 
\draw (A2) edge[] (A3); 
\draw (A4) edge[swap] (A3); 
\draw (A1) edge[swap] (A4);

\begin{scope}[on background layer]
\draw[fill=lightgray, opacity=0.4] (A1)--(A2)--(A3)--(A4)--(A1);
\end{scope}

\end{tikzpicture}}
&
\hspace{1cm}
&
\begin{tikzpicture}[thick,scale=2]

\coordinate (A1) at (0,0); 
\coordinate (A2) at (0,1); 
\coordinate (A3) at (1,1); 
\coordinate (A4) at (1,0); 

\begin{scope}[font=\scriptsize]
\node[label=45:{$a$}] (n2) at (0,-0.1) {};
\node[label=45:{$a$}] (n3) at (-0.1,0) {};
\node[label=-45:{$a$}] (n4) at (0,1.1) {};
\node[label=-45:{$b$}] (n7) at (-0.1,1) {};
\node[label=135:{$a$}] (n5) at (1+0.1,0.1-0.1) {};
\node[label=135:{$b$}] (n6) at (0.9+0.1,0-0.1) {};
\node[label=-135:{$b$}] (n8) at (0.9+0.1,1+0.1) {};
\node[label=-135:{$b$}] (n9) at (1+0.1,0.9+0.1) {};
\end{scope}

\draw (A1) edge["$y$"] (A2); 
\draw (A2) edge["$x$"] (A3); 
\draw (A4) edge["$y$", swap] (A3); 
\draw (A1) edge["$x$", swap] (A4);

\begin{scope}[on background layer]
\draw[fill=lightgray, opacity=0.4] (A1)--(A2)--(A3)--(A4)--(A1);
\end{scope}
\end{tikzpicture}\\
\end{NiceTabular}
\end{center}
\caption{The cube with monomial labelling}
\label{fig:scx-square}
\end{figure}
We write it somewhat more stylistic as in the right of Figure
\ref{fig:scx-square}.
The free cellular resolution of $I$ becomes:
\[ I \vpil S(-2)^4 \vpil S(-3)^4 \vpil S(-4). \]
\end{exam}

\begin{exam} \label{ex:hypercube-tree}
  \[ \ijfr = (x^2, y^2, z^2) + (xy, xz, yz).\]
  This is the same Artin monomial ideal we considered
  in Example \ref{ex:intro-m2} in the Introduction.
The simplicial complex associated to $\ijfr_\Delta$ is
three isolated points, labelled respectively $x,y$ and $z$.
In a polarization  $J(T)$, by Theorem \ref{thm:mid-mid} the cell complex
$C(T)$ will
have three line segments, labelled respectivey $x, y$ and $z$.
The cell complex is r-acyclic and so must be a tree. This gives, up to
isomorphism, the two possibilities given in Figure \ref{fig:scx-trees}.

\begin{figure}[ht]
\begin{center}
  \begin{NiceTabular}{ccc}
     \Block[t,l]{1-1}{
 \begin{tikzpicture}[thick, dot/.style={draw,fill,circle,inner sep=1.5pt},scale=1.5]
\begin{scope}[font=\scriptsize]
\node[dot] (n1) at (0,0) {};
\node[dot, label=90:{$a$}, label=right:{$a$}] (n2) at (1,0) {};
\node[dot, label=left:{$b$}] (n3) at (1 + 0.866, 0 + 0.5) {};
\node[dot, label=135:{$b$}] (n4) at (1 + 0.866, 0 - 0.5) {};

\node[label=135:{$a$}] (n5) at (1+0.05,0-0.1) {};
\node[label=45:{$b$}] (n6) at (0-0.05,0-0.1) {};
\end{scope}

\draw (n1) edge["$x$", swap] (n2);
\draw (n2) edge["$z$"] (n3);
\draw (n2) edge["$y$", swap] (n4);
\end{tikzpicture}} &
  \qquad \quad 
  &     
\begin{tikzpicture}[thick, dot/.style={draw,fill,circle,inner sep=1.5pt},scale=1.5]
\begin{scope}[font=\scriptsize]
\node[dot, label=45:{$b$}] (n1) at (0,0) {}; 
\node[dot, label=45:{$b$}, label=135:{$a$}] (n2) at (1,0) {};
\node[dot, label=45:{$b$}, label=135:{$a$}] (n3) at (2,0) {};
\node[dot, label=135:{$a$}] (n4) at (3,0) {};
\end{scope}

\draw (n1) edge["$x$", swap] (n2);
\draw (n2) edge["$y$", swap] (n3);
\draw (n3) edge["$z$", swap] (n4);
\end{tikzpicture}
\\
\end{NiceTabular}
\end{center}
\caption{The $C(T)$ are trees with labellings}
\label{fig:scx-trees}
\end{figure}
These are embedded on the cube as in Figure \ref{fig:scx-embedcube}.
\begin{figure}[ht]
  \begin{center}
\begin{tikzpicture}[thick, scale=2]
\coordinate (B1) at (0,0,0);
\coordinate (B2) at (0,1,0);
\coordinate[label=right:{$bbb$}] (B3) at (1,1,0);
\coordinate (B4) at (1,0,0);
\coordinate[label=left:{$aaa$}] (B5) at (0,0,1);
\coordinate (B6) at (0,1,1);
\coordinate (B7) at (1,1,1);
\coordinate (B8) at (1,0,1);

\begin{scope}[densely dotted]
\draw[] (B1)--(B2);
\draw[] (B2)--(B3);
\draw[] (B3)--(B4);
\draw[] (B4)--(B1);
\draw[] (B6)--(B7);
\draw[] (B7)--(B8);
\draw[] (B2)--(B6);
\draw[] (B3)--(B7);
\draw[] (B4)--(B8);
\end{scope}
\draw[red, very thick] (B5)--(B6);
\draw[red, very thick] (B8)--(B5);
\draw[red, very thick] (B1)--(B5);
\end{tikzpicture}
\hspace{1cm}
\begin{tikzpicture}[thick, scale=2]
\coordinate (B1) at (0,0,0);
\coordinate (B2) at (0,1,0);
\coordinate[label=right:{$bbb$}] (B3) at (1,1,0);
\coordinate (B4) at (1,0,0);
\coordinate[label=left:{$aaa$}] (B5) at (0,0,1);
\coordinate (B6) at (0,1,1);
\coordinate (B7) at (1,1,1);
\coordinate (B8) at (1,0,1);

\begin{scope}[densely dotted]
\draw[] (B1)--(B2);
\draw[] (B3)--(B4);
\draw[] (B4)--(B1);
\draw[] (B6)--(B7)--(B8)--(B5);
\draw[] (B1)--(B5);
\draw[] (B3)--(B7);
\draw[] (B4)--(B8);
\end{scope}
\begin{scope}[red, very thick]
\draw[] (B2)--(B3);
\draw[] (B2)--(B6);
\draw[] (B5)--(B6);
\end{scope}
\end{tikzpicture}
\end{center}
\caption{The $C(T)$-trees with their embedding on cubes}
\label{fig:scx-embedcube}
\end{figure}
Using Corollary \ref{cor:mid-cover}, 
the corresponding polarizations are then
\begin{align*}
  J_1 & = (x_ax_b, y_ay_b, z_az_b) + (x_ay_a, x_az_a, y_az_a) \\
  J_2 & = (x_ax_b, y_ay_b, z_az_b) + (x_ay_b, x_az_b, y_az_b)
\end{align*}
More simply, using Corollary \ref{cor:mid-cover} to obtain the
polarization $J_2$, amounts to observe for the tree to the right
in Figure \ref{fig:scx-trees}: On the path from edge $x$ to edge $y$
the label on the vertex of $x$ is $a$, and the vertex on $y$ is $b$.
Similarly for the paths form $x$ to $z$ and the path from $y$ to $z$.
See also how to obtain the ideal $J_3$ below from the tree in Figure
\ref{fig:scx-treeB}.

The ideals $J_1$ and $J_2$ are the
ideals given in Example \ref{ex:intro-m2} in
the Introduction.
In Figure \ref{fig:scx-trees} for the edges
labelled $x$ we are free to label its two
vertices with $a$ and $b$ as we like, and
similarly for the $y,z$ edges.
We could have labelled the line graph with $a$ and $b$ as in
Figure \ref{fig:scx-treeB}, which
would have given the embedding on the cube as shown there.

\begin{figure}[ht]
\begin{center}
       \begin{NiceTabular}{ccc}
       \Block[t,l]{1-1}{
\begin{tikzpicture}[thick, dot/.style={draw,fill,circle,inner sep=1.5pt},scale=1.5]
\begin{scope}[font=\scriptsize]
\node[dot, label=45:{$a$}] (n1) at (0,0) {}; 
\node[dot, label=45:{$b$}, label=135:{$b$}] (n2) at (1,0) {};
\node[dot, label=45:{$a$}, label=135:{$a$}] (n3) at (2,0) {};
\node[dot, label=135:{$b$}] (n4) at (3,0) {};
\end{scope}

\draw (n1) edge["$x$", swap] (n2);
\draw (n2) edge["$y$", swap] (n3);
\draw (n3) edge["$z$", swap] (n4);
\end{tikzpicture}} &
\hspace{1cm}
&
\begin{tikzpicture}[thick, scale=2]
\coordinate (B1) at (0,0,0);
\coordinate (B2) at (0,1,0);
\coordinate[label=right:{$bbb$}] (B3) at (1,1,0);
\coordinate (B4) at (1,0,0);
\coordinate[label=left:{$aaa$}] (B5) at (0,0,1);
\coordinate (B6) at (0,1,1);
\coordinate (B7) at (1,1,1);
\coordinate (B8) at (1,0,1);

\begin{scope}[densely dotted]
\draw[] (B1)--(B2);
\draw[] (B2)--(B3);
\draw[] (B3)--(B4);
\draw[] (B5)--(B6);
\draw[] (B6)--(B7);
\draw[] (B8)--(B5);
\draw[] (B1)--(B5);
\draw[] (B2)--(B6);
\draw[] (B3)--(B7);
\end{scope}
\begin{scope}[red, very thick]
\draw[] (B4)--(B1);
\draw[] (B4)--(B8);
\draw[] (B7)--(B8);
\end{scope}
\end{tikzpicture}\\
\end{NiceTabular}
\end{center}
\caption{Alternative labelling of $a$ and $b$, with corresponding
embedding on cube}
\label{fig:scx-treeB}
\end{figure}
This gives a polarization isomorphic to $J_2$:
\begin{equation*}
  J_3  =  (x_ax_b, y_ay_b, z_az_b) + (x_by_b, x_bz_a, y_az_a)
\end{equation*}
In general for the powers $(x_1, \ldots, x_n)^2$, its polarizations
are classified, up to isomorphism, by trees with $n$ edges,
\cite[Section 7]{AFL}. A comprehensive study of ideals obtained
from such polarizations is done in \cite{FO}, see also the present Subsection
\ref{subsec:further-regquot}.
\end{exam}

\begin{exam} \label{ex:cell-lpp}
  \[ \ijfr = (x^2, y^2, z^2, w^2) + (xz, xw, yz,yw, zw). \]
  The simplicial complex $X$ associated to $\ijfr_\Delta $ is given
  in Figure \ref{fig:scx-lpp}.

\begin{figure}[ht]
\begin{center}  
\begin{tikzpicture}[thick, dot/.style={draw,fill,circle,inner sep=1.5pt},scale=1.5]
\node[dot, label=left:{$y$}] (n1) at (0,0) {};
\node[dot, label=left:{$x$}] (n2) at (0,1) {};
\node[dot, label=right:{$w$}] (n3) at (1,0) {};
\node[dot, label=right:{$z$}] (n4) at (1,1) {};

\draw (n1) edge[] (n2);
\end{tikzpicture}
\end{center}
\caption{Simplicial complex: A line segment and two isolated points}
\label{fig:scx-lpp}
\end{figure}
In a polarization $J(T)$, by Theorem \ref{thm:mid-mid} the cell complex
$C(T)$ will then
  have one square, corresponding to the line segment labelled $xy$, and two
  line segments corresponding to the points labelled $z$ and $w$.
  The cell complex $C(T)$ must also be r-acyclic. This gives, up to
  isomorphism, the
  four possibilities for $C(T)$ given in Figure \ref{fig:scx-sll}.
\begin{figure}[ht] 
\begin{center}
\begin{tikzpicture}[thick,scale=1]

\coordinate(A1) at (0,0);
\coordinate(A2) at (0,1);
\coordinate(A3) at (1,1);
\coordinate(A4) at (1,0);
\coordinate(A5) at (1,2);
\coordinate(A6) at (2,1);
\draw[fill=lightgray, opacity=0.4] (A1)--(A2)--(A3)--(A4)--(A1);
\draw (A1) edge["$y$"] (A2);
\draw (A2) edge["$x$"] (A3);
\draw (A3) edge["$y$"] (A4);
\draw (A4) edge["$x$"] (A1);
\draw (A3) edge["$w$", swap] (A5);
\draw (A3) edge["$z$", swap] (A6);
\end{tikzpicture}
\hspace{0.5cm}
\begin{tikzpicture}[thick,scale=1]

\coordinate(A1) at (0,0);
\coordinate(A2) at (0,1);
\coordinate(A3) at (1,1);
\coordinate(A4) at (1,0);
\coordinate(A5) at (2,0);
\coordinate(A6) at (2,1);
\draw[fill=lightgray, opacity=0.4] (A1)--(A2)--(A3)--(A4)--(A1);
\draw (A1) edge["$y$"] (A2);
\draw (A2) edge["$x$"] (A3);
\draw (A3) edge["$y$"] (A4);
\draw (A4) edge["$x$"] (A1);
\draw (A4) edge["$z$", swap] (A5);
\draw (A3) edge["$w$"] (A6);
\end{tikzpicture}
\hspace{0.5cm}
\begin{tikzpicture}[thick,scale=1]

\coordinate(A1) at (0,0);
\coordinate(A2) at (0,1);
\coordinate(A3) at (1,1);
\coordinate(A4) at (1,0);
\coordinate(A5) at (-1,0);
\coordinate(A6) at (2,1);
\draw[fill=lightgray, opacity=0.4] (A1)--(A2)--(A3)--(A4)--(A1);
\draw (A1) edge["$y$"] (A2);
\draw (A2) edge["$x$"] (A3);
\draw (A3) edge["$y$"] (A4);
\draw (A4) edge["$x$"] (A1);
\draw (A1) edge["$z$"] (A5);
\draw (A3) edge["$w$"] (A6);
\end{tikzpicture}
\hspace{0.5cm}
\begin{tikzpicture}[thick,scale=1]

\coordinate(A1) at (0,0);
\coordinate(A2) at (0,1);
\coordinate(A3) at (1,1);
\coordinate(A4) at (1,0);
\coordinate(A5) at (1,2);
\coordinate(A6) at (2,2);
\draw[fill=lightgray, opacity=0.4] (A1)--(A2)--(A3)--(A4)--(A1);
\draw (A1) edge["$y$"] (A2);
\draw (A2) edge["$x$"] (A3);
\draw (A3) edge["$y$"] (A4);
\draw (A4) edge["$x$"] (A1);
\draw (A3) edge["$z$", swap] (A5);
\draw (A5) edge["$w$"] (A6);
\end{tikzpicture}
\end{center}
\caption{Restriction-acyclic $C(T)$ with one square and two line segments}
\label{fig:scx-sll}
\end{figure}
  We thus get, up to isomorphism, four polarizations of $\ijfr $.
  Using Theorem \ref{thm:mid-mid} and Corollary \ref{cor:mid-cover}
  they are seen to be:
  \begin{align*}
    J_1& = (x_ax_b, y_ay_b, z_az_b, w_aw_b)
         + (x_az_a,x_aw_a, y_az_a, y_aw_a, z_aw_a)  \\
    J_2 & =(x_ax_b, y_ay_b, z_az_b, w_aw_b)
         + (x_az_a,x_aw_a, y_bz_a, y_aw_a, z_aw_a)  \\
   J_3 & =(x_ax_b, y_ay_b, z_az_b, w_aw_b)
         + (x_bz_a,x_aw_a, y_bz_a, y_aw_a, z_aw_a)  \\
 J_4 & =(x_ax_b, y_ay_b, z_az_b, w_aw_b)
         + (x_az_a,x_aw_a, y_az_a, y_aw_a, z_bw_a)     
  \end{align*}
 (Again for each variable, the labels $a$ and $b$ may be switched (or not)
  and we will get isomorphic polarizations, which are embedded distinctly
  on the four-dimensional hypercube.) The standard polarization $J_1$
  corresponds to the first cell complex in Figure \ref{fig:scx-sll}.
 Taking the Alexander duals of the above, we get rainbow ideals:
  \begin{align*}
  I_1 & = (x_ay_az_aw_a, x_by_az_aw_a, x_ay_bz_aw_a, x_by_bz_aw_a,        
          x_ay_az_bw_a,x_ay_az_aw_b) \\
  I_2 & = ( x_ay_az_aw_a, x_by_az_aw_a, x_ay_bz_aw_a, x_by_bz_aw_a,        
         x_ay_bz_bw_a,x_ay_az_aw_b) \\
  I_3 & = ( x_ay_az_aw_a, x_by_az_aw_a, x_ay_bz_aw_a, x_by_bz_aw_a,        
        x_by_bz_bw_a,x_ay_az_aw_b) \\
  I_4 & = ( x_ay_az_aw_a, x_by_az_aw_a, x_ay_bz_aw_a, x_by_bz_aw_a,        
        x_ay_az_bw_a,x_ay_az_bw_b)
  \end{align*}
  The minimal free resolution for each of these is the free cellular complex
  given by $C(T)$. It has in each case the form
  \[ S(-4)^6 \vpil S(-5)^6 \vpil S(-6). \]

  The triangulated balls defined by the ideals $J_i$ each have six
  facets which are complements of the variable sets in each of the
  generators of the rainbow ideals
  \begin{align*}
    F1: \quad &  \begin{matrix} \{x_b, y_b, z_b, w_b\} \\ \{x_a, y_b, z_b, w_b\}
    \end{matrix},  \, \,
    \begin{matrix} \{x_b, y_a, z_b, w_b\} \\ \{x_a, y_a, z_b, w_b\}
    \end{matrix}, \, \,
    \begin{matrix}  \{x_b, y_b, z_a, w_b\} \\  \{x_b, y_b, z_b, w_a\}
    \end{matrix} \\
    F2: \quad & \begin{matrix} \{x_b, y_b, z_b, w_b\} \\ \{x_a, y_b, z_b, w_b\}
    \end{matrix}, \, \,
    \begin{matrix} \{x_b, y_a, z_b, w_b\} \\
      \{x_a, y_a, z_b, w_b\}
    \end{matrix}, \, \, 
    \begin{matrix} \{x_b, y_a, z_a, w_b\} \\  \{x_b, y_b, z_b, w_a\} 
    \end{matrix} \\
    F3: \quad & \begin{matrix} \{x_b, y_b, z_b, w_b\} \\ \{x_a, y_b, z_b, w_b\}
    \end{matrix}, \, \,
    \begin{matrix} \{x_b, y_a, z_b, w_b\} \\
      \{x_a, y_a, z_b, w_b\}
    \end{matrix}, \, \, 
    \begin{matrix} \{x_a, y_a, z_a, w_b\} \\  \{x_b, y_b, z_b, w_a\} 
    \end{matrix} \\
    F4:\quad & \begin{matrix} \{x_b, y_b, z_b, w_b\} \\ \{x_a, y_b, z_b, w_b\}
    \end{matrix}, \, \,
    \begin{matrix} \{x_b, y_a, z_b, w_b\} \\
      \{x_a, y_a, z_b, w_b\}
    \end{matrix}, \, \, 
    \begin{matrix} \{x_b, y_b, z_a, w_b\} \\  \{x_b, y_b, z_a, w_a\} 
    \end{matrix}
  \end{align*}

\end{exam}

\begin{rema} \label{rem:cube-acyclic}
  The third cell complex in Figure \ref{fig:scx-sll} can be embedded on the
  $3$-dimensional cube. The $w$-label is changed to $z$. Note however
  that for this embedding the cell complex is {\it not} r-acyclic,
  as restricting to one of the horizontal levels, the cell complex becomes
  two separate points.
  \end{rema}

\begin{exam} 
  \[ \ijfr = (x^2, y^2, z^2, w^2, t^2, u^2) +
    (xz,xw,yz,yw, xt, xu, yt, yu, zt, zu, wt, wu). \]
   The simplicial complex $X$ associated to $\ijfr_\Delta $ is given
  in Figure \ref{fig:scx-threelines}.
\begin{figure}[ht]
\begin{center}  
\begin{tikzpicture}[thick, dot/.style={draw,fill,circle,inner sep=1.5pt},scale=1.5]
\node[dot, label=left:{$y$}] (n1) at (0,0) {};
\node[dot, label=left:{$x$}] (n2) at (0,1) {};

\node[dot, label=left:{$w$}] (n3) at (1,0) {};
\node[dot, label=left:{$z$}] (n4) at (1,1) {};

\node[dot, label=left:{$u$}] (n5) at (2,0) {};
\node[dot, label=left:{$t$}] (n6) at (2,1) {};

\draw (n1) edge[] (n2);
\draw (n3) edge[] (n4);
\draw (n5) edge[] (n6);
\end{tikzpicture}
\end{center}
\caption{Three line segments}
\label{fig:scx-threelines}
\end{figure}
In a polarization $J(T)$, the cell complex $C(T)$ will then
  have three squares, corresponding to the lines labelled $xy, zw$ and
  $tu$.
  The cell complex $C(T)$ must also be \rac. By Corollary
  \ref{lem:mid-dim}, each pair of squares can intersect in at most
  a point. This gives, up to
  isomorphism, the
  three possibilities for $C(T)$ given in Figure \ref{fig:scx-threesquares}.
\begin{figure}[ht]
\begin{center}
\begin{tikzpicture}[thick, font=\footnotesize]
\coordinate (A1) at (0,0);
\coordinate (A2) at (0,1);
\coordinate (A3) at (1,0);
\coordinate (A4) at (1,1);

\draw[fill=lightgray, opacity=0.4] (A1)--(A2)--(A4)--(A3)--(A1);

\draw (A1) edge["$w$"] (A2);
\draw (A2) edge["$z$"] (A4);
\draw (A4) edge["$w$"] (A3);
\draw (A3) edge["$z$"] (A1);

\coordinate (A5) at (-1,0);
\coordinate (A6) at (0,-1);
\coordinate (A7) at (-1,-1);

\draw[fill=lightgray, opacity=0.4] (A1)--(A5)--(A7)--(A6)--(A1);

\draw (A1) edge["$x$", swap] (A5);
\draw (A5) edge["$y$", swap] (A7);
\draw (A7) edge["$x$", swap] (A6);
\draw (A6) edge["$y$", swap] (A1);

\coordinate (A8) at (2,0);
\coordinate (A9) at (1,-1);
\coordinate (A10) at (2,-1);

\draw[fill=lightgray, opacity=0.4] (A3)--(A9)--(A10)--(A8)--(A3);

\draw (A3) edge["$u$", swap] (A9);
\draw (A9) edge["$t$", swap] (A10);
\draw (A10) edge["$u$", swap] (A8);
\draw (A8) edge["$t$", swap] (A3);
\end{tikzpicture}
\begin{tikzpicture}[thick]
\coordinate (A1) at (0,0);
\coordinate (A2) at (0,1);
\coordinate (A3) at (1,0);
\coordinate (A4) at (1,1);

\draw[fill=lightgray, opacity=0.4] (A1)--(A2)--(A4)--(A3)--(A1);

\draw (A1) edge[font=\footnotesize,"$w$"] (A2);
\draw (A2) edge[font=\footnotesize,"$z$"] (A4);
\draw (A4) edge[font=\footnotesize,"$w$"] (A3);
\draw (A3) edge[font=\footnotesize,"$z$"] (A1);

\coordinate (A5) at (-1,0);
\coordinate (A6) at (0,-1);
\coordinate (A7) at (-1,-1);

\draw[fill=lightgray, opacity=0.4] (A1)--(A5)--(A7)--(A6)--(A1);

\draw (A1) edge[font=\footnotesize, "$x$", swap] (A5);
\draw (A5) edge[font=\footnotesize, "$y$", swap] (A7);
\draw (A7) edge[font=\footnotesize,"$x$", swap] (A6);
\draw (A6) edge[font=\footnotesize,"$y$", swap] (A1);

\coordinate (A8) at (2,2);
\coordinate (A9) at (1,2);
\coordinate (A10) at (2,1);

\draw[fill=lightgray, opacity=0.4] (A4)--(A10)--(A8)--(A9)--(A4);

\draw (A4) edge[font=\footnotesize,"$t$", swap] (A10);
\draw (A10) edge[font=\footnotesize,"$u$", swap] (A8);
\draw (A8) edge[font=\footnotesize,"$t$", swap] (A9);
\draw (A9) edge[font=\footnotesize,"$u$", swap] (A4);
\end{tikzpicture}
\begin{tikzpicture}[thick, font=\footnotesize]
\coordinate (A1) at (0,0);
\coordinate (A2) at (0,1);
\coordinate (A3) at (1,0);
\coordinate (A4) at (1,1);

\draw[fill=lightgray, opacity=0.4] (A1)--(A2)--(A4)--(A3)--(A1);

\draw (A1) edge["$w$", swap] (A2);
\draw (A2) edge["$z$"] (A4);
\draw (A4) edge["$w$"] (A3);
\draw (A3) edge["$z$", swap] (A1);

\begin{scope}[rotate=-60]
\coordinate (A5) at (-1,0);
\coordinate (A6) at (0,-1);
\coordinate (A7) at (-1,-1);

\draw[fill=lightgray, opacity=0.4] (A1)--(A5)--(A7)--(A6)--(A1);

\draw (A1) edge["$x$"] (A5);
\draw (A5) edge["$y$", swap] (A7);
\draw (A7) edge["$x$", swap] (A6);
\draw (A6) edge["$y$"] (A1);
\end{scope}

\begin{scope}[rotate=-30]
\coordinate (A8) at (1,0);
\coordinate (A9) at (0,-1);
\coordinate (A10) at (1,-1);

\draw[fill=lightgray, opacity=0.4] (A1)--(A9)--(A10)--(A8)--(A1);

\draw (A1) edge["$u$"] (A9);
\draw (A9) edge["$t$", swap] (A10);
\draw (A10) edge["$u$", swap] (A8);
\draw (A8) edge["$t$"] (A1);
\end{scope}
\end{tikzpicture}
\end{center}
\caption{Restriction-acyclic $C(T)$ of three squares without edge intersections}
\label{fig:scx-threesquares}
\end{figure}
 We thus get, up to isomorphism, three polarizations of $\ijfr $.
  Using Theorem \ref{thm:mid-mid} and Corollary \ref{cor:mid-cover}
  they are seen to be
  (write $M_\square = (x_ax_b, y_ay_b, z_az_b, w_aw_b, t_at_b, u_au_b)$)
  \begin{align*}
    J_1& = M_\square
       + (x_az_b,x_aw_b, y_az_b, y_aw_b, x_at_b, x_au_a, y_at_b, y_au_a,
            z_at_b, z_au_a, w_bt_b, w_bu_a)  \\
   J_2& = M_\square
         + (x_az_b,x_aw_b, y_az_b, y_aw_b, x_at_b, x_au_b, y_at_b, y_au_b,
        z_at_b, z_au_b, w_at_b, w_au_b)  \\
    J_3& = M_\square
         + (x_az_a,x_aw_a, y_az_a, y_aw_a, x_at_a, x_au_a, y_at_a, y_au_a,
            z_at_a, z_au_a, w_at_a, w_au_a) 
  \end{align*}

\end{exam}

  \begin{exam} 
  \[ \ijfr = (x^2, y^2, z^2, w^2) +
    (xz,xw,yw,zw). \]
  The simplicial complex $X$ associated to $\ijfr_\Delta $ is given
in Figure \ref{fig:scx-llp}.
  
\begin{figure}[ht]
\begin{center}  
\begin{tikzpicture}[thick, dot/.style={draw,fill,circle,inner sep=1.5pt},scale=1.5]
\node[dot, label=left:{$y$}] (n1) at (0,0) {};
\node[dot, label=left:{$z$}] (n2) at (0.707, 0.707) {};

\node[dot, label=left:{$w$}] (n3) at (1.414,0) {};
\node[dot, label=left:{$x$}] (n4) at (0.707,-0.707) {};

\draw (n1) edge[] (n2);
\draw (n1) edge[] (n4);

\end{tikzpicture}
\end{center}
\caption{Two line segments and an isolated point}
\label{fig:scx-llp}
\end{figure}
  In a polarization $J(T)$, the cell complex $C(T)$ will then
  have two squares, corresponding to the lines labelled $xy$ and $ yz$,
  and one line segment, labelled $w$. 
  The cell complex $C(T)$ must  be \rac.  This gives, up to
  isomorphism, the
  two possibilities for $C(T)$ given in Figure \ref{fig:scx-ssl}.
\begin{figure}[ht]
\begin{center}
\begin{tikzpicture}[thick, scale=1.5]
\coordinate (A1) at (0,0);
\coordinate (A2) at (0,1);
\coordinate (A3) at (1,0);
\coordinate (A4) at (1,1);

\draw[fill=lightgray, opacity=0.4] (A1)--(A2)--(A4)--(A3)--(A1);

\draw (A1) edge["$y$"] (A2);
\draw (A2) edge["$x$"] (A4);
\draw (A4) edge["$y$", swap] (A3);
\draw (A3) edge["$x$"] (A1);
\coordinate (A5) at (2,0);
\coordinate (A6) at (2,1);

\draw[fill=lightgray, opacity=0.4] (A3)--(A5)--(A6)--(A4)--(A3);

\draw (A3) edge["$z$", swap] (A5);
\draw (A4) edge["$z$"] (A6);
\draw (A5) edge["$y$"] (A6);

\coordinate (A7) at (1 + 0.5, 1 + 0.866);

\draw (A4) edge["$w$"] (A7);
\end{tikzpicture}
\hspace{1cm}
\begin{tikzpicture}[thick, scale=1.5]
\coordinate (A1) at (0,0);
\coordinate (A2) at (0,1);
\coordinate (A3) at (1,0);
\coordinate (A4) at (1,1);

\draw[fill=lightgray, opacity=0.4] (A1)--(A2)--(A4)--(A3)--(A1);

\draw (A1) edge["$y$"] (A2);
\draw (A2) edge["$x$"] (A4);
\draw (A4) edge["$y$", swap] (A3);
\draw (A3) edge["$x$"] (A1);

\coordinate (A5) at (2,0);
\coordinate (A6) at (2,1);

\draw[fill=lightgray, opacity=0.4] (A3)--(A5)--(A6)--(A4)--(A3);

\draw (A3) edge["$z$", swap] (A5);
\draw (A4) edge["$z$"] (A6);
\draw (A5) edge["$y$"] (A6);

\coordinate (A7) at (2 + 0.5, 1 + 0.866);

\draw (A6) edge["$w$"] (A7);
\end{tikzpicture}
\end{center}
\caption{Restriction-acyclic $C(T)$ with two adjacent squares and one line segment}
\label{fig:scx-ssl}
\end{figure}
 We thus get, up to isomorphism, two polarizations of $\ijfr $.
  Using Theorem \ref{thm:mid-mid} and Corollary \ref{cor:mid-cover}
  they are seen to be:
  \begin{align*}
    J_1& = (x_ax_b, y_ay_b, z_az_b, w_aw_b)
         + (x_az_a,x_aw_a,y_aw_a,z_aw_a) \\
     J_2& = (x_ax_b, y_ay_b, z_az_b, w_aw_b)
          + (x_az_a,x_aw_a,y_aw_a,z_bw_a)
  \end{align*}

\end{exam}

   \begin{exam} \label{ex:scx-star}
  \[ \ijfr = (x^2, y^2, z^2, w^2) +
    (yz,yw,zw). \]
  The simplicial complex $X$ associated to $\ijfr_\Delta $ is given
  in Figure \ref{fig:scx-star}.
\begin{figure}[ht]
\begin{center}  
\begin{tikzpicture}[thick, dot/.style={draw,fill,circle,inner sep=1.5pt},scale=1.5]
\node[dot, label=-45:{$x$}] (n1) at (0,0) {};
\node[dot, label=right:{$z$}] (n2) at (1, 0) {};
\node[dot, label=left:{$w$}] (n3) at (-0.5, 0.866) {};
\node[dot, label=left:{$y$}] (n4) at (-0.5,-0.866) {};
\draw (n1) edge[] (n2);
\draw (n1) edge[] (n3);
\draw (n1) edge[] (n4);
\end{tikzpicture}
\end{center}
\caption{The star with three edges}
\label{fig:scx-star}
\end{figure}
  In a polarization $J(T)$, the cell complex $C(T)$ will then
  have three squares, corresponding to the lines labelled $xy, xz$ and $xw$. 
  The cell complex $C(T)$ must  be \rac.  This gives, up to
  isomorphism, the
  two possibilities for $C(T)$ given in Figure \ref{fig:scx-sss}.
\begin{figure}[ht]
\begin{center}
       \begin{NiceTabular}{ccc}
\begin{tikzpicture}[thick, scale=1.5]
\coordinate (A1) at (0,0);
\coordinate (A2) at (0,1);
\coordinate (A3) at (1,0);
\coordinate (A4) at (1,1);
\draw[fill=lightgray, opacity=0.4] (A1)--(A2)--(A4)--(A3)--(A1);

\draw (A1) edge["$x$"] (A2);
\draw (A2) edge["$y$"] (A4);
\draw (A4) edge["$x$", swap] (A3);
\draw (A3) edge["$y$"] (A1);
\coordinate (A5) at (1 + 0.707, 0 + 0.707);
\coordinate (A6) at (1 + 0.707, 1 + 0.707);
\draw[fill=lightgray, opacity=0.4] (A3)--(A5)--(A6)--(A4)--(A3);
\draw (A3) edge["$z$", swap] (A5);
\draw (A4) edge["$z$"] (A6);
\draw (A5) edge["$x$"] (A6);
\coordinate (A7) at (1 + 0.707, 0 - 0.707);
\coordinate (A8) at (1 + 0.707, 1 - 0.707);
\draw[fill=lightgray, opacity=0.4] (A3)--(A7)--(A8)--(A4)--(A3);
\draw (A3) edge["$w$", swap] (A7);
\draw (A4) edge["$w$"] (A8);
\draw (A7) edge["$x$"] (A8);
\end{tikzpicture}&
\hspace{1cm}
&
       \Block[t,l]{1-1}{
\begin{tikzpicture}[thick, scale=1.5]
\coordinate (A1) at (0,0);
\coordinate (A2) at (0,1);
\coordinate (A3) at (1,0);
\coordinate (A4) at (1,1);
\draw[fill=lightgray, opacity=0.4] (A1)--(A2)--(A4)--(A3)--(A1);
\draw (A1) edge["$x$"] (A2);
\draw (A2) edge["$y$"] (A4);
\draw (A4) edge["$x$", swap] (A3);
\draw (A3) edge["$y$"] (A1);
\coordinate (A5) at (2,0);
\coordinate (A6) at (2,1);
\draw[fill=lightgray, opacity=0.4] (A3)--(A5)--(A6)--(A4)--(A3);
\draw (A3) edge["$z$", swap] (A5);
\draw (A4) edge["$z$"] (A6);
\draw (A5) edge["$x$"] (A6);
\coordinate (A7) at (3,0);
\coordinate (A8) at (3,1);
\draw[fill=lightgray, opacity=0.4] (A6)--(A5)--(A7)--(A8)--(A6);
\draw (A5) edge["$w$", swap] (A7);
\draw (A6) edge["$w$"] (A8);
\draw (A7) edge["$x$"] (A8);
\end{tikzpicture}}\\
\end{NiceTabular}
\end{center}
\caption{Restriction-ayclic $C(T)$ with three adjacent squares}
\label{fig:scx-sss}
\end{figure}  
  We thus get, up to isomorphism, two polarizations of $\ijfr $.
They are seen to be:
  \begin{align*}
    J_1& = (x_ax_b, y_ay_b, z_az_b, w_aw_b)
         + (y_az_a,y_aw_a,z_aw_a) \\
     J_2& = (x_ax_b, y_ay_b, z_az_b, w_aw_b)
          + (y_az_b,y_aw_b,z_aw_b)
  \end{align*}
\end{exam}

 \begin{exam} 
  \[ \ijfr = (x^2, y^2, z^2, w^2, t^2) +
    (xz, yw, zt, wx, ty). \]
  
  The simplicial complex $X$ associated to $\ijfr_\Delta $ is given
  in Figure \ref{fig:scx-pentagon}.
\begin{figure}[ht]
\begin{center}  
\begin{tikzpicture}[thick, dot/.style={draw,fill,circle,inner sep=1.5pt},scale=1.5]
\node[dot, label=left:{$w$}] (n1) at (0,0) {};
\node[dot, label=right:{$z$}] (n2) at (1, 0) {};

\node[dot, label=right:{$y$}] (n3) at (1.309,0.951) {};
\node[dot, label=above:{$x$}] (n4) at (0.5, 1.539) {};
\node[dot, label=left:{$t$}] (n5) at (-0.309, 0.951) {};

\draw (n1) edge[] (n2);
\draw (n2) edge[] (n3);
\draw (n3) edge[] (n4);
\draw (n4) edge[] (n5);
\draw (n5) edge[] (n1);
\end{tikzpicture}
\end{center}
\caption{The five-cycle}
\label{fig:scx-pentagon}
\end{figure}
  In a polarization $J(T)$, the cell complex $C(T)$ will then
  have five squares. Up to isomorphism 
  there is only one possible $C(T)$ given in Figure \ref{fig:scx-s5}.
\begin{figure}[ht]
\begin{center}
\begin{tikzpicture}[dot/.style={draw,fill,circle,inner sep=1pt},scale=1.3]
\coordinate (n0) at (1.6087614290087209,0.7933533402912349);
\coordinate (n1) at (-0.737277336810124,0.6755902076156604);
\coordinate (n2) at (1.087155742747658,-0.9961946980917455);
\coordinate (n3) at (-1.722085089822332,0.5019420299487306);
\coordinate (n4) at (0.08715574274765789,-0.9961946980917455);
\coordinate (n5) at (0.608761429008721,0.7933533402912349);
\coordinate (n6) at (-0.8976520102645502,-1.1698428757586754);
\coordinate (n7) at (0.0,0.0);
\coordinate (n8) at (-0.9848077530122081,-0.17364817766692978);
\coordinate (n9) at (1.0,0.0);
\coordinate (n10) at (-0.12851590780140298,1.4689435479068953);
\begin{scope}[thick]
\draw[] (n5) edge["$t$"]       (n0);
\draw[] (n9) edge["$x$", swap] (n0);
\draw[] (n1) edge["$z$", swap] (n3);
\draw[] (n7) edge["$y$"]       (n1);
\draw[] (n1) edge["$x$"]       (n10);
\draw[] (n4) edge["$t$", swap] (n2);
\draw[] (n9) edge["$w$"]       (n2);
\draw[] (n8) edge["$y$"]       (n3);
\draw[] (n4) edge["$z$"]       (n6);
\draw[] (n7) edge["$w$"]       (n4);
\draw[] (n7) edge["$x$"]       (n5);
\draw[] (n5) edge["$y$", swap] (n10);
\draw[] (n8) edge["$w$", swap] (n6);
\draw[] (n7) edge["$z$"]       (n8);
\draw[] (n7) edge["$t$"]       (n9);
\end{scope}
\begin{scope}[on background layer]
\fill[fill=lightgray, fill opacity=0.2] (0.0,0.0)--(0.08715574274765789,-0.9961946980917455)--(1.087155742747658,-0.9961946980917455)--(1.0,0.0)--(0.0,0.0);
\fill[fill=lightgray, fill opacity=0.2] (0.0,0.0)--(-0.9848077530122081,-0.17364817766692978)--(-0.8976520102645502,-1.1698428757586754)--(0.08715574274765789,-0.9961946980917455)--(0.0,0.0);
\fill[fill=lightgray, fill opacity=0.2] (0.0,0.0)--(-0.9848077530122081,-0.17364817766692978)--(-1.722085089822332,0.5019420299487306)--(-0.737277336810124,0.6755902076156604)--(0.0,0.0);
\fill[fill=lightgray, fill opacity=0.2] (0.0,0.0)--(0.608761429008721,0.7933533402912349)--(-0.12851590780140298,1.4689435479068953)--(-0.737277336810124,0.6755902076156604)--(0.0,0.0);
\fill[fill=lightgray, fill opacity=0.2] (0.0,0.0)--(0.608761429008721,0.7933533402912349)--(1.6087614290087209,0.7933533402912349)--(1.0,0.0)--(0.0,0.0);
\end{scope}
\end{tikzpicture}
\end{center}
\caption{Restriction-acyclic $C(T)$ with five squares in a cycle}
\label{fig:scx-s5}
\end{figure}
The corresponding polarization is the standard polarization:
  \begin{align*}
    J& = (x_ax_b, y_ay_b, z_az_b, w_aw_b, t_at_b)
         +  (x_az_a, y_aw_a, z_at_a, w_ax_a, t_ay_a).
  \end{align*}
\end{exam}


In general, if $\ijfr_\Delta$ corresponds to the simplicial complex
which is a cycle graph, then $\ijfr$ has only the standard polarization.
This is a consequence of the following, which is an analog of
\cite[Corollary 2.3]{FO}.

\begin{prop}
  Suppose $\ijfr = \imm_\square + \ijfr_\Delta$, where the simplicial complex
  $X$ corresponding to $\ijfr_\Delta$ is flag, i.e. all generators
  of $\ijfr_\Delta$ are quadratic.
  Then $\ijfr$ has non-standard polarizations if and only if there is a vertex $v$
  on the simplicial complex such that removing $v$ and all its neighbours, the
  induced simplicial complex on the remaining vertices is disconnected.
\end{prop}

\begin{proof} Let the ring be $k[V]$. So the vertices of the simplicial 
complex  is $V$. 
  Let $v \in V$ be a vertex.
  Let $W \subseteq V \backslash \{v\}$ be the vertices which
  are not neighbours of $v$. Then the generators of $\ijfr_\Delta$
  are the $vw$ where $w \in W$ (as $v$ varies).
  We have a non-standard polarization
  involving $v$ iff there is a disjoint union $A \sqcup B = W$ with
  non-empty parts
  such that replacing the generators $vw$ with $v_aw$ for $w \in A$
  and $v_bw$ for $w \in B$, gives an ideal $\ijfr^\prime$ which is a separation of
  $\ijfr$. Pick then  $w^a \in A$ and $w^b \in B$. Then $(v_a - v_b)w^a w^b$
  is zero in the quotient ring $S^\prime/\ijfr^\prime$. As $(v_a - v_b)$ is
  a non-zero divisor, this gives that $w^a w^b \in \ijfr^\prime$, and so
  $w^aw^b \in \ijfr$. But as $w^a \in A$ and $w^b \in B$ were chosen
  arbitrarily, this means that there is no edge between the
  vertices of $A$ and of $B$. Thus we have a separation of $\ijfr$
  iff such a situation occurs. This separation may then further be
  standard polarized. The result will be a non-standard polarizaton,
  as a variable has been separated in the mixed term generators of $\ijfr$. 
\end{proof}

\section{Liaison for squarefree monomial ideals}
\label{sec:link}
Liaison is a standard and well-developed technique in algebraic
geometry, \cite{Mig,HU,PS}. Here we apply it to squarefree monomial
ideals where it takes a simple form.
Our objective is to show that if $T \sus A_1 \times \cdots \times A_n$
is acyclic, then the complement subset $T^c$ is also acyclic.
Our argument for this is entirely algebraic. One could envision
a more topological argument using Alexander duality, but
it is a challenge to construct a suitable retract onto
$C(T)$ or $C(T^c)$. 
In general it would be interesting to understand more about
how the homology of the free cellular complexes $\Ctalg(T)$ and $\Ctalg(T^c)$
are related.

\subsection{Liaison}
We now consider linkage of squarefree monomial ideals with respect to
the complete intersection ideal
\begin{equation} \label{eq:link-ci}
  C = (m(A_1), m(A_2), \ldots, m(A_n)).
\end{equation}
Two ideals $J$ and $K$ containing $C$ are {\it linked} if
their heights (codimensions) are all equal (so is $n$), and
$K = (C : J)$ and $J = (C : K)$.

For simplicial complexes, we have the following
concrete, simple, and general result.
We thank an anonymous referee for the below statement.

\begin{prop} \label{pro:link-xyz}
  Let $X,Y$ be simplicial sub-complexes of a simplicial complex $Z$.
  Then the ideal
  $I_Y = (I_Z : I_X)$ if and only if $Y$ is the simplicial complex
  generated by $Z \backslash X$.

  In particular, under this assumption, for $F$ a facet of $Z$:
  \begin{itemize}
  \item[i.] If $F$ is not in $X$, it is a facet of $Y$,
    \item[ii.] If $F$ is in $X$, it is not a facet of $Y$, 
    \item[iii.] The set of facets of $Y$ are contained in the set of
      facets of $Z$.
\end{itemize}
\end{prop}

\begin{proof} Let $Y^\prime$ be the simplicial complex generated by
  $Z \backslash X$, and $I_Y = (I_Z : I_X)$.
  We show firstly $Y \sus Y^\prime$ and
  secondly $Y^\prime \sus Y$.
  If $F$ is a face of $Y$, the monomial $m(F) \not \in I_Y$.
  Then there is $m(R) \in I_X$ (equivalently $R \not \in X$),
  such that $m(F) \cdot m(R) \not \in I_Z$.
  This monomial has $m(F \cup R)$ as a factor so $F \cup R$ is in $Z$
  and not in $X$. The upshot is that $Y \sus Y^\prime$.

  Conversely let $G$ be a facet of $Y^\prime$. Then $m(G) \in I_X$
  and $m(G) \not \in I_Z$. If $m(G) \in I_Y$, then $m(G) \cdot m(G) \in I_Z$.
  But since $I_Z$ is squarefree, we have $m(G) \in I_Z$, against assumption.
  So $m(G) \not \in I_Y$, giving $G \in Y$. Hence $Y^\prime \sus Y$.

  The three consequences are immediate.
\end{proof}

\begin{coro}
  We have $I_Y = (I_Z : I_X)$ and $I_X = (I_Z : I_Y)$ iff we have
  a disjoint union of facets
   \[ \text{facets}(Z) = \text{facets}(X) \sqcup \text{facets}(Y). \]
  \end{coro}
 
In our specific situation we have:
  
\begin{prop} \label{pro:link-C}
  Let $J$ and $K$ be squarefree monomial ideals
  containing the complete intersection $C$ of \eqref{eq:link-ci},
  and let $X$ and $Y$ be the corresponding simplicial complexes
  via the SR-correspondence
  (which are subcomplexes of $\Delst(\Ab)$).
  Then $J$ and $K$ are linked by $C$ iff the facets of $\Delst(\Ab)$ is the
 disjoint union of the facets of $X$ and $Y$. 
\end{prop}

The following is Claim 1 in Lemma 3.6 in \cite{AN}.

  \begin{coro} \label{cor:link-JK}
   Let  $T \sus  A_1 \times \cdots \times A_n$ and let $J$ be the Alexander dual
    $I(T)^\vee$. Let $K$ be a squarefree monomial ideal containing $C$.
    Then $J$ and $K$ are linked by $C$ iff $K$ is the Alexander dual
    $I(T^c)^\vee$.
  \end{coro}

  \begin{proof}
        If $J$ and $K$ are linked by $C$, the facets of  $\Delst(\Ab)$ is the
        disjoint union of the facets of $X$ and $Y$.
        The facets of $X$ are then $[\bfa]^c$ for $\bfa \in T$, and
        so the facets of $Y$ must then be the $[\bfa]^c$ for $\bfa \in T^c$.
        Hence $K$ is the Alexander dual of $I(T^c)$.

        Conversely if $K$ is the Alexander dual of $I(T^c)$,
        the facets of $Y$ and $X$ are
        the $[\bfa]^c$ for respectively $\bfa \in T^c$ and $\bfa \in T$.
        So $\Delst(\Ab)$ is the disjoint union of these facets, and so
        $J$ and $K$ are linked.
      \end{proof}
     
We have a diagram:
 \[
   \xymatrix{  I(T) \ar@{<->}_{\text {complement rainbow ideals}}[d]
     \ar@{<->}^{\text{Alexander dual}}[rr] & & 
     J(T) \ar@{<->}^{{{\text{\rotatebox{0}{{linked}}}}}}[d] \\
     I(T^c) \ar@{<->}^{\text{Alexander dual}}[rr] &  & K(T)  
      }
 \]

 \begin{coro}
  If the rainbow monomial ideal $I(T)$ has a linear resolution,
  the complement rainbow ideal $I(T^c)$ also has a linear resolution.
  Hence $T$ is acyclic iff the complement $T^c$ is acyclic.
\end{coro}

\begin{proof}
Let $J$ be the Alexander dual of $I(T)$.
It is Cohen-Macaulay since $I(T)$ has a linear resolution.
The linked ideal $K$ of $J$ is also Cohen-Macaulay by \cite[Thm.21.23]{Eis}.
The Alexander dual of $K$ is $I(T^c)$ and so this ideal also has
linear resolution. 
\end{proof}

\begin{coro}
  If $X$ is a full-dimensional triangulated ball on $\Delst(\Ab)$, its
  full-dimensional complement
  $Y = X(T^c)$ on $\Delst(\Ab)$ is also a triangulated ball.
\end{coro}

As stated after Theorem \ref{thm:rain-const}, a full-dimensional CM
subcomplex of $\Delst(\Ab)$ may not be shellable. 
 H.T. Hall \cite{Hall} gives a full-dimensional
 subcomplex $Y$  of a $12$-dimensional hyperoctahedron $H$
 which is shellable, but
 for no facet $F$ outside of $Y$ is $Y \cup \overline{F}$ shellable,
 where $\overline{F}$ is the simplicial complex generated by $F$.
 (One says $H$ is not extendably shellable.)
 The consequence of this fact is the following.

 \begin{prop} \label{pro:link-hall}
  The complement sub-complex $X$ of $Y$, generated by $H \backslash Y$,
  is a full-dimensional CM subcomplex of $H$ which is {\it not shellable}.
  (But by Theorem \ref{thm:rain-const} $X$ {\it is
    constructible}.)
\end{prop}

\begin{proof} By Proposition \ref{pro:link-C} 
  $X$ and $Y$ are linked. Since $Y$ is CM, 
 $X$ is also CM by \cite[Thm.21.23]{Eis}.
  
 If $F_1, \ldots, F_m$ is a shelling of $X$, the simplicial
 complex generated by $X \backslash \overline{F_m}$
 is shellable by $F_1, \ldots, F_{m-1}$ and so is CM.
 Its complement in turn is $Y \cup \overline{F_m}$, which is the linked complex
 and so would also be
 CM. By Lemma \ref{lem:pre-CM}, $Y \cap \overline{F_m}$ would be CM of dimension
 one less than $\overline{F_m}$, and so a pure codimension one
 subcomplex of $\overline{F_m}$. This would
 imply shellability of $Y \cup \overline{F_m}$, contradicting the
 result of \cite{Hall}.
 \end{proof}

\section{The boundary}
\label{sec:bound}

For an acyclic $T \sus A_1 \times \cdots \times A_n$, we have the
partition into two acyclic sets $T \sqcup T^c = A_1 \times \cdots \times A_n$.
Let $X(T)$ and
$X(T^c)$ be the corresponding simplicial balls, whose union $X(T) \cup X(T^c)$
is $\Delst(\Ab)$. Recall that $J$ is the SR-ideal of $X(T)$ and
$K$ is the SR-ideal of $X(T^c)$ (formerly mostly denoted $Y$). 
We also denote $L(T) = I(T^c)$. 

\begin{prop} \label{pro:bound-bound} For acyclic $T$,
  the boundaries of $X(T)$ and $X(T^c)$ are the same:
  \[\pd X(T) = \pd X(T^c) = X(T) \cap X(T^c). \]
  So this is a simplicial sphere with defining ideal $J + K$,
  which is a Gorenstein ideal.
  The Alexander dual of $J+K$ is $I(T) \cap L(T)$, an
  ideal generated in degree $n+1$. 
\end{prop}

\begin{proof}
  Every codimension one
  face in $\Delst(\Ab)$ is on two facets. 
  The boundary of $X(T)$ consists of those codimension one faces
  on only one facet of $X(T)$. But then the other facet in $\Delst(\Ab)$
  containing this codimension one face must be in $X(T^c)$.
  Conversely, let $F$ be a face in $X(T) \cap X(T^c)$. Then there
  are facets $[\bfa]^c$ in $X(T)$ and $[\bfb]^c$ in $X(T^c)$ containing
  $F$. If $\bfa$ and $\bfb$ have $r$ distinct coordinates, consider
  a path from $\bfa$ to $\bfb$ of length $r$ changing one coordinate
  at a time. Then every
  $\bfa_i$ on this path will have $F \sus [\bfa_i]^c$. At some
  stage $\bfa_i \in T$ and $\bfa_{i+1} \in T^c$. Their intersection
  is a codimension one face i) on the boundaries of $X(T^c)$ and $X(T)$, and
  ii) containing $F$. This shows the statement about the boundaries.

  That $J + K$ is the defining ideal is clear
  since this is the ideal of  $X(T) \cap X(T^c)$. Since $X(T)$ and  $X(T^c)$
  are simplicial balls, their boundaries are simplicial spheres
  (which have dimension one less that $X(T)$ and $X(T^c)$). Thus
  its defining ideal is Gorenstein \cite[Corollary 5.6.5]{BH93}.

  The complements of the facets of $X(T) \cap X(T^c)$, give the generators
  of $I(T)  \cap L(T)$ and they must then have degree $n+1$
  (one more than the generators of $I(T)$).
  \end{proof}

  \begin{rema}
    For ideals $J$ and $K$ linked by a complete
    intersection $C$, the sum $J + K$ may
    not be a Gorenstein ideal. But this holds if $J \cap K = C$
    \cite[Thm.2.8]{PS}. This
    is the situation here, since the Alexander dual of $J \cap K$
    is $I(T) + L(T)$, which is the ideal $B$ generated by all rainbow monomials.
    Since $C$ is the Alexander dual of $B$  this follows.
    \end{rema}

    \begin{rema} In general for $Z$ a $d$-dimensional simplicial sphere, and
      $X$ a $d$-dimensional Cohen-Macaulay sub-complex, the complement
      simplicial complex $Y$ generated by $Z \backslash X$ is also
      Cohen-Macaulay, and their boundaries $\pd X = \pd Y = X \cap Y$.
      The latter may be shown using that Cohen-Macaulay simplicial
      complexes are gallery connected (strongly connected),
      \cite[Prop. 9.1.12]{HeHi}.
      \end{rema}
    
    \begin{rema}  A concrete, effective description of the
      SR-ideal of the boundary is given in \cite[Thm.5.2]{DFN}
      for simplicial spheres defined by letterplace ideals.
    \end{rema}

\section{Artin monomial ideals from a box}
\label{sec:Artin}

Subsection \ref{subsec:pol} considered polarizations $J$ of Artin monomial
ideals $\ijfr \sus k[x_1, \ldots, x_n]$ containing the complete intersection
\begin{equation} \label{eq:Artin-ic} \ic = (x_1^{\alpha_1}, \ldots, x_n^{\alpha_n})
\end{equation}
Section \ref{sec:link} considered linkage of two squarefree monomial ideals
$J$ and $K$ in $k[A]$
with respect to the complete intersection
\[ C = (m(A_1), m(A_2), \ldots, m(A_n)). \]
We show this corresponds to linkage downstairs of Artin monomial ideals.

\medskip
For an integer $m \geq 0$ let $[m]_0$ denote the chain
$\{ 0 < 1 <  \cdots < m\}.$
For $i = 1, \ldots, n$ let  $A_i$ be a set with $\alpha_i \geq 1$ elements.

Recall that an up-set $U$ in a poset $P$ is a set such that $x \in U$
and $x \leq y$ implies $y \in U$.  Similarly we have the notion of
down-set.
Consider an up-set $U$ in the product poset
\[ \bfal = [\alpha_1-1]_0 \times [\alpha_2-1]_0 \times
  \cdots \times [\alpha_n-1]_0,\]
We get a monomial ideal $\iifr(U) = (x^\bfb \, |\, \bfb \in U )$.
This is 
an {\it $(\alpha_1 -1, \ldots, \alpha_n - 1)$-determined}
Artin monomial ideal (a notion
introduced in \cite{Mil-ad2}, see also \cite[Sec. 5.20]{Mi-St}).
Let $\iifr \sus k[x_1, \ldots, x_n]$
be the Artin monomial ideal $\ic+\iifr(U)$.

Let $D$ be the complement of $U$ in $\bfal$, so $D$ is a down-set in $\bfal$.
Let $\bfal^\opos$ be the opposite poset of $\bfal$. There is an
isomorphism of posets $\bfal^\opos \iso \bfal$ sending:
\[ \bfi^\opos = (i_1^\opos, i_2^\opos, \ldots, i_n^\opos)
  \mapsto (\alpha_1-1-i_1, \alpha_2 - 1-i_2, \ldots, \alpha_n-1-i_n) =:
  \overline{\bfi}. \]
This sends the the up-set $D^\opos$ in $\bfal^{\opos}$  to an up-set $\oD$
in $\bfal$. 
Let $\ik \sus k[x_1, \ldots, x_n]$
be the Artin ideal $\ic + \iifr(\oD)$, the {\it dual} of $\iifr$
with respect to $\ic$.

\begin{rema} \label{rem:artin-ad}
  The ideal $\iifr(\oD)$ is the {\it Alexander dual of $\iifr(U)$ with respect
    to $(\alpha_1 -1, \alpha_2-1, \ldots, \alpha_n-1)$},
  see \cite[Def. 5.20]{Mi-St} or \cite{Mil-ad}. The definition
in \cite{Mi-St}
  is different
  from the one above (and more general in scope), but Lemma \ref{lem:Artin-link}
  below and \cite[Cor. 5.25]{Mi-St} show they coincide in our case.
\end{rema}

\begin {exam} \label{ex:Artin-ud}
  Let $n = 3$ and $\alpha_1, \alpha_2, \alpha_3 = 2$, so
  $\bfal = \{0,1\}^3$.
  Let the up-set of $\bfal$ be:
  \[ U : 111, 101, 011. \]
  It gives an Artin monomial ideal:
  \[ \iifr = \iifr(U) + \ic =  (x_1^2, x_1x_3, x_2x_3, x_2^2, x_3^2). \]
  The corresponding down-set is:
\[ D : 000, 100, 010, 001, 110,  \]
and so
\[ \oD : 111, 011, 101, 110, 001. \]
The dual Artin monomial of $\iifr$ with respect to $(x_1^2, x_2^2, x_3^2)$
is:
\[ \ik = \iifr(\oD) + \ic = (x_1^2, x_1x_2, x_2^2, x_3). \]
\end {exam}

\begin{lemm} \label{lem:Artin-link}
  The ideals $\iifr$ and $\ik$ are linked with respect to the complete
  intersection $\ic$.
\end{lemm}

\begin{proof}
  Note that $(b_1, \ldots, b_n) \in \alpha$ is in $\oD$ iff
  $(\alpha_1 - 1- b_1, \ldots, \alpha_n -1 - b_n)$ is in $D$, i.e.
  is not in $U$.
  That is, for {\it every} $(a_1, \ldots, a_n)$ in $U$
  there is some $j$ with
  $\alpha_j - 1- b_j < a_j$, or equivalently $a_j + b_j \geq \al_j$.
  
  Now let us show that $\ik = \ic + \iifr(\oD)$ equals
  $(\ic:\iifr)=(\ic:\iifr(U))$.
  That  $x_1^{b_1}\cdots x_n^{b_n}\in \iifr(\oD)$ is by the above 
 equivalent to for {\it every} $x_1^{a_1} \cdots x_n^{a_n}$ in $\iifr(U)$
  there is some $j$ with $a_j + b_j \geq \al_j$.
  But this is equivalent to $x_1^{b_1} \cdots x_n^{b_n} \cdot \iifr(U)$ being
  in $\ic$.
  \end{proof}

 \begin{prop} \label{pro:pol-commute}
     Let $J$ be a polarization of $\ijfr$. Let $K$ be
the link of $J$ by $C$, and $\ik$ the link of $\ijfr$ by $\ic$,
so $\ik$ is the dual of $\ijfr$. Then $K$ is a polarization of $\ik$.
\end{prop}

\begin{proof}
  Set ${A_l}=\{a^l_{1},\ldots,a^l_{\alpha_l}\}$ for every $l=1,\ldots,n$.
We divide the quotient rings $k[A]/J, k[A]/K$ and  $k[A]/C$ by all the variable
differences $a^l_{j}-a^l_{j-1}$ for $j=2,\ldots,\alpha_l$ and $l = 1, \ldots, n$.
We get quotient rings $k[x_1,\ldots,x_n]/\ijfr$, $k[x_1,\ldots,x_n]/\ik^\prime$
and $k[x_1,\ldots,x_n]/\ic$, where $\ijfr$, $\ik^\prime$ and $\ic$ are Artin.
By \cite[Proposition 2.8]{HU} when we cut down by a regular sequence for
$k[A]/C$ the ideals we get are still linked, and it is a regular sequence
also for $k[A]/J$ and $k[A]/K$.
Thus $\ik^\prime$ is linked to $\ijfr$ by $\ic$, so so $\ik^\prime = \ik$ and
$K$ is a polarization of this.
\end{proof}

 To sum up we have a diagram:
 \[
   \xymatrix{ & I(T) \ar@{<->}_{\text{complement rainbow ideals} \hskip 3mm}[ld] \ar@{<->}^{\text{Alexander dual}}[rr] & & 
     J(T) \ar@{<->}_{{{\text{\rotatebox{36}{\hskip -7mm{linked}}}}}}[ld] \ar^{\text{polarization}}[d]\\
     L(T) \ar@{<->}^{\text{Alexander dual}}[rr] & & K(T) \ar_{\text{polarization}}[d] &
     \ijfr. \ar@{<->}^{\text{\rotatebox{36}{\hskip -7mm{  linked}}}}[ld] \\
     & & \ik & }
 \]
 The left-right arrows are correspondences, and the vertical arrows 
 are maps of ideals. But note that the process of polarization
 goes from the lower ideal to the upper ideal.

\section{Squeezed balls and spheres}
\label{sec:squeezed}
This section is somewhat independent of the rest of the article.
We show that the squeezed balls
constructed by G.Kalai \cite{Ka} (their boundaries are squeezed spheres),
have Stanley-Reisner ideals which are
regular quotients by variable differences
of a special class of letterplace ideals.
Letterplace ideals are polarizations of Artin monomial ideals
by Theorem 5.12 and Corollary 2.4 in \cite{FGH}.

In all, this shows that the two classical ways of constructing large numbers
of triangulated balls and spheres, those of Bier balls and spheres \cite{BZPS},
\cite[Cor.5.3]{DFN} and 
squeezed balls and spheres \cite{Ka},
both are part of the setting we give here.
They are actually at the opposite ends, as Bier balls arise from
letterplace ideals of discrete posets (antichains) \cite[Cor.5.3.]{DFN},
  while squeezed balls arise from letterplace ideals of totally ordered
  posets (chains), see Subsection \ref{subsec:app-lp}.

\subsection{Squeezed balls}

Given an integer $e \geq 1$. Consider sequences of natural numbers
$\bfi : 1 \leq i_1 < i_2 < \cdots < i_e$ such that
each difference $i_p - i_{p-1} \geq 2$. Let $F_e$ be the set of
such sequences. There is a partial order on $F_e$ by $\bfi \leq \bfj$
if each $i_p \leq j_p$. Let $F_e(m)$ be such sequences with
$i_e < m$.

  Each such sequence induces a subset of natural numbers of cardinality $2e$:
  \[ \bfio : i_1, i_1+1, i_2, i_2+1, \ldots, i_e, i_e + 1. \]
  Let $I$ be a {\it finite}
  down-set (order ideal) in $F_e$ for its partial order.
  Kalai shows \cite[Sec.4]{Ka} that the simplicial complex $B(I)$ with facets
  $\bfio$ for $ \bfi \in I$,
  is a simplicial ball, termed a {\it squeezed ball}, with boundary 
  a {\it squeezed sphere}. By a counting
  argument, there are so many of these that most of these balls and
  spheres do not have
  convex realization.

\subsection{Poset homomorphisms}
Let $P$ and $Q$ be posets. We consider order-preserving maps $\phi : P \pil Q$,
and denote by $\Hom(P,Q)$ the set of such.
This becomes itself a poset (this is an internal Hom), by
$\phi \preceq \psi$ if
$\phi(p) \leq \psi(p)$ for every $p \in P$.
There are order-preserving maps:
\[ P \times \Hom(P,Q) \pil P \times Q  \quad \text{inducing} \quad
  \Hom(P,Q) \mto{G} \Hom(P, P \times Q), \]
and for $\phi : P \pil Q$ the image of the map $G(\phi)$ is the graph
$\Gamma(\phi) \sus P \times Q$.

Denote by $[n]$ the chain $\{1 < 2 < \cdots < n\}$, and $[n]_0$ the chain
$[n] \cup \{0 \}$. Note that $[n]_0 \iso \Hom([n],\{0 < 1\})$. 
We have a composition:
\begin{align}\label{eq:app-alpha} [n] \times \Hom([n],[e]_0) \pil
  [n] \times [e]_0 & \mto{\alpha} [n+2e]   \\ \notag
  (i,j) & \mapsto i + 2j
\end{align}
giving
\[ \Phi : \Hom([n],[e]_0) \pil \Hom([n],[n+2e]). \]
We also have a composition:
\begin{align*} [e] \times \Hom([e],[n]_0) \pil
  [e] \times [n]_0 & \mto{\beta} [n+2e-1] \\
  (j,i) & \mapsto (2j-1) + i
\end{align*}
giving
\[ \Psi : \Hom([e],[n]_0) \pil \Hom([e],[n+2e-1]). \]

Note that the domains of $\Phi$ and $\Psi$ identify:
\begin{align*} \Hom([n],[e]_0)= & \Hom([n],\Hom([e],\{0<1\})) \iso
                                  \Hom([n]\times [e], \{0 < 1 \}) \\
  \iso & \Hom([e], \Hom([n],\{0<1\}))
         = \Hom([e],[n]_0).
\end{align*}
This is just mapping a partition to its dual partition. Given
$\phi : [n] \pil [e]_0$, it may be characterized by a sequence 
\[ 0 = a_0 \leq a_1 \leq \cdots \leq a_e \leq a_{e+1} = n, \]
such that on the interval $[a_{p}+1 , a_{p+1}]$, the value
of $\phi$ is $p$. These $a$'s give the dual map $\phi^\prime : [e] \pil [n]_0$.

\begin{lemm}
  \label{lem:app-iso} \hskip 1mm

  \begin{itemize}
  \item[a.] The image of $\Psi$ is $F_e(n+2e)$ and this gives an
    isomorphism of posets between $\Hom([e],[n]_0)$ and $F_e(n+2e)$.
  \item[b.] Let $\bfi = \im \Psi(\phi^\prime)$.
    Then $\bfio = (\im \Phi(\phi))^c$.
    \end{itemize}

  \end{lemm}
  
\begin{proof}
  a. Let  $\phi \in \Hom([n],[e]_0)$ and $\phi^\prime \in \Hom([e],[n]_0)$
  correspond, and let $a_p = \phi^\prime(p)$.
The image of the map $\Psi(\phi^\prime)$ is $\bfi$ where $i_p = a_p + 2p-1$
for $p = 1, \ldots, e$.
This gives a bijection between
$\Hom([e],[n]_0)$ and $F_e(n+2e)$, and it is clear
that the partial orders correspond.

b. 
The gaps in the image of $\Phi(\phi)$ in $[n+2e]$ are
as follows: If
\[ a_{p-1} < a_p = \cdots =  a_q < a_{q+1}, \]
there is a gap ranging from
$(a_p + 2(p-1))+1$ to $(a_q +1 +2q)-1$. This
gap has even size $2q-2p+2$. The values for $\bfi$ are
\begin{align*} i_p= a_p + 2p-1, \quad & i_{p+1} = a_{p+1} + 2(p+1)-1 = a_p + 2p + 1,\\
  \ldots  \quad &  i_q = a_q + 2q-1 = a_p + 2q-1,
\end{align*}
with each step in this range of size two. 
This shows the last claim.
  \end{proof}

\subsection{Letterplace ideals} \label{subsec:app-lp}
For a finite set $R$ let $k[x_R]$ be the polynomial ring with variables
$x_r, r \in R$. If $S \sus R$ let $x_S$ be the monomial $\prod_{s \in S}x_s$.

Let $P$ be a finite poset, and
$\cI$ a down-set in (the internal Hom) $\Hom(P,[m])$.
Each $\phi \in \cI$ gives a graph $\Gamma(\phi) \sus P \times [m]$.
The associated
{\it co-letterplace ideal}
$L(P,[m];\cI)$ is the ideal in $k[x_{P \times [m]}]$ generated
by monomials $x_{\Gamma \phi}$ for $\phi \in \cI$.
The associated {\it letterplace ideal} $L(\cI;P,[m])$ is the Alexander
dual of $L(P,[m];\cI)$. The letterplace
ideal is a polarization of an Artin monomial ideal by Corollary 2.4 and
Theorem 5.12 in \cite{FGH}.
Its associated triangulated ball is described in \cite{DFN}. 

Now consider a down-set $\cI \sus \Hom([n],[e]_0)$.
It gives a co-letterplace ideal $L([n],[e]_0;\cI)$ and a letterplace
ideal $L(\cI;[n],[e]_0)$. The co-letterplace ideal is
a rainbow ideal with $n$ color classes each of size $e+1$.
The map $\alpha : [n] \times [e]_0 \pil [n+2e]$ from
\eqref{eq:app-alpha} is given by
$(i,j) \mapsto i+2j$. It gives a map of polynomial rings
$\tilde{\alpha} : k[x_{[n] \times [e]_0}] \pil k[x_{[n+2e]}]$.
We may think of this as collapsing each fiber of $\alpha$ into
a single variable. If a fiber has $t$ elements,
we are dividing by $t-1$ variable differences, and we run through all fibers.
The variable differences are of type $x_{i,j+1} - x_{i+2,j}$, so
$x_{i,j+1}$ and $x_{i+2,j}$ are in different color classes.
Let $L^\alpha(\cI;[n],[e]_0)$ be the image ideal of $L(\cI;[n],[e]_0)$
by $\tilde{\alpha}$, an ideal in the polynomial ring $k[x_{[n+2e]}]$.

  \begin{theo}
    Let the poset ideal $\cI \sus \Hom([n],[e]_0)$ correspond to the
    poset ideal $I \sus F_e(n+2e)$ by the isomorphism $\Psi$.
    The Stanley-Reisner ideal of the squeezed ball $B(I)$ is the
    ideal $L^\alpha(\cI;[n],[e]_0)$.
  \end{theo}

  \begin{proof}
    Each fiber of $\alpha$
    is seen to be a {\it bistrict} chain in $[n] \times [e]_0^\opos$. 
That is, if $(p,u)$ and $(q,v)$ are in the same fiber of $\alpha$,
either $p > q$ and $u < v$, {\it or} $p < q$ and $u > v$.
By Theorem 7.3 in \cite{FGH}, the images by $\tilde{\alpha}$ 
\[ L^\alpha(\cI;[n],[e]_0) \text{ of } L(\cI;[n],[e]_0), \quad
  L^\alpha([n],[e]_0;\cI) \text{ of } L([n],[e]_0;\cI) \]
are still square-free and they are Alexander dual ideals.

The facets of the simplicial complex with Stanley-Reisner ideal
$L^\alpha(\cI;[n],[e]_0)$ are the complement in $[n+2e]$ of the
indices of the minimal generators of $L^\alpha([n],[e]_0;\cI)$.
These index sets are $\im \Phi(\phi)$ for $\phi \in \cI$ and by 
Lemma \ref{lem:app-iso} part b their complements give
the facets of the squeezed ball $B(I)$.
\end{proof}

\section{Further problems} \label{sec:further}
We sketch three open research areas.

\subsection{Acyclic subsets $T \sus \bfA$}
\label{subsec:further-acyclic}
Section \ref{sec:hypercube} illustrates well how one may make
acyclic subsets on a hypercube. When the cardinalities
of the $A_i$ become $\geq 3$, things get more complicated.
\Ign{
The definition of a acyclic   set involves an exchange condition
which has a similarity to the usual way of defining a matroid.
Acyclic sets seem a basic combinatorial object. 

In the case $n  = 2$, so $T \sus A_1 \times A_2$, the acyclic   $T$
are well understood. They are a class of bipartite graphs called
Ferrer's graphs and are studied in \cite{CN}.
They arise as follows: Put a total order on $A_1$ and
a total order on $A_2$. Then $T$ is a down-set for the product order on
$A_1 \times A_2$. (this is an internal Hom) by $\phi \preceq \psi$ if
$\phi(p) \leq \psi(p)$ for every $p \in P$.}

\begin{problem}
 Give machines for constructing acyclic sets $T \sus \bfA$.
\end{problem}

One such machine is the following.

\medskip
\noindent {\bf Co-letterplace ideals.} This is from \cite{FGH}, which
is again inspired by \cite{EHM}.
Recall from Subsection \ref{subsec:app-lp} that for a poset $P$
an order-preserving map $\phi : P \pil [n]$ gives a graph
\[ \Gamma \phi = \{(p, \phi(p)) \, | \, p \in P \} \sus P \times [n]. \] 
The graph may be considered as an element of $[n]^P = \times_{p \in P} [n]$,
the cartesian product of $[n]$'s, one factor for each $p \in P$.
For $\cI$ a down-set in $\Hom(P, [n])$,
let
\[ T = \{ \Gamma \phi \, | \, \phi \in \cI \} \sus [n]^P. \]
The ideal $I(T)$ is the co-letterplace ideal from \cite{FGH}, and has
linear resolution, in fact it has linear quotients.
In \cite{DFN} the simplicial balls $X(T)$ are studied.
Their boundaries, simplicial spheres, have a particularly nice and
compact description, \cite[Sec.5]{DFN}. The classical cases are:
\begin{itemize}
\item Bier balls and spheres \cite{Bier, BZPS}: $P$ is an antichain and
  $n = 2$, and in \cite{Mur11} for any $n$.
\item Squeezed balls and spheres \cite{Ka}: $P$ is a chain. Or rather they are
  derived from this case, see Section \ref{sec:squeezed}.
  \end{itemize}
  D.Lu and Z.Wang, \cite{LW}, extend this. They replace the down-set
  $\cI$ of $\Hom(P, [n])$ with larger classes of
  subsets of $\Hom(P,[n])$.

\medskip
Polarizations of powers of maximal graded ideals $(x_1, \ldots, x_n)^m$
are characterized in \cite{AFL}. When $n = 3$, there is a nice visual
procedure for this. 

\Ign{
For $n=3$ or $m=2$ this may be translated
to simple constructive procedures.

\medskip
\noindent {\bf 2. Trees $m=2$.} \label{rem:res-tree}
Polarizations of $(x_1, x_2, \ldots, x_n)^2$ are (up to isomorphism)
in bijection with trees with $n$ edges. As an example, when $n=3$ 
there are two such trees:
\begin{center}
  \begin{tikzpicture}[scale=0.7]
\draw[thick] (0,0)--(0,-0.8) node[anchor=west] at (0,0.7){};
\draw[thick] (0,0)--(0.75,0.5) node[anchor=north] at (0.5,-0.33)  {};
\draw[thick] (0,0)--(-0.75,0.5) node[anchor=south] at (-0.5,-0.33)  {};
\filldraw (0,0) circle (2pt);
\filldraw (0,-0.8) circle (2pt);
\filldraw (0.75,0.5) circle (2pt);
\filldraw (-0.75,0.5) circle (2pt);

\draw[thick] (3,0)--(4,0) node[anchor=north] at (3.5,0){};
\draw[thick] (4,0)--(5,0) node[anchor=north] at (4.5,0)  {};
\draw[thick] (5,0)--(6,0) node[anchor=north] at (5.5,0)  {};
\filldraw (3,0) circle (2pt);
\filldraw (4,0) circle (2pt);
\filldraw (5,0) circle (2pt);
\filldraw (6,0) circle (2pt) ;
\end{tikzpicture}
\end{center}
giving the two different polarizations of Example \ref{ex:intro-m2}.
  The left tree corresponds to the left hexagon triangulation
  in Figure \ref{fig:HexPol}. The vertices of the trees correspond to the
  triangles of the hexagons.

Given a tree $t$ with edges and vertices  $(E,V)$.
Let $I \sus E \times V$ be the incidence relation between edges and
vertices. We may non-canonically identify $I$ and $E \times \{0,1\}$.
Each vertex $v \in V$ gives a natural section
$E \mto{s_v} I \iso E \times \{0,1\}$
of the projection to $E$, by $e \mapsto (e,w)$ where $w$ is the vertex of $e$
such that $e$ is on the path from $w$ to $v$.
This gives an element $\im p_2 \circ s_v$  of
$\{0,1\}^E$. As $v$ varies, these elements give a subset $T$ which is
acyclic.
}

\medskip
\noindent {\bf Planar triangulations $n=3$.}  The generators of $(x_1,x_2,x_3)^m$ corresponds
to the set $\Delta_3(m)$ consisting of all triples $(a_1,a_2,a_3) \in
{\mathbb N}_0^3$ with $a_1 + a_2 + a_3=m$. Figure \ref{fig:problem-down}
displays in the
case $m=3$ these ten triples as vertices.
The figure has three ``down-triangles''.
Polarizations of $(x_1,x_2,x_3)^m$ are in bijection with removals
of either zero or one edge from each down-triangle.
\begin{figure}
\begin{tikzpicture}[scale=0.5]
\draw (-3,0)--(3,0);
\draw (-2,1.6)--(2,1.6);
\draw (-1,3.2)--(1,3.2);
\draw (-3,0)--(0,4.8);
\draw (3,0)--(0,4.8) ;
\draw (-1,0)--(1,3.2);
\draw (1,0)--(-1,3.2);
\draw (-1,0)--(-2,1.6);
\draw (1,0)--(2,1.6);
\filldraw[black] (0,4.8) circle (2pt)  node[anchor=south] at (0,4.9){(3,0,0)};
\filldraw[black] (-3,0) circle (2pt)  node[anchor=east] at (-3.1,0){(0,3,0)};
\filldraw[black] (3,0) circle (2pt)  node[anchor=west] at (3.1,0){(0,0,3)};
\filldraw[fill opacity=0.1] (-2,1.6)--(0,1.6)--(-1,0)--(-2,1.6);
\filldraw[fill opacity=0.1] (-1,3.2)--(1,3.2)--(0,1.6)--(-1,3.2);
\filldraw[fill opacity=0.1] (2,1.6)--(0,1.6)--(1,0)--(2,1.6);
\end{tikzpicture}
  
\caption{Vertices correspond to generators of $(x_1,x_2,x_3)^3$}
\label{fig:problem-down}
\end{figure}

\subsection{Regular quotients of rainbow ideals}
\label{subsec:further-regquot}
If the simplicial ball $X(T)$ is a one-dimensional ball,
then $\Delst(\Ab)$ must be one-dimensional. This occurs only for
$n = 2$ and $|A_1|= |A_2| = 2$, when $\Delst(\Ab)$ is a square (perimeter),
or $n=1$ and $|A_1|=3$, when it is a triangle (perimeter).
Hence the only one-dimensional balls $X(T)$ are:

\begin{center}
\begin{tikzpicture}[dot/.style={draw,fill,circle,inner sep=1pt},scale=0.9]
\fill (0,0)  circle (.08);
\fill (1,0)  circle (.08);

\draw[very thick] (0,0)--(1,0);
\end{tikzpicture}
\quad
\begin{tikzpicture}[dot/.style={draw,fill,circle,inner sep=1pt},scale=0.9]
\fill (0,0)  circle (.08);
\fill (1,0)  circle (.08);
\fill (2,0)  circle (.08);

\draw[very thick] (0,0)--(2,0);
\end{tikzpicture}
\quad
\begin{tikzpicture}[dot/.style={draw,fill,circle,inner sep=1pt},scale=0.9]
\fill (0,0)  circle (.08);
\fill (1,0)  circle (.08);
\fill (2,0)  circle (.08);
\fill (3,0)  circle (.08);
\draw[very thick] (0,0)--(3,0);
\end{tikzpicture}
\end{center}

In particular the line graphs with $\geq 5$ vertices:

\begin{center}
\begin{tikzpicture}[dot/.style={draw,fill,circle,inner sep=1pt},scale=0.9]
\fill (0,0)  circle (.08);
\fill (1,0)  circle (.08);
\fill (2,0)  circle (.08);
\fill (3.5,0)  circle (.08);
\fill (4.5,0)  circle (.08);
\draw[very thick] (0,0)--(2.3,0);
\draw[very thick,dotted] (2.3,0)--(3.2,0);
\draw[very thick] (3.2,0)--(4.5,0);
\end{tikzpicture}
\end{center}
are not triangulated balls we can obtain by our constructions here.

Concerning two-dimensional $X(T)$,
only triangulations of $m$-gons where $m \leq 6$ are obtained this way.
From this perspective our setting seems somewhat limited. However this
is {not} the case: By \cite{FO} {\it every} line graph or triangulation of a polygon
may be obtained from $X(T)$'s by {\it joining} variables:
Divide $k[A]/J(T)$ by a regular sequence of variable differences.

In our article here, the regular sequences we divide out
by to get the Artin monomial ideals,
consists of variable differences $a^i_{j} - a^{i^\prime}_{j^\prime}$
{where} $i = i^\prime$ (the terms are in same color class $A_i$).
However we may also divide out by variable differences where $i \neq i^\prime$.
Starting form $k[A]/J$ we then get a ring $k[A^\prime]/J^\prime$, coming
from a surjection $A \pil A^\prime$. This is how we obtain any line graph
and any triangulation of a polygon in \cite{FO}.
We thus have the processes:
\[ \xymatrix{{}  & J(T) \ar^{\text{joining variables}}[rd] & \\
      & & J^\prime. \\
      \ijfr \ar^{\text{separating variables}}[ruu] & &}
    \]
    Each of the arrows above preserves the graded Betti numbers.
    Hence every $J^\prime$
    has the same graded Betti numbers as $\ijfr$.  
In Section \ref{sec:squeezed} we
obtained in this way the squeezed balls and spheres of G.Kalai \cite{Ka}.
\begin{problem} Given a subset $T \sus \bfA$. Which sets of
  variable differences $\{ a^i_{j}-a^{i^\prime}_{j^\prime} \}$ are such
  that:
  \begin{itemize}
  \item[i.] they are a regular sequence for $k[A]/J$, where
    $J = J(T)$ is the Alexander dual of the rainbow ideal $I(T)$?
  \item[ii.] in addition the resulting ideal $J^\prime$ is square free?
  \item[iii.] in addition the ideal $J^\prime$ defines a simplicial ball?
  \end{itemize}
  \end{problem}
  \cite{FO} achieves a full understanding of this for acyclic $T$
  constructed from trees (see Example \ref{ex:hypercube-tree}).
  Consequences in \cite{FlPart} are relations between partitions of vertices
  and facets of classes of simplicial
  complexes. General as the results of \cite{FlPart,FO} may seem, they
  consider only a tiny fraction of Artin monomial ideals:
  the second powers $(x_1, \ldots, x_n)^2$.
  In \cite{MPol} they do the above for standard polarizations $I^\st$
  of any monomial ideal $I$. They look at $J^\prime$ whose
  standard polarization $J^{\prime \st}$ is isomorphic to $I^\st$. They descibe
  how to give such $J^\prime$,  and they give 
  a poset structure on the family of $J^\prime$s.
  
\subsection{Triangulations and convex polytopes}
There are many combinatorial constructions of simplicial or simple polytopes.
(The dual of a simple polytope is a simplicial polytope.)
For instance the nestohedra of \cite{AA,Pos} are simple polytopes.

Our results here construct simplicial balls whose boundaries are
simplicial spheres.
These make  a much wider class as most of these simplicial spheres do
not have convex realizations \cite{BZPS}. Furthermore there is an
additional aspect: our constructions also triangulate the
interior of the sphere.

\begin{problem}
 Which simplicial balls $X(T)$ have boundaries $\pd X(T)$ corresponding
to the polytopes or duals of polytopes mentioned above \cite{AA,Pos},
or to other simple or simplicial polytopes?
\end{problem}

Two classes of such are Bier spheres of threshold simplicial complexes
and of simplicial complexes of ``short sets'' \cite{JTR}.

An $X(T)$ can only be a triangulated {\it polygon} of size $n$ for $n \leq 6$.
Triangulated polygons of sizes $n \geq 7$ are not of this type. (However
their ideals may be separated to give SR-ideals of $X(T)$'s, \cite{FO}.)
    For generalized permutahedra \cite{Pos}, one may also only
    construct $n$-gons for $n \leq 6$, but not for $n \geq 7$.
    Does this indicate some connection?

    \medskip
    There is a vast literature on triangulations of polyhedral complexes,
    see \cite{DLRS, LeSa}.
    In particular we mention regular triangulations \cite{BrGu, GBCP}.

    \begin{problem}
      For classes of triangulations described in the literature,
      which can be described as $X(T)$'s or regular quotients of these?
\end{problem}

\bibliographystyle{amsplain}
\bibliography{polarArtin-alco}

@article{PS,
  title={Liaison des vari{\'e}t{\'e}s alg{\'e}briques. I},
  author={Peskine, Christian and Szpiro, Lucien},
  journal={Inventiones mathematicae},
  volume={26},
  number={4},
  pages={271--302},
  year={1974},
  publisher={Springer}
}

@book{Eis,
  title={Commutative algebra: with a view toward algebraic geometry},
  author={Eisenbud, David},
  volume={150},
  year={2013},
  publisher={Springer Science \& Business Media}
}

@article {AFL,
    AUTHOR = {Almousa, Ayah and Fl{\o}ystad, Gunnar and Lohne, Henning},
     TITLE = {Polarizations of powers of graded maximal ideals},
   JOURNAL = {J. Pure Appl. Algebra},
  FJOURNAL = {Journal of Pure and Applied Algebra},
    VOLUME = {226},

YEAR = {2022},
    NUMBER = {5},
     PAGES = {Paper No. 106924, 33},
      ISSN = {0022-4049},
   MRCLASS = {13F20 (05E40 13F55 55U10)},
  MRNUMBER = {4328650},
       DOI = {10.1016/j.jpaa.2021.106924},
       URL = {https://doi.org/10.1016/j.jpaa.2021.106924},
}

@article{DFN,
  title={{Resolutions of co-letterplace ideals and generalizations of Bier spheres}},
  author={D’Al{\`\i}, Alessio and Fl{\o}ystad, Gunnar and Nematbakhsh, Amin},
  journal={Transactions of the American Mathematical Society},
  volume={371},
  number={12},
  pages={8733--8753},
  year={2019}
}

@book{HeHi,
  title={Monomial ideals},
  author={Herzog, J{\"u}rgen and Hibi, Takayuki},
  booktitle={Monomial Ideals},
  year={2011},
  publisher={Springer}
}

@book{Mi-St,
  title={Combinatorial commutative algebra},
  author={Miller, Ezra and Sturmfels, Bernd},
  volume={227},
  year={2005},
  publisher={Springer Science \& Business Media}
}

@article{FGH,
  title={Letterplace and co-letterplace ideals of posets},
  author={Fl{\o}ystad, Gunnar and Greve, Bj{\o}rn M{\o}ller and Herzog, J{\"u}rgen},
  journal={Journal of Pure and Applied Algebra},
  volume={221},
  number={5},
  pages={1218--1241},
  year={2017},
  publisher={Elsevier}
}

@article{FlPart,
  title={Partitions of Vertices and Facets in Trees and Stacked Simplicial Complexes},
  author={Fl{\o}ystad, Gunnar},
  journal={Graphs and Combinatorics},
  volume={40},
  number={4},
  pages={79},
  year={2024},
  publisher={Springer}
}

@article{Ka,
  title={Many triangulated spheres},
  author={Kalai, Gil},
  journal={Discrete \& Computational Geometry},
  volume={3},
  number={1},
  pages={1--14},
  year={1988},
  publisher={Springer}
}

@article{BZPS,
  title={Bier spheres and posets},
  author={Bj{\"o}rner, Anders and Paffenholz, Andreas and Sj{\"o}strand, Jonas and Ziegler, G{\"u}nter M},
  journal={Discrete \& Computational Geometry},
  volume={34},
  number={1},
  pages={71--86},
  year={2005},
  publisher={Springer}
}

@article{NR09,
  title={Betti numbers of monomial ideals and shifted skew shapes},
  author={Nagel, Uwe and Reiner, Victor},
  journal={{The Electronic Journal of Combinatorics}},
  volume={16},
  number={2},
  pages={3},
  year={2009}
}

@article{Mur11,
  title={Spheres arising from multicomplexes},
  author={Murai, Satoshi},
  journal={Journal of Combinatorial Theory, Series A},
  volume={118},
  number={8},
  pages={2167--2184},
  year={2011},
  publisher={Elsevier}
}

@article{ER98,
  title={{Resolutions of Stanley-Reisner rings and Alexander duality}},
  author={Eagon, John A and Reiner, Victor},
  journal={Journal of Pure and Applied Algebra},
  volume={130},
  number={3},
  pages={265--275},
  year={1998},
  publisher={Elsevier}
}

@incollection {Bj95,
    AUTHOR = {Bj\"{o}rner, A.},
     TITLE = {Topological methods},
 BOOKTITLE = {Handbook of combinatorics, {V}ol. 1, 2},
     PAGES = {1819--1872},
 PUBLISHER = {Elsevier Sci. B. V., Amsterdam},
      YEAR = {1995}
}

@book {BH93,
    AUTHOR = {Bruns, Winfried and Herzog, J\"{u}rgen},
     TITLE = {Cohen-{M}acaulay rings},
    SERIES = {Cambridge Studies in Advanced Mathematics},
    VOLUME = {39},
 PUBLISHER = {Cambridge University Press, Cambridge},
      YEAR = {1993},
     PAGES = {xii+403},
}

@article{AN,
  title={Linear strands of edge ideals of multipartite uniform clutters},
  author={Nematbakhsh, Amin},
  journal={Journal of Pure and Applied Algebra},
  volume={225},
  number={10},
  pages={106690},
  year={2021},
  publisher={Elsevier}
}

@article{St,
  title={{Balanced Cohen-Macaulay complexes}},
  author={Stanley, Richard P},
  journal={Transactions of the American Mathematical Society},
  volume={249},
  number={1},
  pages={139--157},
  year={1979}
}

@article{FO,
  title={{Triangulations of polygons and stacked simplicial complexes: separating their Stanley--Reisner ideals}},
  author={Fl{\o}ystad, Gunnar and Orlich, Milo},
  journal={Journal of Algebraic Combinatorics},
  volume={57},
  number={3},
  pages={659--686},
  year={2023},
  publisher={Springer}
}

@article{AA,
    AUTHOR = {Aguiar, Marcelo and Ardila, Federico},
     TITLE = {Hopf monoids and generalized permutahedra},
   JOURNAL = {Mem. Amer. Math. Soc.},
  FJOURNAL = {Memoirs of the American Mathematical Society},
    VOLUME = {289},
      YEAR = {2023},
    NUMBER = {1437},
     PAGES = {vi+119},
}

@article{Pos,
  title={Permutohedra, associahedra, and beyond},
  author={Postnikov, Alexander},
  journal={International Mathematics Research Notices},
  volume={2009},
  number={6},
  pages={1026--1106},
  year={2009},
  publisher={OUP}
}

@article{HU,
  title={The structure of linkage},
  author={Huneke, Craig and Ulrich, Bernd},
  journal={Annals of Mathematics},
  volume={126},
  number={2},
  pages={277--334},
  year={1987},
  publisher={JSTOR}
}

@book{Mu,
  title={Topology},
  author={Munkres, James R},
  volume={2},
  year={2000},
  publisher={Prentice Hall Upper Saddle River}
}

@article{Olt-Const,
  title={Constructible ideals},
  author={Olteanu, Anda},
  journal={Communications in Algebra},
  volume={37},
  number={5},
  pages={1656--1669},
  year={2009},
  publisher={Taylor \& Francis}
}

@article{Mil-ad2,
  title={{The Alexander duality functors and local duality with
  monomial support}},
  author={Miller, Ezra},
  journal={Journal of Algebra},
  volume={231},
  number={1},
  pages={180--234},
  year={2000},
  publisher={Elsevier}
}

@article{Bier,
author={Bier, Thomas},
title={{A remark on Alexander duality and the disjunct join}},
journal={preprint},
pages={8 pages},
year={1992}
}

@book{Mat,
  title={{Using the Borsuk-Ulam theorem: lectures on topological methods in combinatorics and geometry}},
  author={Matou{\v{s}}ek, Ji{\v{r}}{\'\i}},
  volume={2003},
  year={2003},
  publisher={Springer}
}

@article{MoKA,
  title={{On vertex decomposable simplicial complexes and their Alexander duals}},
  author={Moradi, Somayeh and Khosh-Ahang, Fahimeh},
  journal={Mathematica scandinavica},
  pages={43--56},
  year={2016},
  publisher={JSTOR}
}

@article{Rec-GS,
  title={Rectilinear Crossings in Complete Balanced d-Partite d-Uniform Hypergraphs},
  author={Gangopadhyay, Rahul and Shannigrahi, Saswata},
  journal={Graphs and Combinatorics},
  volume={36},
  number={3},
  pages={905--911},
  year={2020},
  publisher={Springer}
}

@inproceedings{Mat-DKPU,
  title={Effective heuristics for matchings in hypergraphs},
  author={Dufoss{\'e}, Fanny and Kaya, Kamer and Panagiotas, Ioannis and U{\c{c}}ar, Bora},
  booktitle={International Symposium on Experimental Algorithms},
  pages={248--264},
  year={2019},
  organization={Springer}
}

@book{GC,
  title={Graphs and cubes},
  author={Ovchinnikov, Sergei},
  year={2011},
  publisher={Springer Science \& Business Media}
}

@article{CN,
  title={{Monomial and toric ideals associated to Ferrers graphs}},
  author={Corso, Alberto and Nagel, Uwe},
  journal={Transactions of the American Mathematical Society},
  volume={361},
  number={3},
  pages={1371--1395},
  year={2009}
}

@book{Hat,
  title={Algebraic topology},
  author={Hatcher, Allen},
  year={2005},
  publisher={Cambridge University Press}
}

@article{BerComb,
  title={{Homological and combinatorial aspects of virtually Cohen--Macaulay sheaves}},
  author={Berkesch, Christine and Klein, Patricia and Loper, Michael C and Yang, Jay},
  journal={Transactions of the London Mathematical Society},
  volume={8},
  number={1},
  pages={413--434},
  year={2021},
  publisher={Wiley Online Library}
}

@article{BESVirt,
  title={Virtual resolutions for a product of projective spaces},
  author={Berkesch, Christine and Erman, Daniel and Smith, Gregory G},
  journal={Algebraic Geometry},
  volume={7},
  number={4},
  pages={460--481},
  year={2020},
  publisher={European Mathematical Society Publishing House}
}

@article{AlBa,
  title={{Symbolic powers of multipartite hypergraphs and
  Waldschmidt constant}},
  author={Alilooee, Ali and Banerjee, Arindam},
  journal={arXiv preprint arXiv:2103.06468},
  year={2021}
}

@phdthesis{Beck,
  title={Matchings and flows in hypergraphs},
  author={Beckenbach, Isabel Leonie},
  journal={Ph.D.thesis, FU-Berlin},
  year={2019}
}

@article{BeSc,
  title={Perfect f-matchings and f-factors in hypergraphs—A combinatorial approach},
  author={Beckenbach, Isabel and Scheidweiler, Robert},
  journal={Discrete Mathematics},
  volume={340},
  number={10},
  pages={2499--2506},
  year={2017},
  publisher={Elsevier}
}

@article{HLSR,
  title={Polarizations of equigenerated strongly stable ideals},
  author={Hao, Hanson and Luo, Yuyuan and Saint Rain, Sterling},
  booktitle={Joint Mathematics Meetings (JMM 2023)},
  organization={AMS},
  year={2023}
}

@book{Mig,
  title={Introduction to liaison theory and deficiency modules},
  author={Migliore, Juan C},
  volume={165},
  year={1998},
  publisher={Springer Science \& Business Media}
}

@article{GaHiPe,
  title={Resolutions of a-stable ideals},
  author={Gasharov, Vesselin and Hibi, Takayuki and Peeva, Irena},
  journal={Journal of Algebra},
  volume={254},
  number={2},
  pages={375--394},
  year={2002},
  publisher={Elsevier}
}

@article{AB,
  title={Polarizations and hook partitions},
  author={Almousa, Ayah and VandeBogert, Keller},
  journal={Journal of Pure and Applied Algebra},
  volume={226},
  number={9},
  pages={107056},
  year={2022},
  publisher={Elsevier}
}

@book{DLRS,
  title={Triangulations: structures for algorithms and applications},
  author={De Loera, Jes{\'u}s and Rambau, J{\"o}rg and Santos, Francisco},
  volume={25},
  year={2010},
  publisher={Springer Science \& Business Media}
}

@book{GBCP,
  title={Gröbner bases and convex polytopes},
  author={Sturmfels, Bernd},
  volume={8},
  year={1996},
  publisher={American Mathematical Soc.}
}

@article{EHM,
  title={{Monomial ideals and toric rings of Hibi type arising from a finite poset}},
  author={Ene, Viviana and Herzog, J{\"u}rgen and Mohammadi, Fatemeh},
  journal={European Journal of Combinatorics},
  volume={32},
  number={3},
  pages={404--421},
  year={2011},
  publisher={Elsevier}
}

@incollection{LeSa,
  title={Subdivisions and triangulations of polytopes},
  author={Lee, Carl W and Santos, Francisco},
  booktitle={Handbook of discrete and computational geometry},
  pages={415--447},
  year={2017},
  publisher={Chapman and Hall/CRC}
}

@book{BrGu,
  title={{Polytopes, rings, and K-theory}},
  author={Bruns, Winfried and Gubeladze, Joseph},
  year={2009},
  publisher={Springer Science \& Business Media}
}

@article{Mil-ad,
  title={Alexander duality for monomial ideals and their resolutions},
  author={Miller, Ezra},
  journal={arXiv preprint math/9812095},
  year={1998}
}

@book{Hall,
  title={Counterexamples in discrete geometry},
  journal={PhD-thesis},
  author={Hall, Huntington Tracy},
  year={2004},
  publisher={PhD-thesis, University of California, Berkeley}
}

@article{HT-cones,
  title={Resolutions by mapping cones},
  author={Herzog, J{\"u}rgen and Takayama, Yukihide},
  journal={Homology, Homotopy and Applications},
  volume={4},
  number={2},
  pages={277--294},
  year={2002},
  publisher={International Press of Boston}
}

@article{JTR,
  title={{Polytopal Bier spheres and Kantorovich--Rubinstein polytopes of weighted cycles}},
  author={Jevti{\'c}, Filip D and Timotijevi{\'c}, Marinko and {\v{Z}}ivaljevi{\'c}, Rade T},
  journal={Discrete \& Computational Geometry},
  volume={65},
  pages={1275--1286},
  year={2021},
  publisher={Springer}
}

@article{MPol,
  title={Polarization and depolarization of monomial ideals with application to multi-state system reliability},
  author={Mohammadi, Fatemeh and Pascual-Ortigosa, Patricia and S{\'a}enz-de-Cabez{\'o}n, Eduardo and Wynn, Henry P},
  journal={Journal of Algebraic Combinatorics},
  volume={51},
  number={4},
  pages={617--639},
  year={2020},
  publisher={Springer}
}

@article{LW,
  title={The resolutions of generalized co-letterplace ideals and their powers},
  author={Lu, Dancheng and Wang, Zexin},
  journal={Journal of Algebra},
  volume={673},
  pages={321--350},
  year={2025}
}

@article{Ya99,
  title={{Alexander duality for Stanley--Reisner rings and squarefree
  Nn-graded modules}},
  author={Yanagawa, Kohji},
  journal={Journal of Algebra},
  volume={225},
  number={2},
  pages={630--645},
  year={2000},
  publisher={Elsevier}
}

\end{document}